\documentclass[reqno,onefignum,onetabnum]{siamart171218}
\usepackage{subdepth}
\usepackage{overpic}
\usepackage{amsmath,amstext,amsbsy,amssymb,mathdots}
\usepackage{booktabs,bm}
\usepackage{mathrsfs}
\usepackage{tikz,cite}
\usetikzlibrary{intersections}

\usepackage{courier}
\usetikzlibrary{positioning}
\usetikzlibrary{shapes,arrows}
\usetikzlibrary{decorations.markings}
\usetikzlibrary{arrows.meta}
\usetikzlibrary{decorations.pathreplacing}
\usepackage[fleqn,tbtags]{mathtools}
\usepackage{amsfonts}
\usepackage{relsize}
\usepackage{xcolor,colortbl}
\usepackage{psfrag}
\usepackage{alltt}
\usepackage{arydshln}
\usepackage{enumerate}
\usepackage{epstopdf}
\usepackage{enumitem}
\usepackage{multirow}
\usepackage{framed}
\usepackage{cancel}
\usepackage{bbm}
\usepackage{relsize}
\usepackage{tikz}
\usepackage{tkz-euclide}
\usepackage{graphicx}
\usepackage{subcaption}
\usepackage[font=small]{caption}

\usepackage{algorithm}
\usepackage[noend]{algpseudocode}

\usepackage{hyperref}

\newcommand*{\CopyCounter}[2]{%
  \expandafter\def\csname c@#2\endcsname{\csname c@#1\endcsname}%
  \expandafter\def\csname p@#2\endcsname{\csname p@#1\endcsname}%
  \expandafter\def\csname the#2\endcsname{\csname the#1\endcsname}}

\CopyCounter{Theorem}{ProposedProblem}
\CopyCounter{Theorem}{Proposition}
\CopyCounter{Theorem}{Property}
\CopyCounter{Theorem}{Claim}
\CopyCounter{Theorem}{Lemma}
\CopyCounter{Theorem}{Corollary}
\CopyCounter{Theorem}{Conjecture}
\CopyCounter{Theorem}{Definition}
\CopyCounter{Theorem}{Example}
\CopyCounter{Theorem}{Remark}
\CopyCounter{Theorem}{Question}
\CopyCounter{Theorem}{Condition}
\CopyCounter{Theorem}{Criterion}
\CopyCounter{Theorem}{Observation}
\theoremstyle{plain}

\newcommand{\domain}{\Omega}

\newcommand{\boundary}{\partial \domain}

\newcommand{\x}{\bm{x}}

\newcommand{\ba}{\bm{a}}
\newcommand{\balpha}{\bm{\alpha}}

\newcommand{\bb}{\bm{b}}

\newcommand{\Bf}{\bm{f}}

\newcommand{\bk}{\bm{k}}

\newcommand{\w}{\bm{w}}
\newcommand{\y}{\bm{y}}

\newcommand{\J}{{\cal{J}}}

\newcommand{\R}{\mathbb{R}}
\newcommand{\E}{\mathbb{E}}

\newcommand{\M}{\mathcal{M}}
\newcommand{\A}{\mathcal{A}}
\newcommand{\Aext}{\A_e}
\newcommand{\I}{\mathcal{I}}

\newcommand{\mS}{\mathcal{S}}
\newcommand{\sbar}{\bar{s}}
\newcommand{\uhat}{\hat{u}}

\newcommand{\what}{\hat{w}}

\newcommand{\Vhat}{\widehat{V}}

\newcommand{\What}{\widehat{W}}
\newcommand{\Winterp}{W}
\newcommand{\PP}{\mathbb{P}}

\DeclareMathOperator*{\argmin}{arg\,min}
\DeclareMathOperator*{\argmax}{arg\,max}

\usepackage{xargs}
\usepackage{comment}
\usepackage[colorinlistoftodos,prependcaption,textsize=tiny]{todonotes}
\newcommandx{\avtodo}[2][1=]{\todo[linecolor=red,backgroundcolor=red!25,bordercolor=red,#1]{#2}}
\newcommandx{\mgtodo}[2][1=]{\todo[linecolor=blue,backgroundcolor=blue!25,bordercolor=blue,#1]{#2}}

\colorlet{lightgray}{gray!40}

\newcolumntype{"}{@{\hskip\tabcolsep\vrule width 1pt\hskip\tabcolsep}}
\setlist[description]{font=\normalfont\space}

\usepackage{lipsum}
\ifpdf
  \DeclareGraphicsExtensions{.eps,.pdf,.png,.jpg}
\else
  \DeclareGraphicsExtensions{.eps}
\fi

\crefname{subsection}{section}{sections}

\newsiamremark{remark}{\bf Remark}
\newsiamremark{hypothesis}{Hypothesis}
\crefname{hypothesis}{Hypothesis}{Hypotheses}

\headers{Uncertainty in Piecewise-Deterministic Markov Processes}{Cartee, Farah, Nellis, Van Hook, Vladimirsky}

\title{Quantifying and managing uncertainty in piecewise-deterministic Markov processes\thanks{Submitted to the editors August 2nd, 2020.
\funding{This work was started in an REU program partially supported by the NSF-RTG award DMS-1645643.
The first and last authors' work was supported in part by the NSF award DMS-1738010.  
The last author's work was also supported by the Simons Foundation Fellowship and by the NSF award DMS-2111522.}}}
\author{Elliot Cartee\thanks{Department of Mathematics, University of Chicago, Chicago, IL 60637
  (\email{ecartee@uchicago.edu}).} 
  \and
  Antonio Farah\thanks{Department of Mathematics, University of Texas at Austin, Austin, TX 78712
  (\email{antoniofarah@utexas.edu}).}
  \and 
  April Nellis\thanks{Department of Mathematics, University of Michigan, Ann Arbor, MI, 48109
  (\email{nellisa@umich.edu}).}
  \and 
  Jacob Van Hook\thanks{Department of Mathematics, University of Pennsylvania, Philadelphia, PA 19104
  (\email{jvanhook@sas.upenn.edu}).} 
  \and 
  Alexander Vladimirsky\thanks{Department of Mathematics and Center for Applied Mathematics, Cornell University, Ithaca, NY 14853
  (\email{vladimirsky@cornell.edu}).}
}

\usepackage{amsopn}

\begin{document}
\maketitle
\begin{abstract}
In piecewise-deterministic Markov processes (PDMPs) the state of a finite-dimensional system evolves continuously,
but the evolutive equation may change randomly as a result of discrete switches.  A running cost is integrated along the corresponding piecewise-deterministic trajectory up to the termination to produce the {\em cumulative cost} of the process.
We address three natural questions related to uncertainty in cumulative cost of PDMP models: (1) how to compute the Cumulative Distribution Function (CDF) of the cumulative cost when the switching rates are fully known; (2) how to accurately bound the CDF when the switching rates are uncertain; and (3) assuming the PDMP is controlled, how to select a control to optimize that CDF.
In all three cases, our approach requires posing a 
system of suitable hyperbolic partial differential equations, which are then solved numerically on an augmented state space. 
We illustrate our method using simple examples of trajectory planning under uncertainty for several 1D and 2D first-exit time problems.
In the Appendix, we also apply this method to a model of fish harvesting in an environment with random switches in carrying capacity. 
\end{abstract}

\begin{keywords}
piecewise-deterministic process, optimal control, Hamilton-Jacobi PDEs, uncertainty quantification, robustness 
\end{keywords}

\begin{AMS}
60J27, 49L20, 35F61, 60K40, 65M06, 65C40
\end{AMS}


\section{Introduction}
\label{sec:Intro}

Piecewise-deterministic Markov processes (PDMPs) provide a powerful formalism for modeling 
discrete random changes in a global environment.  That formalism is particularly useful when the number of {\em deterministic modes} of the global environment is relatively small and there is a high fidelity statistical characterization of mode-to-mode {\em switching} rates. 
Such processes arise in a broad range of applications, especially in the biological sciences \cite{RudnickiBioModels}. For example, they can be used to model keratin network formation \cite{BeilKeratin}, SIRS epidemic spread \cite{LiSIRSEpidemic}, genetic networks \cite{ZeiserGeneticNetworks},
and predator-prey systems \cite{benaim2016lotka, costa2016piecewise, HeningStrickler2019}. 
In other disciplines, applications of PDMPs include models of fatigue crack growth \cite{ChiquetFatigueCrack}, financial contagion \cite{DassiosContagion},
manufacturing processes \cite{akella1986optimal, bielecki1988optimality, olsder1980time, sethi2012hierarchical},
sustainable development, economic growth \& climate change  \cite{haurie2005multigenerational, haurie2006stochastic},
and path-planning under uncertainty \cite{andrews2013deterministic, ShenVlad, CarteeVlad_Poaching, GeeVladimirsky2020}.  

In this paper we focus on a computational framework for quantifying uncertainty in outcomes of PDMPs due to random switching times and possible uncertainty in switching rates.  If a PDMP system is controlled in real time, we also show that this uncertainty of outcomes can be actively {\em managed}. 

In our PDMP models, the full state of the system is described by a continuous component $\x \in \domain \subset \R^d$ and a 
discrete component $i \in \M = \{1, \ldots, M \}$ that represents the current deterministic ``mode''.
Starting from the initial configuration $(\x,i)$, the evolution of continuous component $\y(t)$ is defined by a (mode-dependent) ODE
\begin{align}
\label{eq:dynamics}
\y'(t)  &=  \Bf(\y(t),m(t)) \, = \, \Bf_{m(t)}(\y(t)),\\
\nonumber
\y(0) &= \x \in \domain,\\
\nonumber
m(0) &= i \in \M,
\end{align}
while the switches in mode $m(t)$ are based on a continuous-time Markov process on $\M$.  Using $\lambda_{ij}$ to denote the rate of $(i \to j)$ switching, we can write
\begin{align}
\label{eq:switching_rate}
\lim_{\tau\to 0}\frac{\PP(m(t+\tau)=j \, | \, m(t) = i)}{\tau}  & = \lambda_{ij}, 
\qquad \forall t \geq 0, \, i \in \M, \, j \in \M \backslash \{i\}.
\end{align}
Here, we focus on {\em exit-time} problems, in which the process stops as soon as the system reaches a compact exit set $Q \subset \domain$.  
Due to the random mode-switches, the exit-time 
$T_{\x, i} = \min \{t \geq 0 \, | \, \y(t) \in Q\}$ is also random, which makes it somewhat harder to approximate the distribution for our main object of study -- the {\em cumulative cost} of the PDMP $\J(\x,i)$. 

In addition to mode-dependent dynamics $\Bf: \domain \times \M \rightarrow \R^d$, we 
also include a mode-dependent
{\em running cost} $C:\domain \times \M \rightarrow (0,+\infty)$
and {\em exit cost} $q: Q \times \M \rightarrow [0,+\infty)$.
To simplify the notation, 
we will also sometimes use the mode as a subscript:
\begin{equation*}
C_i(\x) = C(\x,i), \quad
\Bf_i(\x) = \Bf(\x,i), \quad
q_i(\x) = q(\x,i), \quad
\text{etc.}
\end{equation*}
We will 
assume that $q_i$'s are continuous in $\x$, while $C_i$'s and $\Bf_i$'s are bounded and piecewise Lipschitz continuous.
The cumulative cost is then formally defined as
\begin{align}\label{eq:J}
\J_i(\x) &= \J(\x,i)    
 \; = \; \int_{0}^{T_{\x, i}} C \Big(\y(t), m(t) \Big) \, dt \, + \, 
q \Big( \y\left( T_{\x, i} \right)\!, m\left(T_{\x, i} \right)\! \Big).
\end{align}
We will generally assume that $\domain$ is a closed set and the process can continue on $\boundary \backslash Q,$
but if the dynamics forces us to leave $\domain$ before reaching $Q$, this will result in $\J = +\infty.$
We note that the notion of cumulative cost is much more common in {\em controlled} PDMPs,
where it is used to select criteria for control optimization. 
But we also consider $\J$ in this simpler uncontrolled case to focus on a single measurable outcome of the process. 

We develop our approach in this general setting, but our numerical experiments highlight that studying $\J$ is far from trivial even if 
$C \equiv 1,\, q \equiv 0,$ and $Q = \boundary$, yielding $\J(\x,i) =  T_{\x, i},$ the time until we reach the boundary. 
For a motivating example, consider a ``sailboat'' traveling with unit speed on an interval $\domain = [0,1]$ and subject to random mode (wind direction) switches.  
We will assume that it is moving rightward in mode 1 and leftward in mode 2,
the time intervals between mode switches are independent exponentially distributed random variables with rate $\lambda,$ 
and the process terminates as soon as the boat reaches $Q =\{0,1\}$.
While we describe this example in terms of sailboat navigation, similar ``velocity jump processes'' are also often used to model dispersal in biological systems \cite{othmer1988models, HillenDiffLimit}.  
But in contrast to our approach, the main focus there is on equations describing the evolving density of dispersing cells or organisms
rather than on the distribution of some performance measure $\J$ for individual organisms.
Another distinction is our assumption that each individual path terminates on reaching some exit set $Q$ -- this introduces additional structure, 
which we later leverage to obtain efficient numerical methods.

Throughout the paper, we take an 
exploratory approach, focusing on derivation of equations and numerical methods as well as instructive test problems
rather than proofs of convergence or realistic applications.
To streamline the presentation, we illustrate our methods on simple ``first-exit time'' problems\footnote{
To ensure computational reproducibility, our full code for all examples is available at\\ \url{https://github.com/eikonal-equation/UQ_PDMP}} 
in 1D and 2D similar to the sailboat example described above.
But in the Appendix we show how the same approach is useful more broadly (with general $C_i$'s and $\Bf_i$'s) 
by considering fish harvesting in an environment with random switches in carrying capacity. 

In~\cref{sec:compute}, we explain how the CDF for $\J$ can be computed by solving a system of coupled linear PDEs.
Our equations can be interpreted as a PDMP-adapted version of the Kolmogorov Backward Equation generalized to handle arbitrary running costs rather than just time.  
Another related approach is the previous development of numerical methods for the Liouville-Master Equation in \cite{annunziato2008analysis}. We also derive simpler recursive difference equations to compute the CDF for a discrete analog of our setting -- a random route-switching process on a graph.

In most real world applications, all switching rates $(\lambda_{ij})$ will be known only approximately and it is necessary to bound the results of this modeling uncertainty.  
In~\cref{sec:bounds}, we show how bounds on these switching rates can be used to bound the CDF of $\J.$
Interestingly, it turns out that 
it is easier to compute tight bounds if the switching rates are not assumed to be constant in time.

In many applications, the focus is on optimally {\em controlling} PDMPs (affecting the dynamics in each deterministic mode), with the notion of optimality typically based on the average-case outcomes (e.g., minimizing the expected total cost).  
Once a control is fixed, the same uncertainty quantification tools covered in sections \ref{sec:compute} and \ref{sec:bounds} become relevant.  
Moreover, the control can also be selected to {\em manage} the uncertainty, providing some robustness guarantees or minimizing the probability of undesirable outcomes.  
Following the latter idea, we introduce a method for optimizing the CDF of controlled PDMP models in~\cref{sec:optimize}.
We conclude by discussing further extensions and limitations of our approach in~\cref{sec:Conclude}.

\label{sec:intro}

\section{Computing the CDF}
\label{sec:compute}

Before discussing the methods for approximating the CDF for the randomly switching process described in \cref{sec:Intro}, we first consider the same challenge for  Markov-style switching on a graph in \cref{subsec:discrete}, turning to a continuous version in \cref{subsec:continuous}.
Numerical methods for the latter are then described in \cref{subsec:cdf_numerics} and illustrated by computational experiments in \cref{subsec:cdf_experiments}.

\subsection{Discrete PDMPs}
\label{subsec:discrete}

We start by reviewing a simple model of deterministic routing
on a directed graph with a finite node set $X = \{ \x_1, \ldots, \x_N \}$,
a set of directed edges $E \subset X \times X,$
and a target set $Q \subset X.$   
We will assume that $K: X \times X \to (0, +\infty]$ specifies the known cost of
possible ``steps'' (i.e., node-to-node transitions) with $K(\x,\x') = +\infty$ iff $(\x,\x') \not \in E$.
A {\em route} on this graph can be specified in {\em feedback form} by a mapping 
$F: X \to X$ such that $(\x,F(\x)) \in E \; \forall \x \in X.$
Given a starting position $\y_0 = \x \in X$, a path can be defined by a sequence
$\y_{n+1} = F(\y_n)$, terminating as soon as $\y_n \in Q$.  
We will further assume that the terminal cost charged at that point is specified by $q:Q \to [0, +\infty).$
If the path enters $Q$ after $\bar{n}(\x)$ steps,
its cumulative cost can be expressed as 
\begin{equation*}
\J(\x) \; = \; \sum\limits_{n=0}^{\bar{n}(\x) - 1} K \left(\y\strut_n, \y\strut_{n+1} \right) \, + \, q \left(\y\strut_{\bar{n}(\x)}\right),
\end{equation*}
with $\J(\x) = +\infty$ if the path remains forever in $X \setminus Q,$
which can happen if a route specified by $F$ contains loops.
The recursive relationship among $\J$ values makes it easy to recover all of them by solving a linear system
\begin{align} \label{eq:disc_determ_DP}
\nonumber
\J(\x) \; = \; &K \left(\x, F(\x) \right) \, + \, \J\left(F(\x)\right),  &&\qquad \forall \x \in X \setminus Q; \\
\J(\x) \; = \; &q(\x),   &&\qquad \forall \x \in Q. 
\end{align}

We will now consider a version of the problem with a total of $M$ different routes $F_1, \ldots, F_M$, each of them with its own
pair of running and terminal costs $(K_i, q_i)$ defined on the same graph. These routes are equivalent to the modes in a PDMP.
To simplify the notation, we will use $K_i(\x)$ as a shorthand for $K_i \left(\x, F_i(\x) \right).$
We define a random route-switching process by assuming that there
is a chance of switching to another route after each step.  That is, if the current route is $F_i$, the probability $p_{ij}$ of switching to $F_j$ after the next step 
 is known a priori for all $i,j \in \M= \{1, \ldots, M\}.$  The number of steps is now a random variable, along with the cost paid for all future steps.  In defining the new random cumulative cost $\J_i(\x)$, we note that the subscript only encodes the initial route used in the first step as we depart from $\x.$ 
It is easy to see that $u_i(\x) = \E\left[\J_i(\x)\right]$ should satisfy a recursive relation
\begin{align} \label{eq:disc_expect_DP}
\nonumber
u_i(\x) \; = \; &K_i(\x) \, + \, \sum\limits_{j=1}^M p_{ij} u_j\left(F_i(\x)\right),   &&\qquad \forall \x \in X \setminus Q, \, i \in \M;\\
u_i(\x) \; = \; &q_i(\x),   &&\qquad \forall \x \in Q, \, i \in \M. 
\end{align}
It is worth noting that this system of $MN$ linear equations lacks the nice causal property that we enjoyed in the deterministic (single route) case.  
There we knew that a finite $\J(\x)$ implied that the path from $\x$ prescribed by $F$ included no loops and reached $Q$ in a finite number of steps.
As a result, the part of system \cref{eq:disc_determ_DP} corresponding to such finite $\J$'s  was always triangular up to a permutation.
The same is clearly not true for the multi-route case of \cref{eq:disc_expect_DP}, where loops can easily arise as a result of random route-switching.

\begin{figure}
\begin{center}
\begin{tikzpicture}[x=0.68pt,y=0.68pt,yscale=1,xscale=1]
	\begin{scope}[decoration={markings, mark=at position 1 with {\arrow{Triangle}}}]
	
	\tikzset{endpoint/.style = {circle, fill=lightgray, draw=black, inner sep=0pt, minimum size=27pt}}
	\tikzset{midpoint/.style = {circle, draw=black, inner sep=0pt, minimum size=27pt}}
	\tikzset{arrow1/.style = {-{Latex[length=6pt,width=6pt]}}}
	\tikzset{arrow2/.style = {{Latex[length=6pt,width=6pt]}-{Latex[length=6pt,width=6pt]}}}
	
	\draw (0,200) node {Mode $1$};
	\draw (0,0) node {Mode $2$};
	\draw [decorate,decoration={brace,amplitude=10pt}] (250,30) -- (-250,30);
	\draw [decorate,decoration={brace,amplitude=10pt}] (-250,175) -- (250,175);
	
	\draw (-250,150) node[endpoint] (x11) {$\x\strut_1^1$};
	\draw (-150,150) node[midpoint] (x21) {$\x\strut_2^1$};
	\draw (-50,150)   node[midpoint] (x31) {$\x\strut_3^1$};
	\draw (50,150)    node[midpoint] (x41) {$\x\strut_4^1$};
	\draw (150,150)  node[midpoint] (x51) {$\x\strut_5^1$};
	\draw (250,150)  node[endpoint] (x61) {$\x\strut_6^1$};
	
	\draw (-250,50)   node[endpoint] (x12) {$\x\strut_1^2$};
	\draw (-150,50)   node[midpoint] (x22) {$\x\strut_2^2$};
	\draw (-50,50)     node[midpoint] (x32) {$\x\strut_3^2$};
	\draw (50,50)      node[midpoint] (x42) {$\x\strut_4^2$};
	\draw (150,50)    node[midpoint] (x52) {$\x\strut_5^2$};
	\draw (250,50)    node[endpoint] (x62) {$\x\strut_6^2$};
	
	\draw[arrow1] (x21.east) -- (x31.west) node[midway,anchor=south] {$p\strut_{11}$};
	\draw[arrow1] (x31.east) -- (x41.west) node[midway,anchor=south] {$p\strut_{11}$};
	\draw[arrow1] (x41.east) -- (x51.west) node[midway,anchor=south] {$p\strut_{11}$};
	\draw[arrow1] (x51.east) -- (x61.west) node[midway,anchor=south] {$p\strut_{11}$};
	
	\draw[arrow1] (x52.west) -- (x42.east) node[midway,anchor=north] {$p\strut_{22}$};
	\draw[arrow1] (x42.west) -- (x32.east) node[midway,anchor=north] {$p\strut_{22}$};
	\draw[arrow1] (x32.west) -- (x22.east) node[midway,anchor=north] {$p\strut_{22}$};
	\draw[arrow1] (x22.west) -- (x12.east) node[midway,anchor=north] {$p\strut_{22}$};
	
	\draw[arrow1] (x21.east) -- (x32.north) node[midway,anchor=west] {$p\strut_{12}$};
	\draw[arrow1] (x31.east) -- (x42.north) node[midway,anchor=west] {$p\strut_{12}$};
	\draw[arrow1] (x41.east) -- (x52.north) node[midway,anchor=west] {$p\strut_{12}$};
	\draw[arrow1] (x51.east) -- (x62.north) node[midway,anchor=west] {$p\strut_{12}$};
	
	\draw[arrow1] (x22.west) -- (x11.south) node[midway,anchor=east] {$p\strut_{21}$};
	\draw[arrow1] (x32.west) -- (x21.south) node[midway,anchor=east] {$p\strut_{21}$};
	\draw[arrow1] (x42.west) -- (x31.south) node[midway,anchor=east] {$p\strut_{21}$};
	\draw[arrow1] (x52.west) -- (x41.south) node[midway,anchor=east] {$p\strut_{21}$};

	\end{scope}
\end{tikzpicture}
\end{center}
\caption{
Fully discrete PDMP with $M=2$ modes and $N=6$ nodes.  In mode 1 the motion is always to the right; in mode 2 the motion is always to the left.  
The exit set is $Q= \{\x\strut_1, \x\strut_6\}.$}\label{fig:discrete}
\end{figure}
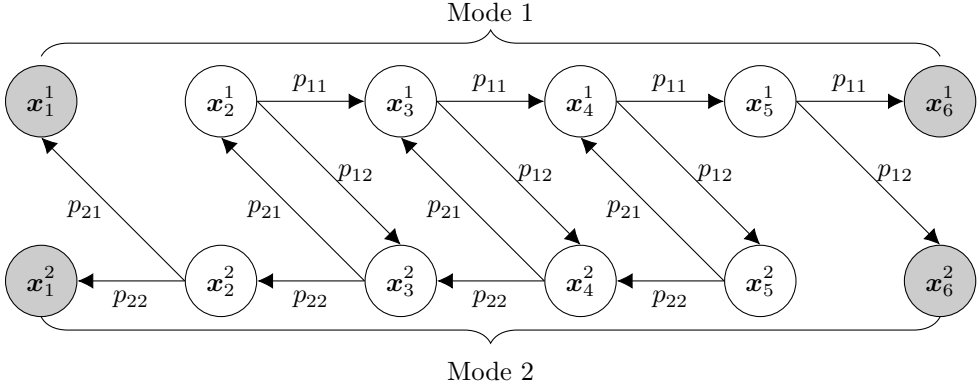


We note that this process can be also interpreted as a Markov chain on an extended graph. 
One would create $M$ copies of the original graph (on the nodes $\x^i_n$) with each route (or mode) $F_i$ represented as  a separate ``layer'' and inter-layer transitions governed by $p_{ij}$'s. \Cref{fig:discrete} illustrates one such example with two modes and associated probabilities $p_{11}, p_{12}, p_{21},$ and $p_{22}$. 
In the special case of $K_i \equiv 1$ and $q_i \equiv 0$ for all $i \in \M$, the above equations for $u_i$'s are simply describing the {\em mean hitting time} for the set $Q \times \M$.  However, we are interested in more general costs and would also like to compute the full CDFs $w_i(\x,s) = \PP(\J_i(\x) \leq s)$ for each $\J_i.$
It is easy to show that these functions must satisfy a recursive relationship
 \begin{equation}  \label{eq:discCDF}
w_i(\x, s) \; = \; 
   \sum\limits_{j=1}^{M} p_{ij} \, w_j\Big(F_i(\x), \, s -  K_i(\x) \Big),
   \qquad \qquad 
   \forall \x \notin Q, \, i \in \M, \, s > 0;
 \end{equation}
with the initial and boundary conditions
\begin{equation} \label{eq:discBC}
w_i(\x, s) = 
    \begin{cases}
    0, & \text{if } \left(\x \notin Q, \, s \leq 0 \right) \text{ or } \left( \x \in Q, \, s < q_i(\x) \right); \\
   1, & \text{if }\x \in Q, \, s \geq q_i(\x).
    \end{cases} 
 \end{equation}
We will assume that the range of $s$ values of interest is $\mS = [0,S]$, where $S$ is some constant specified in advance.

Based on the general properties of CDFs, all $w_i$'s are monotone non-decreasing and upper-semicontinuous
\footnote{
In addition, the finite size of $X$ guarantees that all $w_i$'s are piecewise-constant in $s$.  This can be used to construct a finite
time algorithm for solving \eqref{eq:discCDF} exactly despite the fact that $s$ is a continuous variable.  
We do not include this algorithm here due to space constraints and to keep the focus 
on the continuous-time setting, 
where discretizing $s$ and approximating $w_i$'s is generally unavoidable.}
 in $s$.
Moreover, the positivity of $K_i$'s ensures the explicit causality of this system:
in \cref{eq:discCDF} each $w_i(\x, s)$ can only depend on $w_i(\x', s')$ if $s'<s$.
Thus, the system can be solved in a single sweep (from the initial conditions at $s=0$, ``upward'' in $s$).

Still, it can be useful to precompute $s_i^0(\x) = \inf \{ s \mid w_i(\x,s) > 0 \}$
and $w_i^0(\x) = w_i \left(\x, s_i^0(\x) \right)$ by computations on $X$ alone. 
Intuitively, $s_i^0(\x)$ can be thought of as the minimum attainable cost starting in mode $i$ at position $\x$, and $w_i^0(\x)$ is the probability of attaining said cost.
It is easy to see that $s_i^0$ satisfies the recursive system:
\begin{align} \label{eq:disc_best_case_DP}
\nonumber
s_i^0(\x) \; = \; & K_i(\x) + \min\limits_{\mathop{j \in \M}\limits_{\scriptstyle\text{s.t. } p_{ij} > 0}} \left\{ s_j^0\left(F_i(\x)\right) \right\},  &&\qquad\forall \x \in X \setminus Q, i \in \M; \\
s_i^0(\x) \; = \; & q_i(\x),   &&\qquad \forall \x \in Q, i \in \M.
\end{align}
solvable by the standard Dijkstra's method in $O(MN \log (MN))$ operations.

The values of $w_i^0(\x)$ can also be found in the process of computing $s_i^0(\x)$. If $\I(\x) \subset \M$ is the $\argmin$ set in \cref{eq:disc_best_case_DP},
then
\begin{align}\label{eq:disc_best_case_probability}
\nonumber
w_i^0(\x) \; = \; & \sum\limits_{j \in \I(\x)} p_{ij} w_j^0 \left( F_i(\x) \right),  &&\qquad \forall \x \in X \setminus Q, \, i \in \M; \\
w_i^0(\x) \; = \; & 1,   &&\qquad \forall \x \in Q, \, i \in \M. 
\end{align}
Numerically solving \Cref{eq:disc_best_case_DP,,eq:disc_best_case_probability} can be advantageous because they are computed on the lower-dimensional domain $X\times \M$ instead of $X\times \M \times \mS$.
This information can then be used as initial/boundary conditions to solve \cref{eq:discCDF} on a smaller subset of $X \times \M \times \mS$.


\subsection{Continuous PDMPs}
\label{subsec:continuous}

We are now interested in extending our results from the discrete case to continuous settings.
The PDMP model described in \cref{sec:Intro} is based on continuous in time and space evolution of the state $\y(t)$
and continuous in time Markov chain governing the changes in mode $m(t).$   Here, we start with a somewhat simpler version, in which
this Markov chain is discretized in time, while the state evolution is continuous. Choosing some small fixed time interval $\tau > 0,$ 
we  assume that the system starting in mode $m(0) = i \in \M$ and state $\y(0) =\x \in \Omega \setminus Q$ evolves
according to an ODE 
$\y'(t)  =  \Bf_i(\y(t))$ with no random switches until the time
\begin{equation*}
\tau_{\x,i} \; = \; \min \left( \tau, \; \min \left\{ t \, | \, \y(t) \in Q \right\} \right),
\end{equation*}
at which point a switch to another mode may occur.   The process is repeated (starting from $\tilde{\x} = \y(\tau_{\x,i})$ and 
a possibly new mode $j,$ integrating the ODE over the time interval of length $\tau_{\tilde{\x},j},$ etc.) and the running cost is accumulated until
$\y(t)$ enters a compact exit set $Q.$

We define natural analogs for operators used to pose the graph routing problem in the previous subsection: 
\begin{align}
\label{eq:continuous_analog1} F_i(\x) \, &= \, \x + \int_0^{\tau_{\x,i}} \Bf_i(\y(t)) \, dt \, = \, \y\left(\tau_{\x,i}\right), \\
\label{eq:continuous_analog2} K_i(\x) \, &= \, 
\int_0^{\tau_{\x, i}} C_i(\y(t)) \, dt,
\end{align}
where $C_i:\domain \rightarrow (0,+\infty)$
is the running cost for that mode.  
We define the probability of switching to each mode $j$ at the end of time interval of length $\tau_{\x,i}$
by requiring consistency with the continuous in time Markov process described in \eqref{eq:switching_rate}.
In the latter, there could be multiple mode transitions over the time $\tau_{\x,i},$ and here we simply use 
the probability of {\em finishing} this time interval in mode $j$:  
 \begin{equation*}
 p_{ij}(\tau_{\x,i}) \; = \; \PP(m(\tau_{\x,i})=j \, | \, m(0)=i).
 \end{equation*} 
We compute these probabilities using a transition rate matrix $\Lambda = \left(\lambda_{ij}\right),$  where $\lambda_{ij}$'s encode the rate of $(i \rightarrow j)$ switching for $i \neq j$, while
the diagonal elements are defined by $\lambda_{ii} = - \sum_{j \neq i} \lambda_{ij}$.
The evolution of the probability matrix $P(t)=\left(p_{ij}(t)\right)$ is then given by an ODE
\begin{equation*}
\frac{d}{dt}P(t) = P(t)\Lambda, \qquad P(0)=I,
\end{equation*} and it follows that $P(\tau_{\x,i}) = \exp(\Lambda \tau_{\x,i})$.
Finally, if $\y(\tau_{\x,i}) \in Q$, we assume that the PDMP will immediately terminate with an exit cost of $q_j\left(\y(\tau_{\x,i})\right)$, where $j$ is the final mode 
after a possible last transition.

With this notation in hand, we can define the same functions characterizing the random cumulative cost: $u_i$,  $w_i$, $s_i^0$, and $w_i^0$
will all satisfy the same recursive formulas already defined on a graph in the previous subsection.  
The only caveat is that $p_{ij}$'s will need to be replaced by $p_{ij}(\tau_{\x,i})$. 
Since $\tau$ and $\tau_{\x,i}$ are equivalent except on a small neighborhood of $Q$, in the following sections we will slightly abuse the notation by referring to $\tau$ to simplify the formulas.

The original setting of  \cref{sec:Intro} (with continuous in time Markov chain for mode switching) can be obtained from the above in the limit by letting $\tau \rightarrow 0$.
A standard argument based on a Taylor series expansion 
 shows that the expected costs $u_i(\x) = \mathbb{E}[\J_i(\x)]$ formally satisfy a system of 
linear PDEs:
\begin{equation}\label{eq:expected_pde}
\nabla u_i(\x)\cdot\Bf_i(\x) \, + \, C_i(\x) \, + \, \sum_{j \not= i}\biggl[\lambda_{ij}\bigl(u_j(\x) - u_i(\x)\bigr)\biggr] \; = \; 0
\end{equation}
with boundary conditions 
$
u_i(\x) = q(\x,i)
$
on $Q \times \M.$
We omit the derivation of \eqref{eq:expected_pde} for the sake of brevity but  
use a similar approach below to derive a system of PDEs satisfied by the cumulative distribution functions $w_i(\x,s).$
The first order approximations of the transition probabilities are:
\begin{equation}\label{eq:probabilities}
\begin{aligned}
p_{ij}(\tau) &= 1 - e^{-\lambda_{ij}\tau} + o(\tau) = \lambda_{ij}\tau + o(\tau), &j \not= i \\
p_{ii}(\tau) &= 1 - \sum_{j\not= i}\lambda_{ij}\tau + o(\tau).
\end{aligned}
\end{equation}
The first-order approximation of the dynamics in \cref{eq:continuous_analog1} is
\begin{equation}\label{eq:approx_dynamics}
F_i(\x) = \x + \tau \Bf_i(\x) + o(\tau),
\end{equation}
and the first-order approximation of the running cost in \cref{eq:continuous_analog2} is
\begin{equation}\label{eq:approx_cost}
K_i(\x) = \int_0^\tau C_i(\y(t))dt = \tau C_i(\x) + o(\tau),
\end{equation}

Plugging in our approximations \eqref{eq:probabilities}, \eqref{eq:approx_dynamics}, and \eqref{eq:approx_cost} into the recursive relationship in Equation \eqref{eq:discCDF} and then Taylor expanding $w_i$ gives:
\begin{equation}\label{eq:approx_recursive}
w_i(\x,s) = \bigl(1 - \sum_{j\not=i}\lambda_{ij}\tau\bigr)w_i\bigl(F_i(\x),s - \tau C_i(\x)\bigr) + \sum_{j\not=i} \lambda_{ij}\tau w_j\bigl(F_i(\x),s - \tau C_i(\x)\bigr) + o(\tau)
\end{equation}
\begin{equation*}
w_i(\x,s) = w_i\biggl(F_i(\x), \, s - \tau C_i(\x)\biggr) + \tau\sum_{j\not=i} \lambda_{ij}\left[w_j(\x,s)-w_i(\x,s)\right] + o(\tau)
\end{equation*}
\begin{equation*}
0 = \tau\nabla w_i(\x,s)\cdot\Bf_i(\x) - \tau C_i(\x) \frac{\partial w}{\partial s}(\x,s) + \tau\sum_{j\not=i} \lambda_{ij}\left[w_j(\x,s)-w_i(\x,s)\right] + o(\tau),
\end{equation*}
where $\nabla = \left( \frac{\partial}{\partial x_1}, \frac{\partial}{\partial x_2}, ..., \frac{\partial}{\partial x_d}\right)$ denotes the gradient in the spatial coordinates.
Dividing both sides by $\tau$ and then taking the limit as $\tau \to 0$,
we obtain a linear PDE for each mode $i$:
\begin{equation}\label{eq:cdf_pde}
\nabla w_i(\x,s)\cdot\Bf_i(\x) - C_i(\x) \frac{\partial w_i}{\partial s}(\x,s) + \sum_{j\not=i} \lambda_{ij}\left[w_j(\x,s)-w_i(\x,s)\right] \; = \; 0.
\end{equation}
The above derivation is only formal since it assumes that $w_i$'s are sufficiently smooth.  In reality, they will be often non-differentiable and even discontinuous at isolated points;
nevertheless, these value functions can be still interpreted as weak (viscosity) solutions \cite{DavisFarid1999}, which can be approximated numerically by discretizing \cref{eq:approx_recursive}.\\
This system of 
PDEs satisfies the initial/boundary conditions:
\begin{equation}\label{eq:cdf_ic}
w_i(\x,0) =
\begin{cases}
1, &\forall \x \in Q \text{ s.t. } q(\x,i) = 0, \\
0, &\text{otherwise,}
\end{cases}
\end{equation}
\begin{equation}\label{eq:cdf_bc}
w_i(\x,s) = 
\begin{cases}
1, &\forall \x \in Q \text{ s.t. } q(\x,i) \le s, \\
0, &\forall \x \in Q \text{ s.t. } q(\x,i) > s.
\end{cases}
\end{equation}
The above conditions are sufficient when $Q = \boundary$ or if $\domain$ is invariant under all vector fields $\Bf_i.$  All of our examples considered in the next sections fall in this category.  But more generally, if vector fields are such that a trajectory might leave $\domain$ prior to reaching $Q$, one could treat this event as an immediate failure, essentially imposing 
$w_i(\x, s) = 0$ for all $\x \not \in \domain$ and all $s \in \R.$

As in the discrete case of \cref{subsec:discrete}, it can be useful to precompute the minimum attainable cost to use as initial/boundary conditions when solving \cref{eq:cdf_pde}.
From the discrete case we recall that $s^0_i(\x) = \inf\{s \; | \; w_i(\x,s) > 0\}$ denotes the minimal cost possible when starting from position $\x$ in mode $i$ assuming that transitions between modes can occur whenever desired.
In the continuous case these transitions can occur without delay, and therefore $s^0_i(\x) = s^0_j(\x)$ for all $i$ and $j$ in $\M$, so we will replace all of these with $s^0(\x)$.
(Also, unlike in the discrete case, it is entirely possible that $w_i(\x, s^0(\x)) = 0$ for all $i$.  The cost of $s^0(\x)$ might be attainable only through perfectly timed transitions, which in the continuous case would happen with probability zero.) 
A formal Taylor series expansion of \cref{eq:disc_best_case_DP} yields the following differential equation and boundary conditions for $s^0(\x)$:
\begin{align}\label{eq:cont_best_case}
\nonumber &\min_i \left\{C_i(\x) + \nabla s^0(\x)\cdot \Bf_i(\x)\right\} = 0, \qquad\qquad &\x \in \Omega \setminus Q; \\
&s^0(\x) = \min_i \left\{ q_i(\x) \right\}, \qquad\qquad &\x \in Q.
\end{align}

We are also interested in the probability $w^0_i(\x)$ of attaining that minimal cost $s_0(\x)$ when starting from mode $i$ and position $\x$.
If we denote the argmin set of \cref{eq:cont_best_case} as $\I(\x)$, then $w^0_i(\x)$ formally satisfies the following system:
\begin{align}\label{eq:cont_best_case_prob}
\nonumber &0 = \nabla w_i^0(\x)\cdot \Bf_i(\x) + \sum_{j \not= i} \lambda_{ij} \left[ w_j^0(\x) - w_i^0(\x)\right], &&\x \in \Omega \setminus Q, i \in \I(\x); \\
\nonumber &w_i^0(\x) = 1, &&\x \in Q, i \in \I(\x); \\
&w_i^0(\x) = 0, &&\x \in \Omega, i \not\in \I(\x).
\end{align}
Once $s^0(\x)$ and $w^0_i(\x)$'s are known, the computation of $w_i$'s can be restricted to $\big\{ (\x,s) \, \mid \, s \in (s_0(\x), S] \big\},$
solving PDEs \cref{eq:cdf_pde} with ``initial'' conditions $w_i(\x, s^0(\x)) = w_i^0(\x)$.

\begin{remark} [{\bf Related work on Liouville-Master Equation}]
\label{rem:Liouville}
An approach similar to the one presented in this section can be used to derive
PDEs for the time-dependent joint PDMP-state CDF on $\domain \times \M$.
The initial conditions to those PDEs would be based on 
a specific initial configuration $(\x_0, i_0)$ or, more generally, on a specific initial joint CDF on $\domain \times \M$.
This is precisely the setting in \cite{annunziato2008analysis}, where a finite-difference numerical method for 
the ``Liouville-Master Equation'' was developed and tested for the special case of $d=1.$
If one is willing to increase the dimension of the problem, this can be viewed as a more general approach 
than ours (since $\J$ can be viewed as just another component of the continuous state variable).
But the need to solve PDEs separately for different $(\x_0, i_0)$-specific initial conditions is a serious drawback.
Moreover, computing the time-dependent joint CDF seems more suitable for {\em finite-horizon} PDMPs
(where the process terminates after a pre-specified time $T$) 
rather than in our setting (where the process terminates as soon as it reaches $Q \subset \domain$).
\end{remark}

\subsection{Numerics for CDF computation}
\label{subsec:cdf_numerics}

We will approximate the domain $\Omega$ with a rectangular grid of points $\{ \x_{\bk} \}$ with grid spacing $\Delta x$, where $\bk=(k_1, \ldots, k_d)$ is a multi-index and $\x_{\bk}  = (k_1 \Delta x, \ldots, k_d \Delta x).$
We will also approximate the second argument of the CDF with regularly spaced points $s_n = n\Delta s$.

We will derive equations for a 
grid-function
$W_{i,\bk}^n \approx w_i(\x_{\bk}, s_n),$
with $W_{i,\bk}^0$ values determined by the initial conditions \cref{eq:cdf_ic}. 
To simplify the discussion, we assume that both $\boundary$ and $Q$ are grid-aligned,
with boundary values prescribed by 
\cref{eq:cdf_bc}.

Equation \cref{eq:discCDF} is then naturally interpreted as a recipe for a semi-Lagrangian discretization using a pseudo-timestep of length $\tau$. 
To obtain the first-order scheme, we can use the linear approximations (\ref{eq:approx_dynamics}-\ref{eq:approx_cost})
in 
formula \cref{eq:approx_recursive}, yielding
the following equation at each gridpoint $\x_{\bk} \in \domain$, mode $i \in \M$, and cost threshold $s_n$:
\begin{equation}\label{eq:cdf_update}
W_{i,\bk}^n = \sum_{j =1}^M p_{ij}(\tau) \Winterp_{j} \big(\x_{\bk} + \tau \Bf_i(\x_{\bk}), \, s_n - \tau C_i(\x_{\bk}) \big), 
\end{equation}
where $\Winterp_j: \Omega \times \R \rightarrow \R$ is 
the result of interpolating the grid-function $W_{j,\bk}^n$ in both $\x$ and $s$ variables, and the $p_{ij}$'s are defined as in \cref{eq:probabilities}.
In our implementation, 
all $\Winterp_j$'s are defined by multi-linear interpolation, but more sophisticated interpolation techniques (e.g., based on ENO/WENO \cite{shu1998essentially}) may be used instead
to decrease the numerical viscosity.  More accurate approximations of $F_i$ and $K_i$ could be also employed to increase the formal order of accuracy of the discretization.  
For fully deterministic processes, similar semi-Lagrangian schemes have been proven to converge under the grid refinement to 
a discontinuous viscosity solution
on all compact sets not containing the discontinuity \cite{BardiFalconeSoravia1999}.  While we do not attempt to prove this here, our numerical experiments indicate that the same holds true in piecewise-deterministic problems.

Our update formula \cref{eq:cdf_update} is only valid when $\x_{\bk} + \tau \Bf_i(\x_{\bk})$ remains in $\domain$.
With grid-aligned $\boundary,$ a rather conservative sufficient condition for this is 
\begin{equation}\label{eq:tau_ineq1}
\tau \cdot \max_i \bigl\{ \max_{\x} \bigl\{ \left| \Bf(\x,i) \right| \bigr\} \bigr\} \le \Delta x.
\end{equation}
Furthermore, we would like to ensure that our updates are \emph{causal}, that is the right hand side of \cref{eq:cdf_update} depends only upon the 
$W^{n'}$ values with $n' < n.$
While not strictly necessary, this ensures that the updates for each mode are uncoupled, speeding up the computation.
A sufficient condition for this is 
\begin{equation}\label{eq:tau_ineq2}
\tau \cdot \min_i \bigl\{ \min_{\x} \bigl\{ C(\x,i) \bigr\} \bigr\} \ge \Delta s.
\end{equation}

The inequality \cref{eq:tau_ineq1} is only needed if we want to use the same $\tau$ at all grid points instead of selecting a smaller time step near $\boundary$ only.
But if this $\tau$-uniformity is desired, satisfying both \cref{eq:tau_ineq1} and \cref{eq:tau_ineq2} requires
\begin{equation}\label{eq:cfl}
\frac{\Delta s}{\min\{C\}} \le \frac{\Delta x}{\max\{|\Bf|\}}.
\end{equation}
We note that, even though the above restriction looks similar to a 
Courant-Friedrichs-Lewy (CFL) condition,
it is not needed to guarantee the stability (semi-Lagrangian discretizations are unconditionally stable), but simply to ensure the causality (and hence the efficiency) of our discretization.  

Under certain conditions, \cref{eq:cdf_update} may be also re-interpreted as a finite differences discretization of the PDE \cref{eq:cdf_pde}.
To give a concrete example, suppose that $d=1$, and the domain $\Omega = [0,1]$ is approximated by a grid of regularly spaced points denoted $x_k = k\Delta x$. 
Furthermore, suppose that there is a mode $i$ where $C_i\equiv1$, and $f_i(x_k) = f_{ik} > 0$. 
If we choose $\tau = \Delta s$, then \cref{eq:cdf_update} for $n+1$ becomes:
\begin{align*}
W_{i,k}^{n+1} &= \sum_{j =1}^M p_{ij}(\Delta s) \Winterp_{j}(x_k + f_{ik}\Delta s, s_n) \\
&= \Winterp_i(x_k +f_{ik}\Delta s,s_n) + \sum_{j\not= i} \lambda_{ij}\Delta s \left[ \Winterp_j(x_k +f_{ik}\Delta s,s_n) - \Winterp_i(x_k +f_{ik}\Delta s,s_n)\right]\\
&= W_{i,k}^n + \frac{f_{ik}\Delta s}{\Delta x}\left(W_{i,k+1}^n-W_{i,k}^n\right) + \sum_{j\not= i} \lambda_{ij}\Delta s \left(\Winterp_j-\Winterp_i\right)(x_k + f_{ik}\Delta s,s_n);
\end{align*}
\begin{equation}\label{eq:eulerian}
f_{ik} \left[ \frac{W^n_{i,k+1} - W^n_{i,k}}{\Delta x}\right] - \left[ \frac{W_{i,k}^{n+1}-W_{i,k}^n}{\Delta s}\right]  + \sum_{j \not= i} \lambda_{ij}\left(\Winterp_j-\Winterp_i\right)(x_k +f_{ik}\Delta s,s_n) = 0,
\end{equation}
which is a consistent first-order finite differences discretization of \cref{eq:cdf_pde}.
Furthermore, in this 1D example, the CFL condition for this discretization is exactly \cref{eq:cfl}. The scheme \cref{eq:eulerian} is monotone (and thus stable \cite{crandall1980monotone}) whenever this CFL condition is satisfied.
It is important to note that the summands in \cref{eq:eulerian} are evaluated at $x_k + f_{ik}\Delta s$ (and therefore are convex combinations of $W^n$ values at 
$x_k$ and $x_{k+1}$).
Evaluating those terms at the naive choice of $x_k$ would result in a non-monotone discretization, which is in fact unstable.

To compute the minimum attainable cost $s^0(\x)$ and the probability $w^0(\x)$ of attaining it, we use first-order semi-Lagrangian discretizations of \cref{eq:cont_best_case} and \cref{eq:cont_best_case_prob}.
For $d=1$, 
the discretized equations for $s^0(x)$ are
\begin{align}\label{eq:disc_s0}
\nonumber s^0(x_k) &= \min_i \left\{ C_i(x_k)\frac{\Delta x}{|f_i(x_k)|} + s^0\left(x_{k'}\right)\right\}, &&x_k \not\in Q; \\
s^0(x_k) &= \min_i \left\{ q_i(x_k) \right\}, && x_k \in Q;
\end{align}
where
\begin{equation*}
k' = \begin{cases} k+1, &f_i(x_k) > 0; \\ k-1, &f_i(x_k) < 0. \end{cases}
\end{equation*}
In 1D, this system of equations can be solved efficiently with two iterative ``sweeps'' -- first increasing and then decreasing in $k.$
In higher space dimensions, it can be solved in $O(MN\log (N))$ time using a Dijkstra-like method.

In the process of solving for $s^0(x)$, we also solve \cref{eq:cont_best_case_prob} using a first-order semi-Lagrangian scheme.
Using $\I(x_k)$ to denote the argmin set of \eqref{eq:disc_s0}, the values of $w_i^0$ are initialized according to
\begin{equation*}
w_i^0(x_k) = \begin{cases} 1, &x_k \in Q, i \in \I(x_k); \\ 0, &\text{otherwise.} \end{cases}
\end{equation*}
Whenever the value of $s^0(x_k)$ is updated, we simultaneously update $w_i^0$ according to
\begin{equation}\label{eq:disc_w0}
w_i^0(x_k) = \begin{cases} w_i^0\left(x_{k'}\right) + \frac{\Delta x}{\left|f_i\left(x_{k}\right)\right|} \sum_{j \not= i} \lambda_{ij}\left[w_j^0\left(x_{k'}\right) - w_i^0\left(x_{k'}\right)\right], &i \in \I(x_k); \\ 0, &i \not\in \I(x_k). \end{cases}
\end{equation}
These values of $s^0(x)$ and $w_i^0(x)$ are then used as initial/boundary conditions\footnote{
Since the graph of $s^0(x)$ is generally not grid aligned in $\domain \times \mS$, such a domain restriction requires 
either a use of ``cut cells'' just above $s = s^0(x)$ or a conservative ``rounding up''  of $s^0$ values.  
Our implementation relies on the latter, which introduces additional $O(\Delta s)$ errors.}
 for computing $w_i(x)$.
This provides a speed improvement and also reduces the smearing of $w_i$'s discontinuities due to numerical viscosity.

\subsection{Experimental Results}
\label{subsec:cdf_experiments}



We illustrate our approach with three examples of uncontrolled PDMPs on $\R$ and $\R^2$.
In all of these, we assume $Q = \boundary,$ $C \equiv 1,$ and $q \equiv 0,$ ensuring that the cumulative cost $\J$ corresponds to the time to $\boundary$.
For simplicity, we will also assume uniform transition rates; i.e., $\lambda_{ij}=\lambda > 0$ for all $i \neq j.$

\begin{figure}
	$\begin{array}{cccc}
	\includegraphics[width=0.22\textwidth]{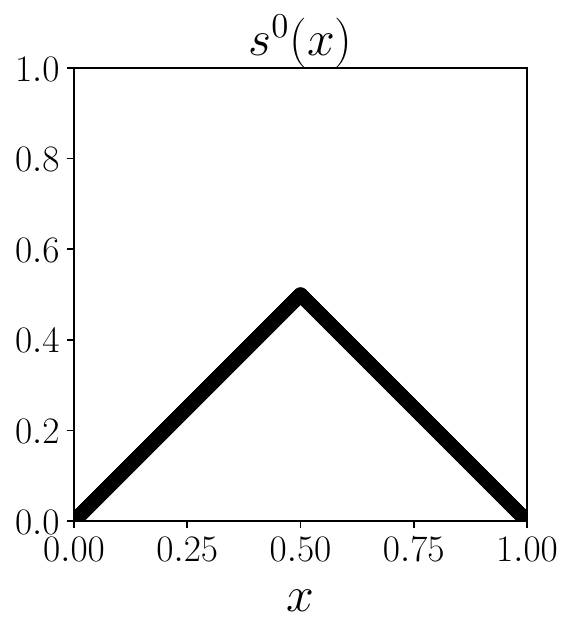} & \includegraphics[width=0.22\textwidth]{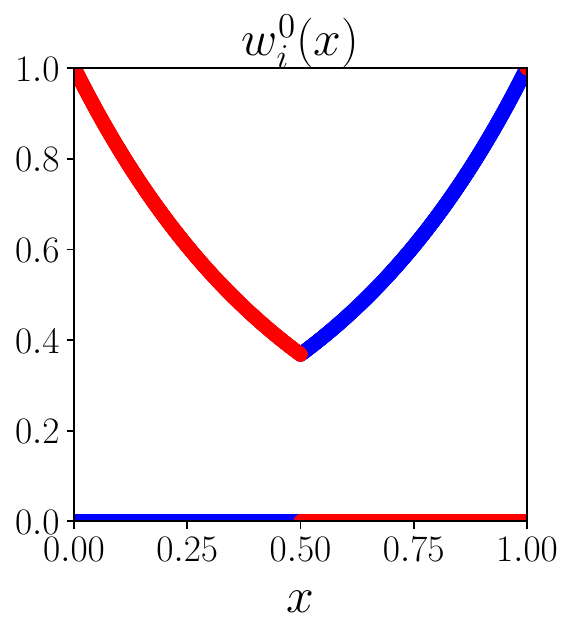} & \includegraphics[width=0.22\textwidth]{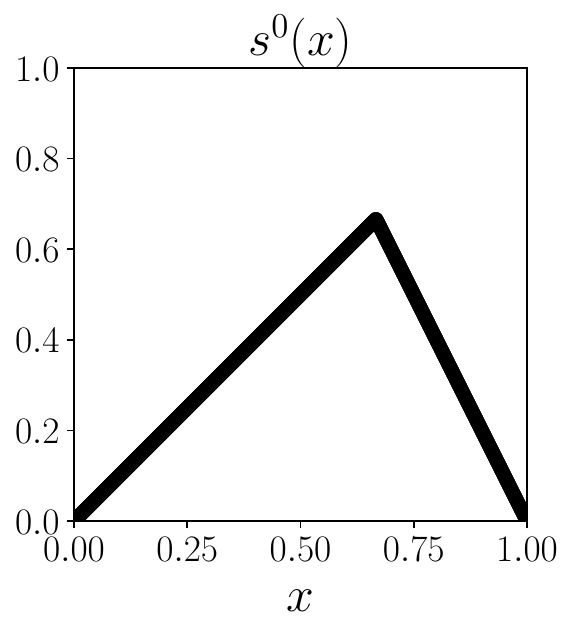} & \includegraphics[width=0.22\textwidth]{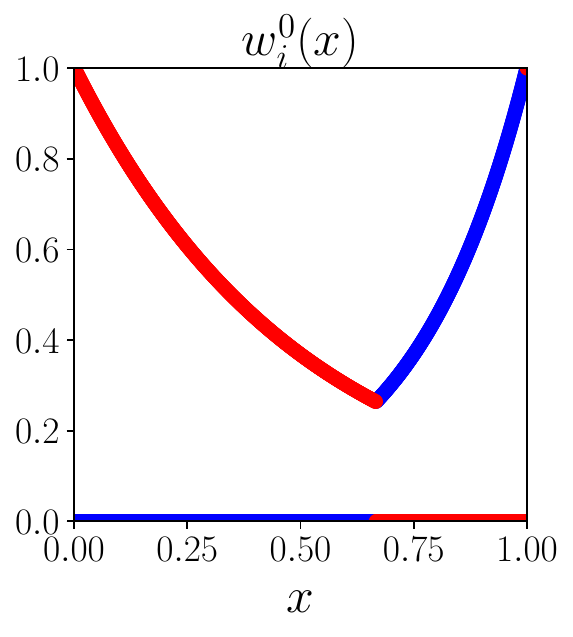} \\
	(A) & (B) & (C) & (D)
	\end{array}$
	\caption{Minimum cost $s^0(x)$ and probability $w^0_i(x)$ of attaining that minimum cost. Subfigures (A) and (B) are for Example 1 and (C) and (D) are for Example 2.
	Graphs of $w^0_i$ are shown in blue for $i=1$ and in red for $i=2.$}
	\label{fig:best_case}
\end{figure}

{\bf Example 1: }
We start by considering a ``sailboat'' test problem described in the introduction with $\Omega = [0,1], \, Q = \{0,1\}, \, M = 2, \, f_i(x) = (-1)^{i+1},$
and  symmetric transition rates $\lambda_{12}=\lambda_{21} = 2.$
For a fixed number $N$ of gridpoints, we set $\Delta s = \Delta x = \frac{1}{N-1}$, as this is the largest value of $\Delta s$ that satisfies \cref{eq:cfl}.
Moreover, this guarantees that no actual interpolation is necessary in \cref{eq:cdf_update}, as $\Winterp_j$ is only evaluated at gridpoints.
We note that solving these discretized equations is equivalent to finding the CDF of a discrete PDMP such as the one pictured in \cref{fig:discrete}, except with a larger 
number $N$ of nodes.  We solve this problem for $s \in [0,1]$, but also precompute $s^0(x)$ and $w^0_i(x)$  (see \cref{fig:best_case}(A-B)) to reduce the computational domain for $w_i$'s.

\begin{figure}
	$\begin{array}{cccc}
	\includegraphics[width=0.22\textwidth]{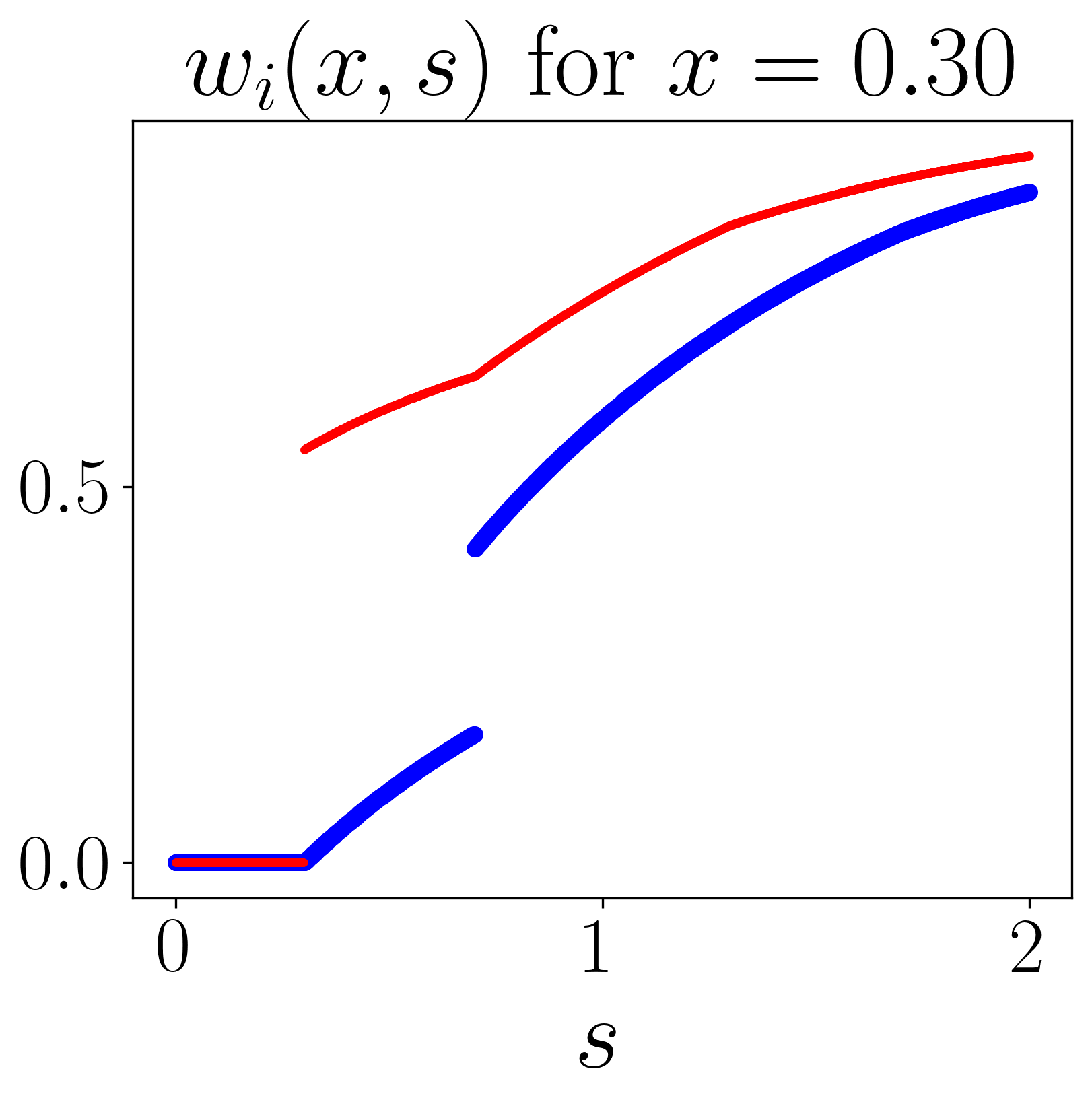} & \includegraphics[width=0.22\textwidth]{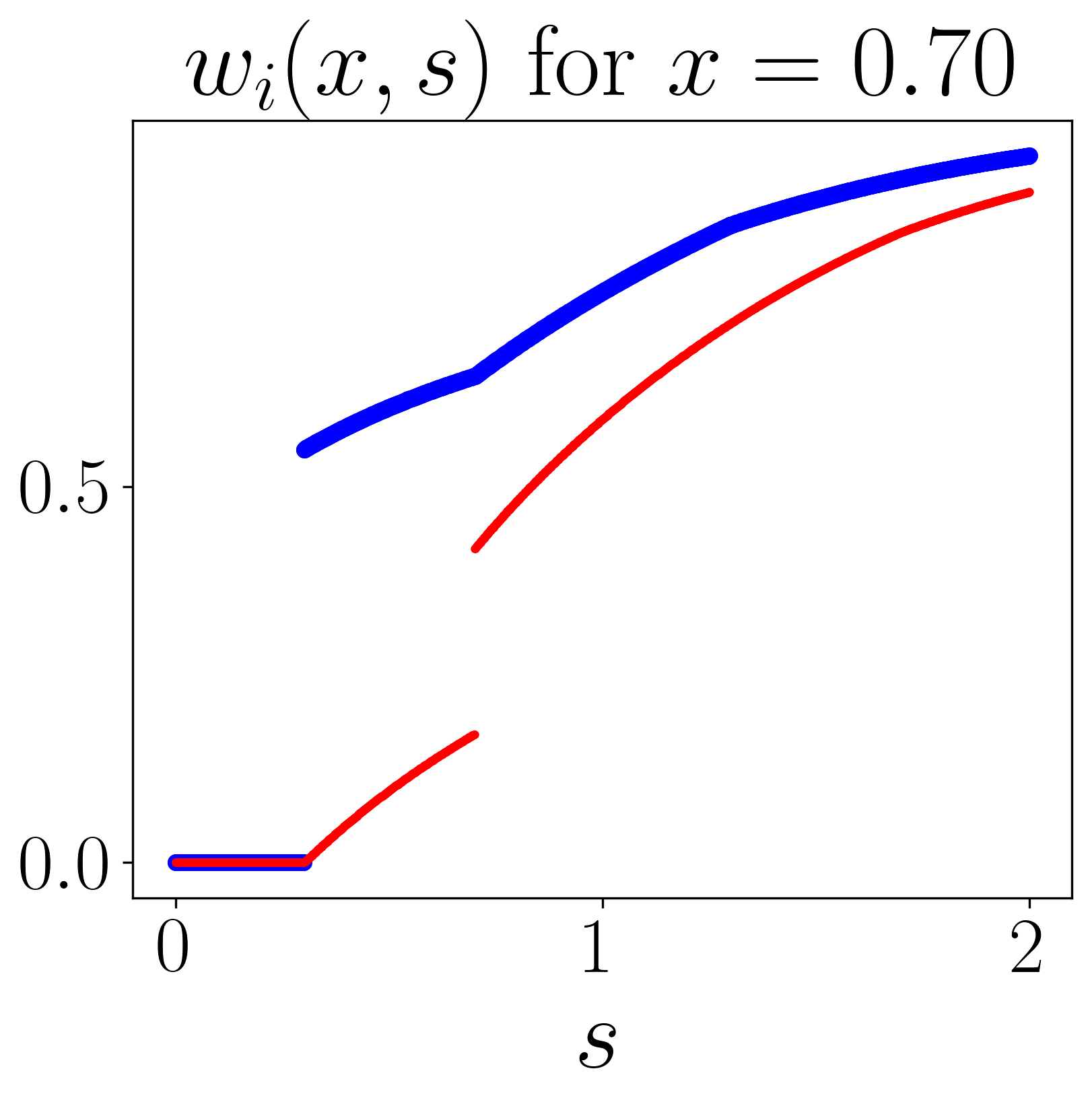} & \includegraphics[width=0.22\textwidth]{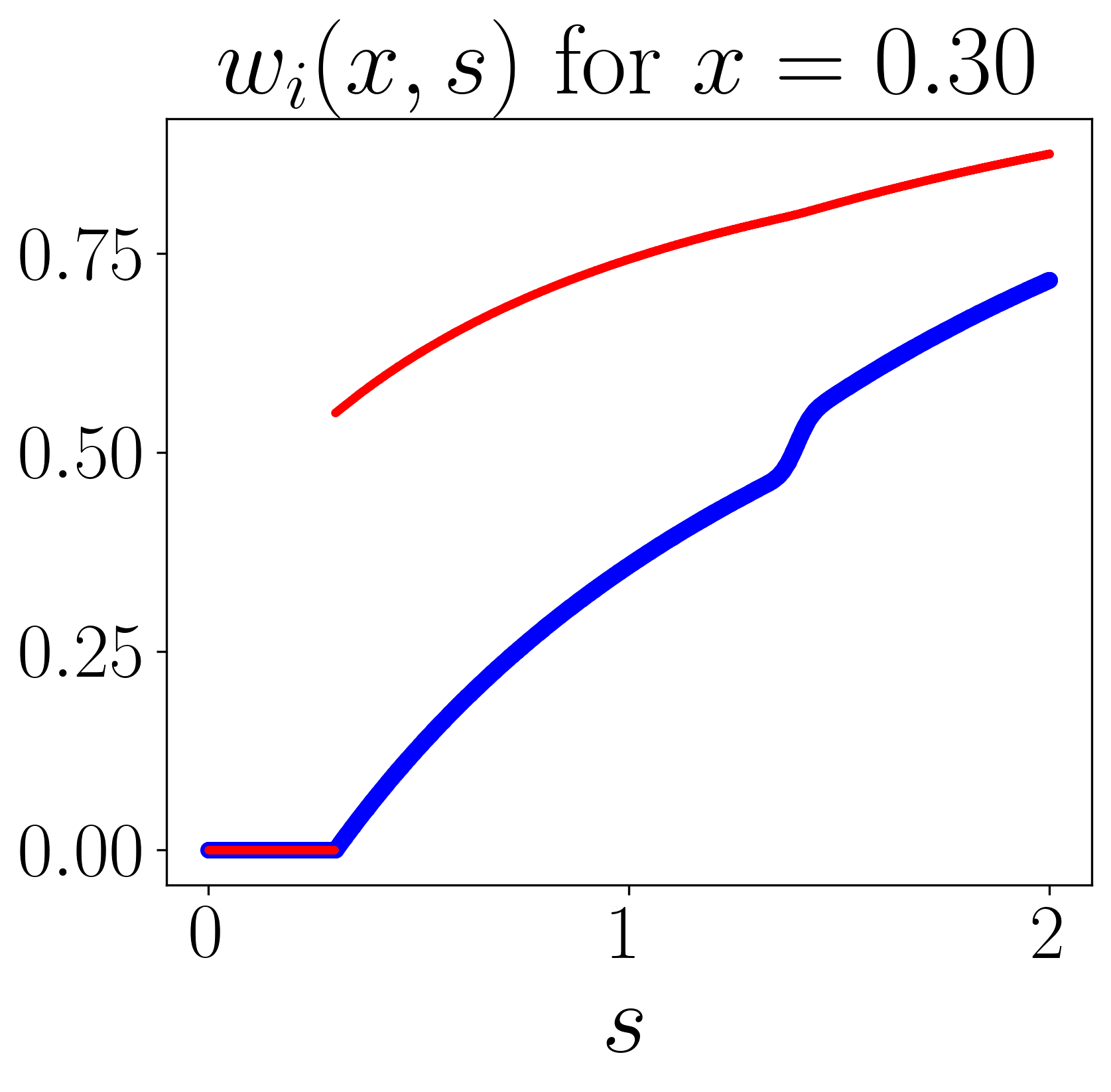} & \includegraphics[width=0.22\textwidth]{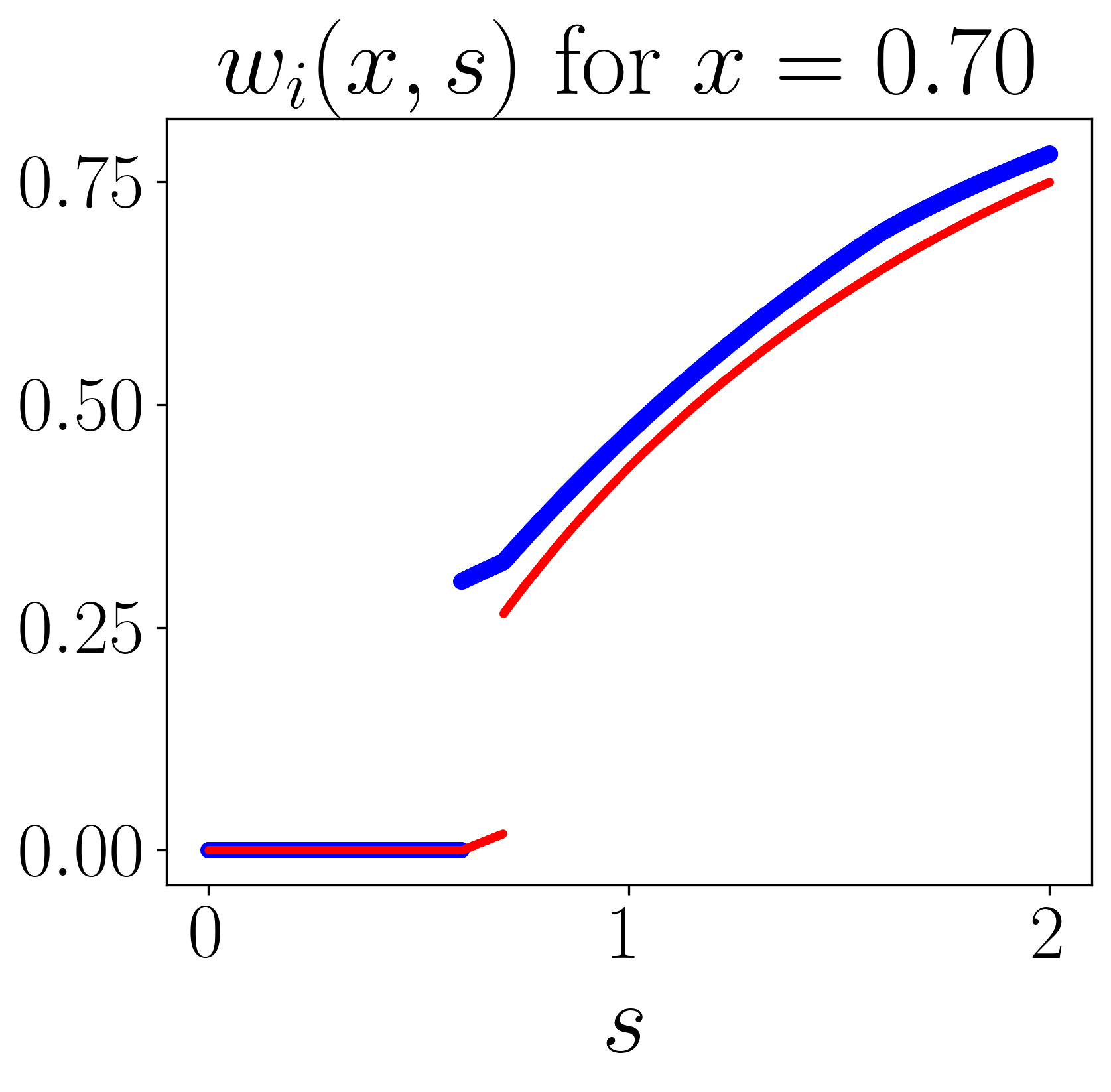} \\
	(A) & (B) & (C) & (D)
	\end{array}$
	\caption{CDF for a particle starting at initial position $x$, in Mode 1 (blue) and Mode 2 (red). In Subfigures (A) and (C), the initial condition is $x=0.30$ while in Subfigures (B) and (D) the initial condition is $x=0.70$. Subfigures (A) and (B) are for Example 1 and (C) and (D) are for Example 2.}
	\label{fig:cdfs}
\end{figure}

\begin{figure}
	\centering
	\begin{subfigure}{0.24\textwidth}
		\includegraphics[width=\textwidth]{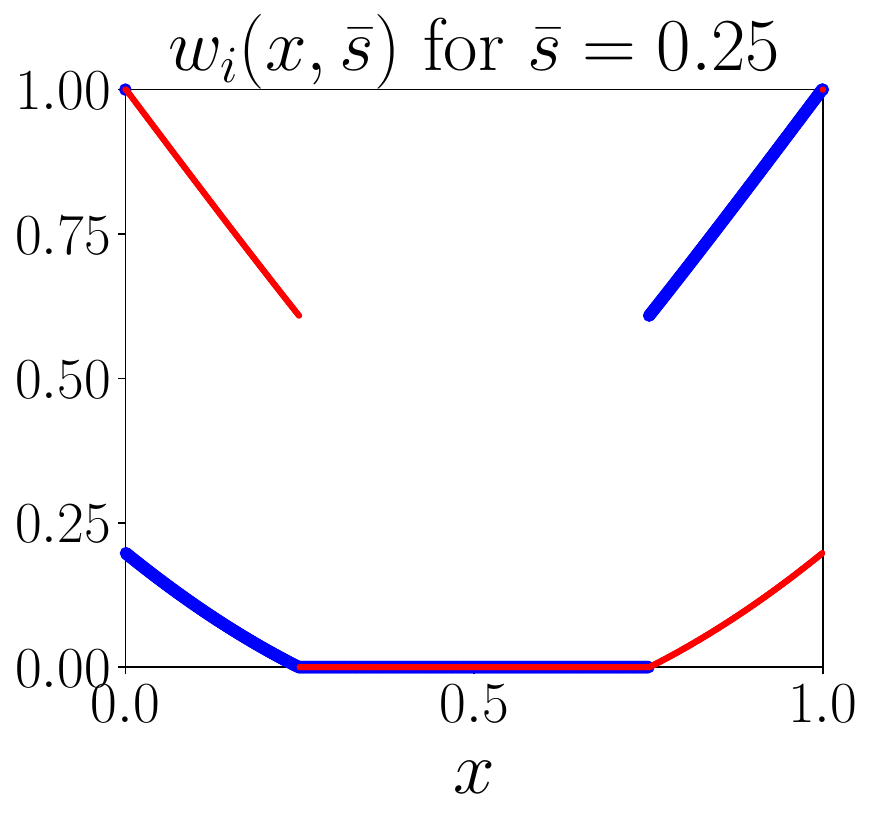}
	\end{subfigure}
	\begin{subfigure}{0.24\textwidth}
		\includegraphics[width=\textwidth]{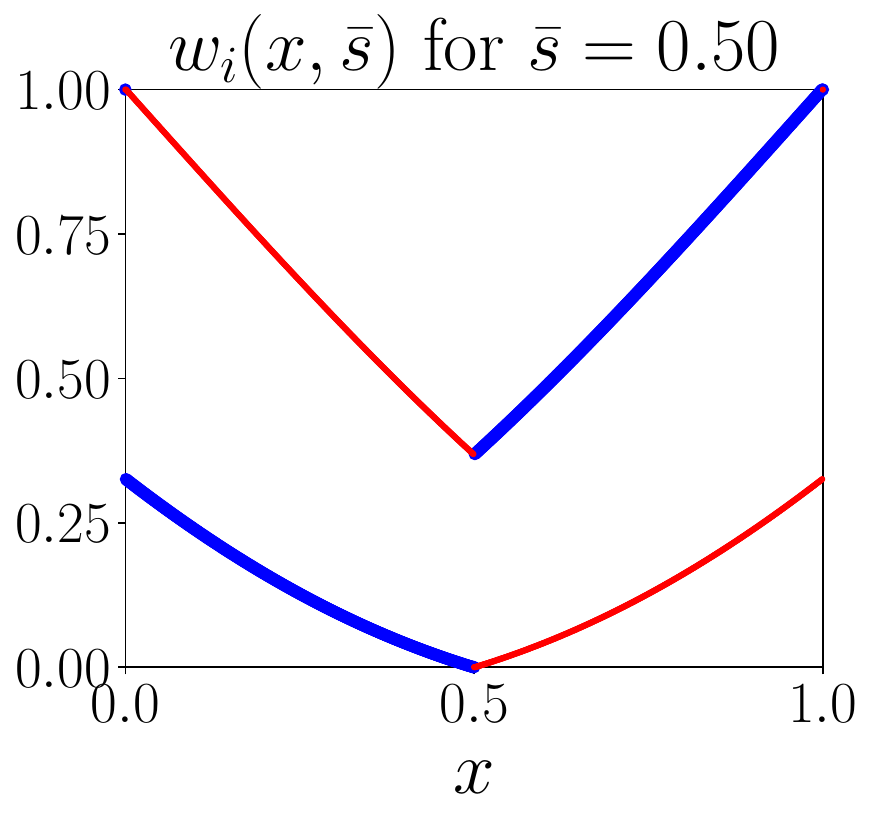}
	\end{subfigure}
	\begin{subfigure}{0.24\textwidth}
		\includegraphics[width=\textwidth]{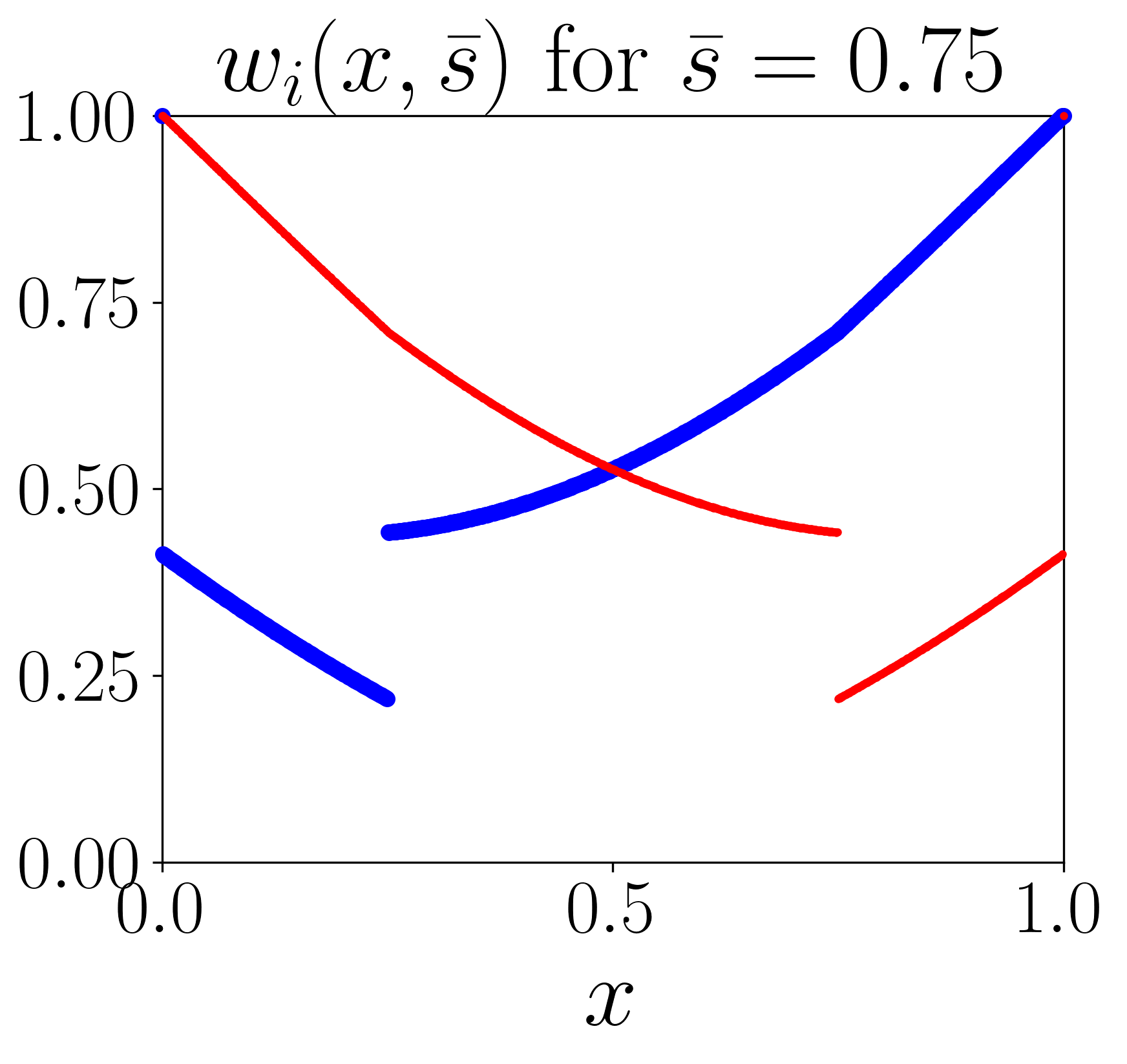}
	\end{subfigure}
	\begin{subfigure}{0.24\textwidth}
		\includegraphics[width=\textwidth]{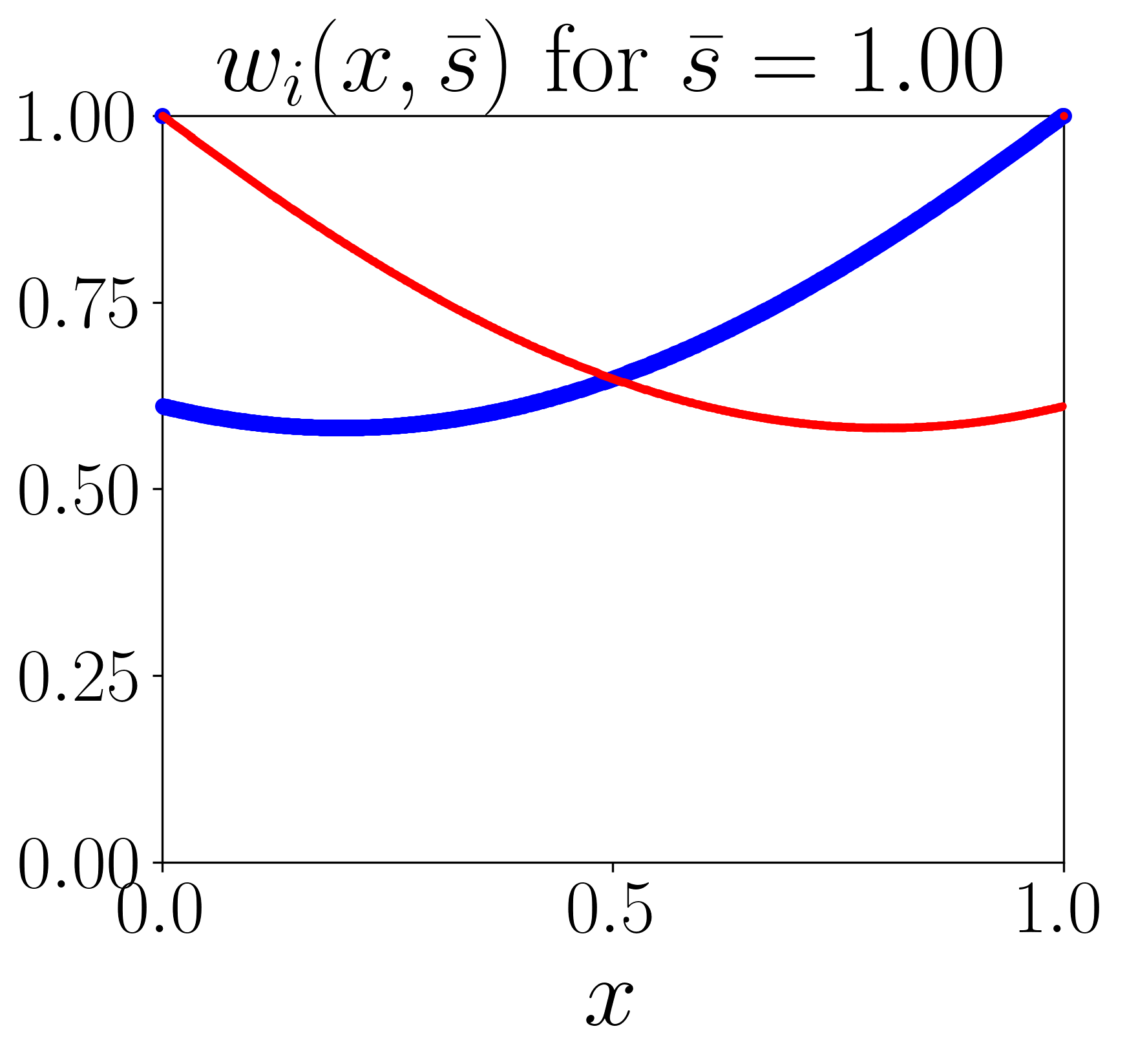}
	\end{subfigure}
	\caption{Example 1: Equal speeds, with a symmetric transition rate $\lambda = 2$. Each subplot is a snapshot of $w_i(x,s) = \PP(\J_i(\x) \le s)$ for a specific value of $s$. In Mode 1 (blue), the particle moves to the right with speed 1. In Mode 2 (red), the particle moves to the left with speed $1.$ Computed on $\Omega \times \mS = [0,1]^2$ with $\Delta x = \Delta s = 0.001$.}
	\label{fig:uniform}
\end{figure}

The key advantage of our approach is that it approximates the distribution $\J$ for all starting configurations simultaneously.
Once $w_i$'s are computed, we can freeze $(x,i)$ and vary $s$ to study the CDF.  In \cref{fig:cdfs}(A-B) this is illustrated for two starting locations $x=0.3$ and $x=0.7.$
But it might be even more revealing to fix a particular deadline $s$ and consider the probability of meeting it from all possible initial configurations.
In \cref{fig:uniform} we show such graphs of $w_i(x,s)$ for four different $s$ values.
Geometric properties of these functions have a natural interpretation, which we highlight focusing on mode 1 and $s=0.25$  
(the blue plot in the first subfigure). First, regardless of mode switches, $s=0.25$ is not enough time to exit if we start too far from $Q$; so, $w_1 = w_2 = 0$ for all $\x \in (0.25, 0.75).$  Second, starting from $x=0.75$ and moving right with speed one we will have just enough time to reach $Q$ provided we experience no mode switches, and if any switches occur the resulting time to target will be higher. So, the jump discontinuity at $x=0.75$ is precisely the probability of zero mode switches occurring in $s=0.25$ time units. (We note that this discontinuity disappears in the last subfigure since $s=1.00$ is enough time to reach $Q$ with no mode switches starting from any 
$(x,i) \in \domain \times \M.$) Finally, a similar argument explains the behavior for starting positions on $x \in (0, 0.25)$.  Since we start in mode 1, the only hope of meeting the $s=0.25$ deadline is a quick switch to mode 2.  Starting from $x=0.25,$ a timely arrival would require an immediate mode switch, and since this happens with probability zero, $w_1$ is continuous at this point. 

Of course, the probability of meeting a deadline is also significantly influenced by the switching rates.  While we do not illustrate this here, the same example is repeated with a range of symmetric and asymmetric rates in Figure \ref{fig:bounds_unequal2} of section \ref{sec:bounds}.


\begin{figure}[hbt]
	\centering
	\begin{subfigure}{0.24\textwidth}
		\includegraphics[width=\textwidth]{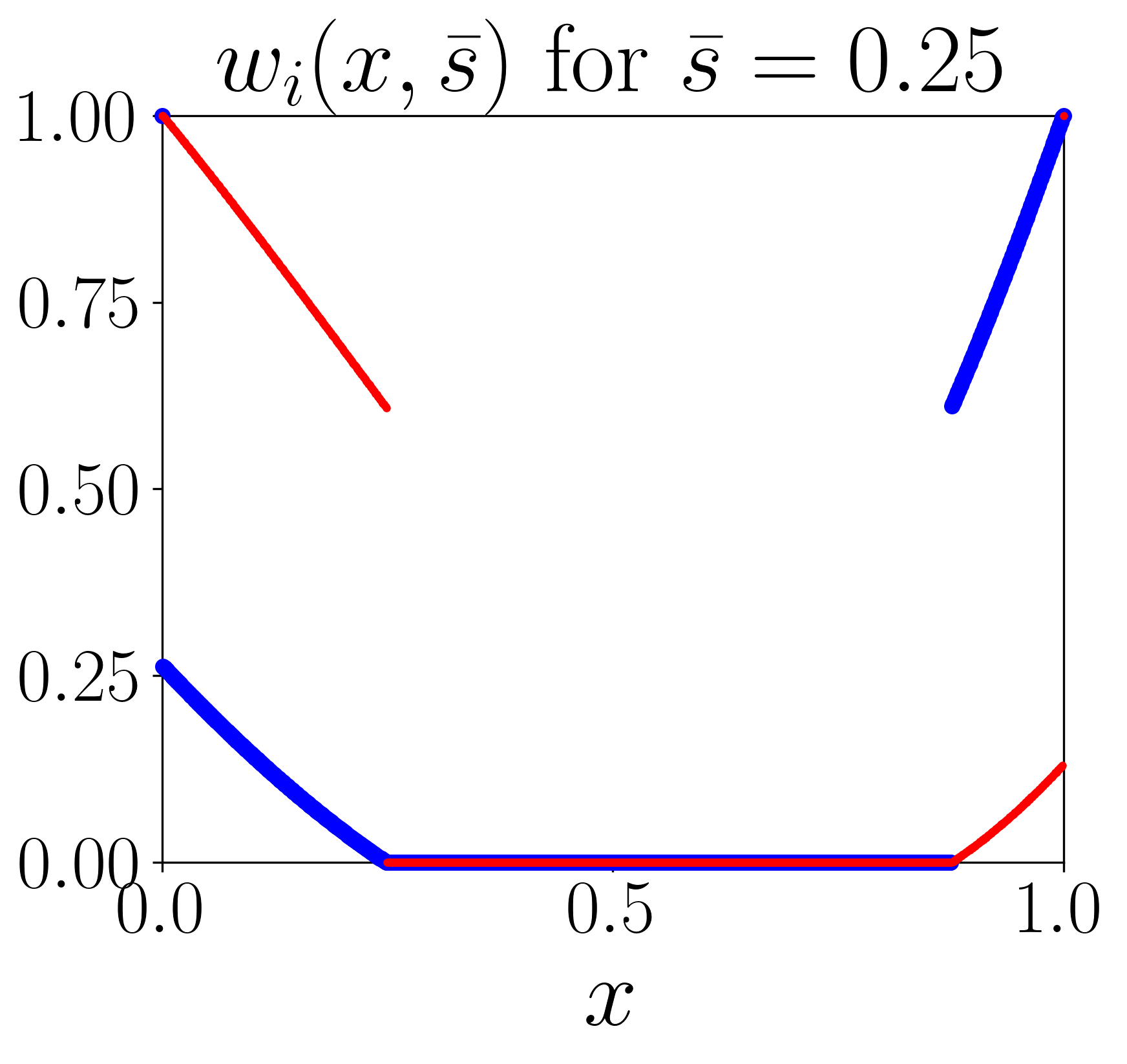}
	\end{subfigure}
	\begin{subfigure}{0.24\textwidth}
		\includegraphics[width=\textwidth]{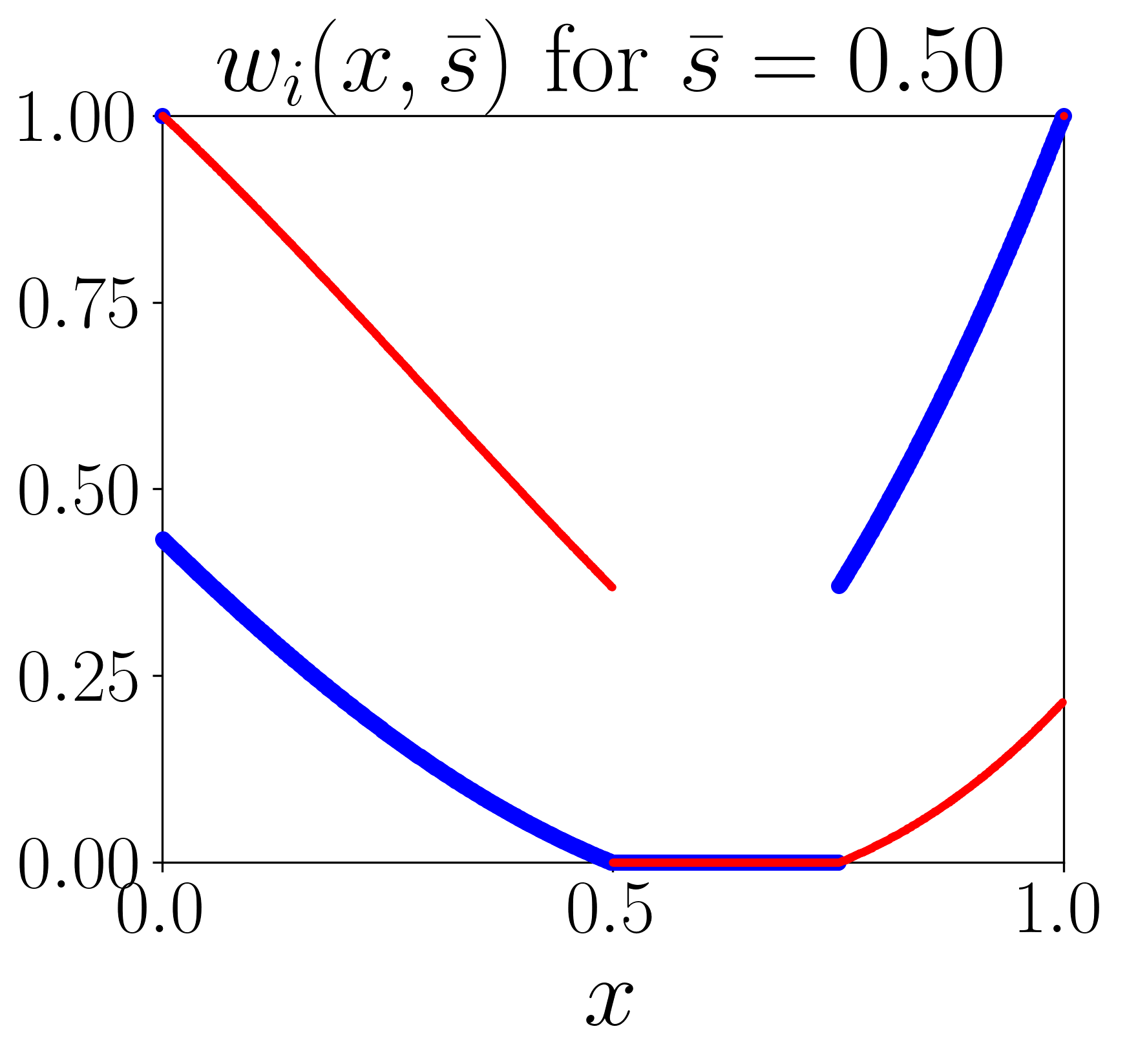}
	\end{subfigure}
	\begin{subfigure}{0.24\textwidth}
		\includegraphics[width=\textwidth]{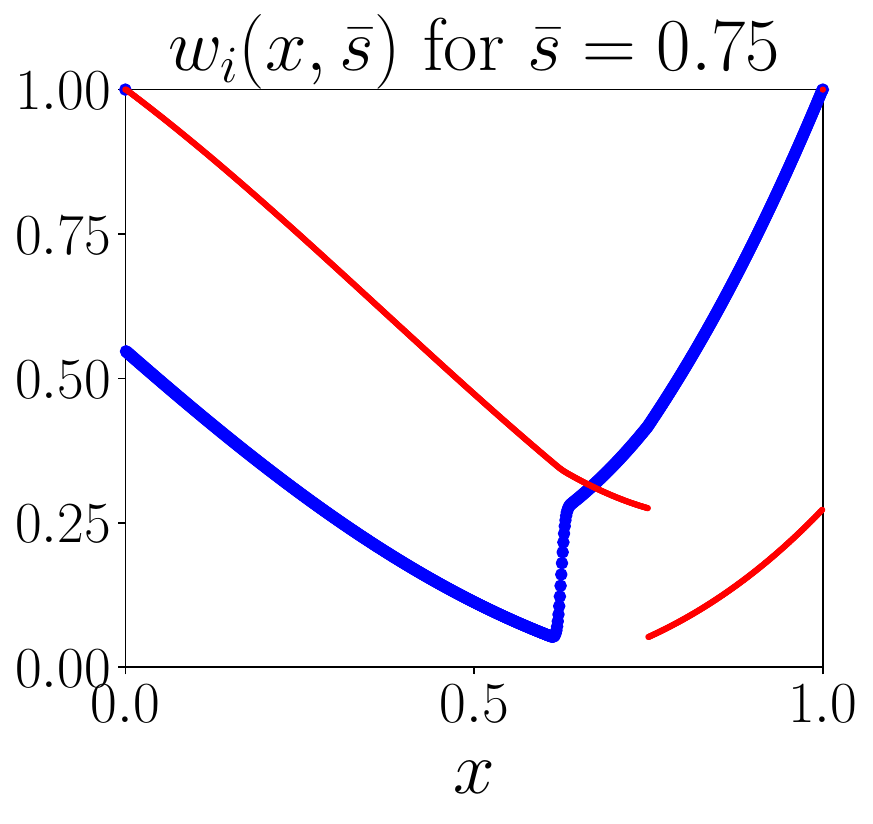}
	\end{subfigure}
	\begin{subfigure}{0.24\textwidth}
		\includegraphics[width=\textwidth]{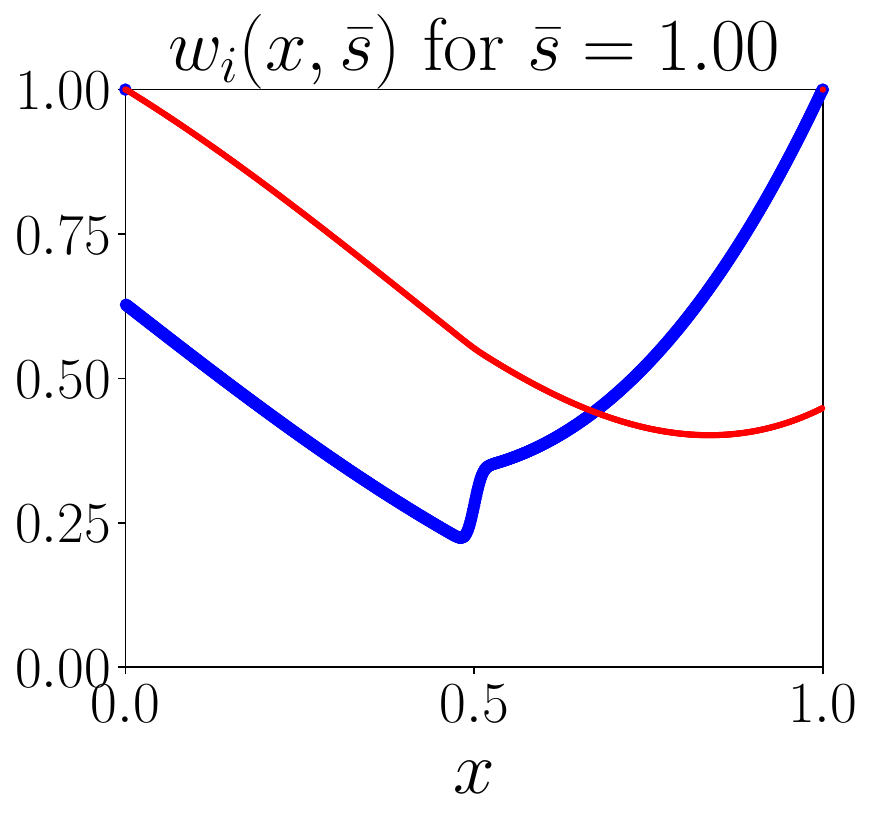}
	\end{subfigure}
	\caption{Example 2: Unequal speeds, with a transition rate $\lambda = 2$. Each subplot is a snapshot of $w_i(x,s) = \PP(\J_i(\x) \le s)$ for a specifc value of $s$. In Mode 1 (blue), the particle moves to the right with speed $1/2$. In Mode 2 (red), the particle moves to the left with speed 1. Computed on $\Omega \times \mS = [0,1]^2$ with $\Delta x = \Delta s = 0.001$.}
	\label{fig:unequal}
\end{figure}

{\bf Example 2: }
We modify the previous example by considering unequal speeds of motion in different modes: $f_1 = 0.5$ and $f_2 = -1.$ 
The CDFs for two starting locations $x=0.3$ and $x=0.7$ are shown in \cref{fig:cdfs}(C-D) while the plots of $w_i(x,s)$ for four different values of $s$ can be found in \cref{fig:unequal}.
We note that in Mode 1, interpolation is now necessary in \cref{eq:cdf_update}, which results in numerical diffusion smoothing out discontinuities, as can be seen in the right two subfigures. 
The absence of such artifacts in the first two subfigures is an additional benefit of pre-computing $s^0(x)$ and $w^0_i(x)$, shown in \cref{fig:best_case}(C-D), to reduce the computational domain for $w_i$'s.
E.g., for the first subfigure, all $x \in [0.25, 0.875)$ have $s^0(x) > s = 0.25$ and so are assigned an exit probability of 0, removing the need of interpolating across discontinuities. 
In contrast, at $s = 0.75$ all $x$ have a nonzero probability of exiting, so interpolation across the discontinuity at $x = 0.625$ is unavoidable.

\begin{figure}
\begin{tabular}{>{\centering\arraybackslash}m{0.5in}>{\centering\arraybackslash}m{0.85in}>{\centering\arraybackslash}m{0.85in}>{\centering\arraybackslash}m{0.85in}>{\centering\arraybackslash}m{0.85in}}
& $w_1  \quad \leftarrow$ & $w_2 \quad \uparrow$ & $w_3 \quad \rightarrow$ & $w_4 \quad \downarrow$ \\
s = 0.25 & \includegraphics[height=0.16\textwidth]{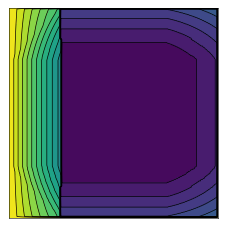} & \includegraphics[height=0.16\textwidth]{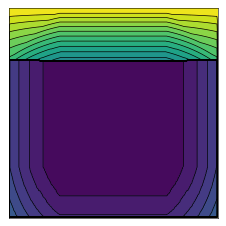} & \includegraphics[height=0.16\textwidth]{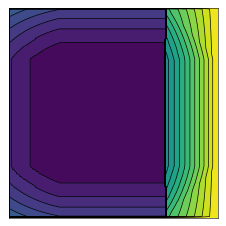} & \includegraphics[height=0.16\textwidth]{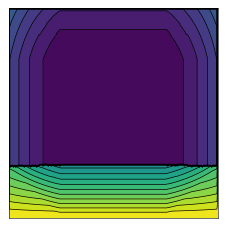} \\
s = 0.50 & \includegraphics[height=0.16\textwidth]{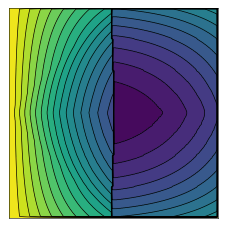} & \includegraphics[height=0.16\textwidth]{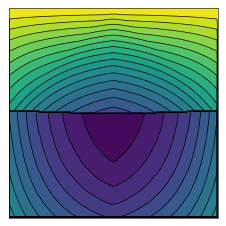} & \includegraphics[height=0.16\textwidth]{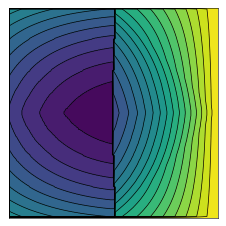} & \includegraphics[height=0.16\textwidth]{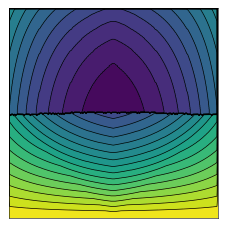} \\
s = 0.75 & \includegraphics[height=0.16\textwidth]{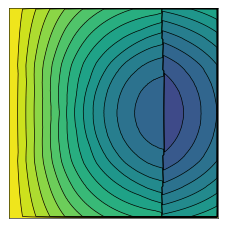} & \includegraphics[height=0.16\textwidth]{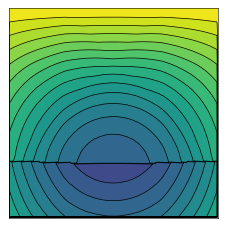} & \includegraphics[height=0.16\textwidth]{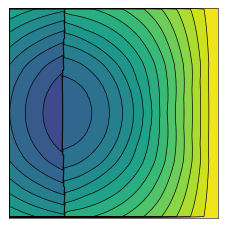} & \includegraphics[height=0.16\textwidth]{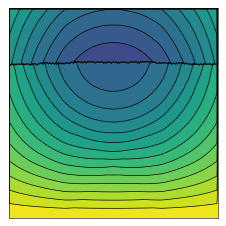} \\
s = 1.00 & \includegraphics[height=0.16\textwidth]{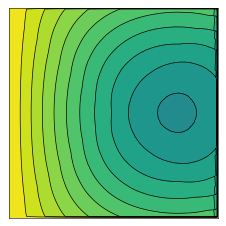} & \includegraphics[height=0.16\textwidth]{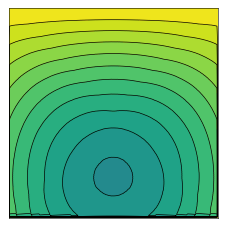} & \includegraphics[height=0.16\textwidth]{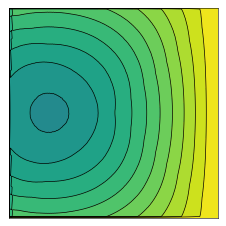} & \includegraphics[height=0.16\textwidth]{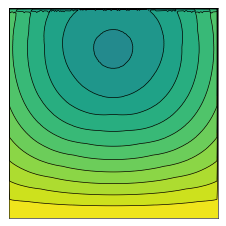}
\end{tabular}
\begin{center}
\qquad\includegraphics[width=0.5\textwidth]{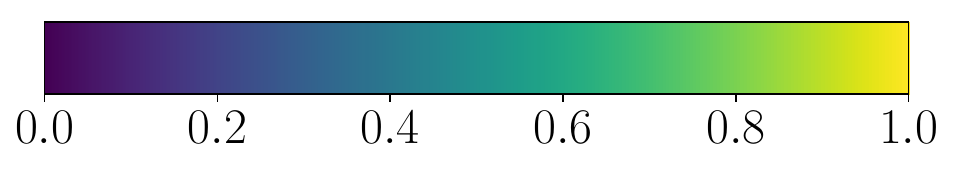}
\end{center}
\caption{Example 3: Mode switching in 2D with transition rates $\lambda = 1$. The particle moves $\leftarrow$ in Mode 1, $\uparrow$ in Mode 2, $\rightarrow$ in Mode 3, and $\downarrow$ in Mode 4. Computed on $\domain \times \mS = [0,1]^3$ with $\Delta x = \Delta y = \Delta s = 0.01$.}
\label{fig:two_dimensional}
\end{figure}

{\bf Example 3: }
We now consider a 2D version of Example 1, with $\domain = [0,1] \times [0,1],$ $Q = \boundary,$ $M=4,$ and $\lambda = 1$.
In all modes, the motion is with speed $|\Bf| = 1$, but the directions of motion differ:
$\leftarrow,\uparrow,\rightarrow,$ and $\downarrow$ in modes $1, \ldots, 4$ respectively.
Numerical approximations of $w_i$'s for different values of $s$ are shown in \cref{fig:two_dimensional}. 
The distinct delineations between darker and lighter regions are analogous to the discontinuities in the earlier one-dimensional cases. 
For example, given $s < 0.5$ and starting positions along the line $y = 0.5$, a timely exit is only possible to the left (via Mode 1) or to the right (via Mode 3). 
Therefore, cross sections of $w_1$ and $w_3$ along $y = 0.5$ at $s = 0.25$ in \cref{fig:two_dimensional} coincide with the one-dimensional graphs for $s = 0.25$ in \cref{fig:uniform}. 
However, as we move closer to the corners of the domain, all four modes have an effect on the probability of exit. 
For example, the region along the diagonal near the top right corner of the $w_1$ graph has higher exit probabilities than surrounding regions because there are multiple possible timely-exit strategies.  These 2D phenomena become prevalent for higher $s$ values.

\begin{remark} [{\bf Related models in biology}]
\label{rem:bio_dispersal}
As noted in the introduction, similar ``velocity jump processes'' are also used to model dispersal in  biological systems \cite{othmer1988models, HillenDiffLimit}.
In that context, all dispersing agents perform long runs with constant velocity but occasionally switch modes/directions.  The usual approach is to derive a system of PDEs governing the evolution 
of agent densities $\rho_i(\x,t)$ in corresponding modes $i \in \M.$  
The symmetric unbounded case in 1D (i.e., $\domain = \R, \, M =2,$ and $\lambda_{12} = \lambda_{21}$) is particularly well-studied, with 
the overall density $\rho = \rho_1 + \rho_2$ evolving according to the ``telegraph equation'' \cite{goldstein1951diffusion}.  
Taking $\lambda_{12} \neq \lambda_{21},$ one can similarly model chemotaxis.  If $\domain = \R^2$ or $\R^3$, one could use a larger number of modes to describe many possible directions of motion, with  $\lambda_{ij}$ chosen to reflect a possible bias in switching (e.g., giving preference to new directions more closely aligned with the preceding run -- as is the case for {\em E. coli} bacteria).  Letting $M \to \infty,$ one can also directly model all possible directions of motion 
by switching to integro-differential equations \cite{othmer1988models}.

While our focus on a single performance measure $\J$ might be restrictive for many of these applications, there are also some settings where it can be advantageous.
For example, if one assumes that agents are removed upon reaching $Q$,  the number of them still remaining by the time $t$ 
could be in principle computed as 
\begin{equation}
\label{eq:R_naive}
R(t) \; = \; \sum\limits_{i \in \M} \int_{\domain \backslash Q} \rho_i(\x,t) \, \bm{dx}.
\end{equation}
But any change in $\rho_i(\x,0)$ would make it necessary to re-solve a system of PDEs for $\rho_i(\x,t)$'s before reusing \eqref{eq:R_naive}.
Here we can offer a much more efficient method by setting $C \equiv 1$ and $q \equiv 0,$
computing $w_i$'s from \eqref{eq:cdf_pde} only once, and then using an alternative formula that works for all initial densities
\begin{equation}
\label{eq:R_smart}
R(t) \; = \; \sum\limits_{i \in \M} \int_{\domain \backslash Q} \rho_i(\x,0) \left( 1 - w_i(\x,t) \right) \, \bm{dx}.
\end{equation}
\end{remark}

\section{Bounds on CDF}
\label{sec:bounds}
We now turn to PDMPs with parameter uncertainty -- in addition to the inherent aleatoric uncertainty due to mode switches. In \cref{sec:compute}, the uncertainty of the outcome could be fully characterized by its CDF computed based on the known transition rates between modes, $\lambda_{ij}$'s. 
Here, however, we consider the case where we only know a {\em range} of potential $\lambda_{ij}$ values. There are two natural models of epistemic uncertainty in this situation, and it is meaningful to consider the upper and lower bounds on the CDF with each of them. We focus on a case where the true transition rates are free to fluctuate within the given range and may take on different values at different times. The upper and lower bounds on the CDF can be then found by considering a nonlinear version of the coupled PDEs seen in \cref{sec:compute}. 
This can also be viewed as an optimal control problem, where the controller is either helping or hindering the particle's exit by choosing the transition rates adaptively. 

The alternative model of epistemic uncertainty is to assume that all transition rates remain fixed (though unknown) throughout the process.  We provide some experimental results for this case as well, though do not propose any computationally efficient methods for finding sharp CDF bounds.

\subsection{Deriving PDEs}
\label{subsec:bounds_pde}

We now extend the results of section \ref{subsec:continuous} by considering the case in which the transition rate matrix \(\Lambda = (\lambda_{ij})\) is not necessarily constant.
Suppose there are known \(a_{ij}, b_{ij}\) for each pair \(i \neq j\), such that each \(\lambda_{ij}\) may vary in the interval \(0 \leq a_{ij} \leq \lambda_{ij} \leq b_{ij}\) throughout the process.
If $L$ is the set of possible transition matrices satisfying these constraints, we will assume that $\Lambda$ might be changing but remains in $L$ throughout the process. 
In this section, we will use $\Lambda_i$ to denote the $i$-th row of $\Lambda$, specifying all transition rates from mode $i$. 
We will also use $L_i$ to denote the set of all allowable $i$-th rows satisfying the above constraints.

We compute an upper bound for \(w_i\), denoted \(w_i^+\), by taking its initial and boundary conditions to be the same as \(w_i\), and 
adaptively selecting the \(\Lambda_i \in L_i\) which maximizes \(w_i^+(\x, s)\). 
Similarly, for the lower bound \(w_i^-\) we take the \(\Lambda_i \in L_i\) which minimizes \(w_i^-(\x, s)\). 
Hence, for each mode $i$, the bounds 
$w_i^+$ and $w_i^-$ satisfy the PDEs:
\begin{equation}
\label{eq:w_plus}
\nabla w_i^+(\x,s)\cdot\Bf_i(\x) - C_i(\x) \frac{\partial w_i^+}{\partial s}(\x,s) + \max_{\Lambda_i \in L_i}  \Big\{ \sum_{j\neq i}  \lambda_{ij}\left[w_j^+(\x,s)-w_i^+(\x,s)\right]\Big\} \; = \; 0;
\end{equation}
\begin{equation}
\label{eq:w_minus}
    \nabla w_i^-(\x,s)\cdot\Bf_i(\x) - C_i(\x) \frac{\partial w_i^-}{\partial s}(\x,s) + \min_{\Lambda_i \in L_i}\Big\{ \sum_{j\neq i}  \lambda_{ij}\left[w_j^-(\x,s)-w_i^-(\x,s)\right]\Big\} \; = \; 0
\end{equation}
with initial and boundary conditions  \eqref{eq:cdf_ic} and \eqref{eq:cdf_bc}.
We note that, if $a_{ij}=b_{ij}$ for all $i \neq j,$ there is no parametric uncertainty, each $L_i$ is a singleton, and the above equations reduce to \eqref{eq:cdf_pde}.

As with the computation of $w_i$'s through \eqref{eq:cdf_pde}, it can be helpful to precompute the minimal attainable cost \(s^0(\x)\) and the probability \(w_i^0(\x)\) of achieving such a cost. 
Since \(\Lambda\) is not constant, instead of having a probability of attaining the cost \(s^0(\x)\) we have a lower bound \(w^{0,-}_i(\x)\) and an upper bound \(w^{0,+}_i(\x)\) for that probability.
We may compute \(s^0(\x)\) in precisely the same way as in \eqref{eq:cont_best_case} since that formula does not depend on \(\Lambda\) at all.
On the other hand, to compute \(w_i^{0, +}(\x)\) we must modify \eqref{eq:cont_best_case_prob} to account for the unknown (and possibly changing) \(\Lambda\):
\begin{align}\label{eq:cont_best_case_prob_upper}
\nonumber &0 = \nabla w_i^{0,+}(\x)\cdot \Bf_i(\x) + \max_{\Lambda_i \in L_i}  \Big\{ \sum_{j \not= i}\lambda_{ij} \left[ w_j^{0,+}(\x) - w_i^{0,+}(\x)\right]\Big\}, &&\x \in \Omega \setminus Q, i \in \I(\x); \\
\nonumber &w_i^{0,+}(\x) = 1, &&\x \in Q, i \in \I(\x); \\
&w_i^{0,+}(\x) = 0, &&\x \in \Omega, i \not\in \I(\x).
\end{align}
Similarly, to compute \(w_i^{0,-}(\x)\) we have:
\begin{align}\label{eq:cont_best_case_prob_lower}
\nonumber &0 = \nabla w_i^{0,-}(\x)\cdot \Bf_i(\x) +\min_{\Lambda_i \in L_i}   \Big\{ \sum_{j \not= i}\lambda_{ij} \left[ w_j^{0,-}(\x) - w_i^{0,-}(\x)\right]\Big\}, &&\x \in \Omega \setminus Q, i \in \I(\x); \\
\nonumber &w_i^{0,-}(\x) = 1, &&\x \in Q, i \in \I(\x); \\
&w_i^{0,-}(\x) = 0, &&\x \in \Omega, i \not\in \I(\x).
\end{align}
In both cases, \(\I(\x)\) is the argmin set of \cref{eq:cont_best_case} as in \cref{eq:cont_best_case_prob}.

\subsection{Calculating Bounds}
\label{subsec:bounds_calculation}

For numerical computations of the bounds described in \cref{subsec:bounds_pde}, we rely on a discretization similar to that presented in \cref{subsec:cdf_numerics}. When $\tau$ is small enough that the approximations in \eqref{eq:probabilities} can be made, this optimization can be written recursively. 
For the CDF lower bound $w^-_i(\x,s)$, the semi-Lagrangian scheme is
\begin{equation}\label{eq:linear_program_lower}
    w^-_i(\x,s) = w^{-}_i(\tilde{\x},\tilde{s}) + \min_{\Lambda_i \in L_i} \Big\{ \tau\sum_{j\not=i} \lambda_{ij}\left[w^-_j(\tilde{\x},\tilde{s})-w^-_i(\tilde{\x},\tilde{s})\right]\Big\},
\end{equation}
where $\tilde{\x} = F_i(\x)$ and  $\tilde{s} = s - \tau C_i(\x)$.
Recall from \cref{eq:tau_ineq2} that $C_i(\x) > 0$ so it is always the case that $\tilde{s} < s$. Therefore, $ \w^{-}_i(\tilde{\x},\tilde{s})$ has already been calculated and so can be used in the computation of $ \w^{-}_i(\x,s)$. \\

For an efficient implementation of \cref{eq:linear_program_lower}, the optimal $\Lambda^* \in L$  can be found explicitly.
When minimizing $w^-_i(\x, s)$,  we would naturally like to subtract as much as possible and add as little as possible to $w^{-}_i\big(F_i(\x), \, s - \tau C_i(\x)\big) = w^{-}_i(\tilde{\x},\tilde{s})$. Therefore, 
\begin{align}
    \nonumber
    \left[w^-_j(\tilde{\x},\tilde{s})-w^-_i(\tilde{\x},\tilde{s})\right] \leq 0 &\implies \lambda^*_{ij} = b_{ij}; \\
    \label{eq:optimization_lower}
    \left[w^-_j(\tilde{\x},\tilde{s})-w^-_i(\tilde{\x},\tilde{s})\right] > 0 &\implies \lambda^*_{ij} = a_{ij}. 
\end{align}

For the the CDF upper bound $w^+_i,$ the scheme is similar modulo replacing $\min$ with $\max$ in \cref{eq:linear_program_lower}
and flipping the signs of inequalities in \cref{eq:optimization_lower}.

\subsection{Experimental Results}
\label{subsec:bounds_results}


\subsubsection*{Example 4: CDF bounds and comparison to fixed-$\Lambda$ CDFs}
\begin{figure}
	\centering
	\begin{subfigure}{0.48\textwidth}
		\includegraphics[width=\textwidth]{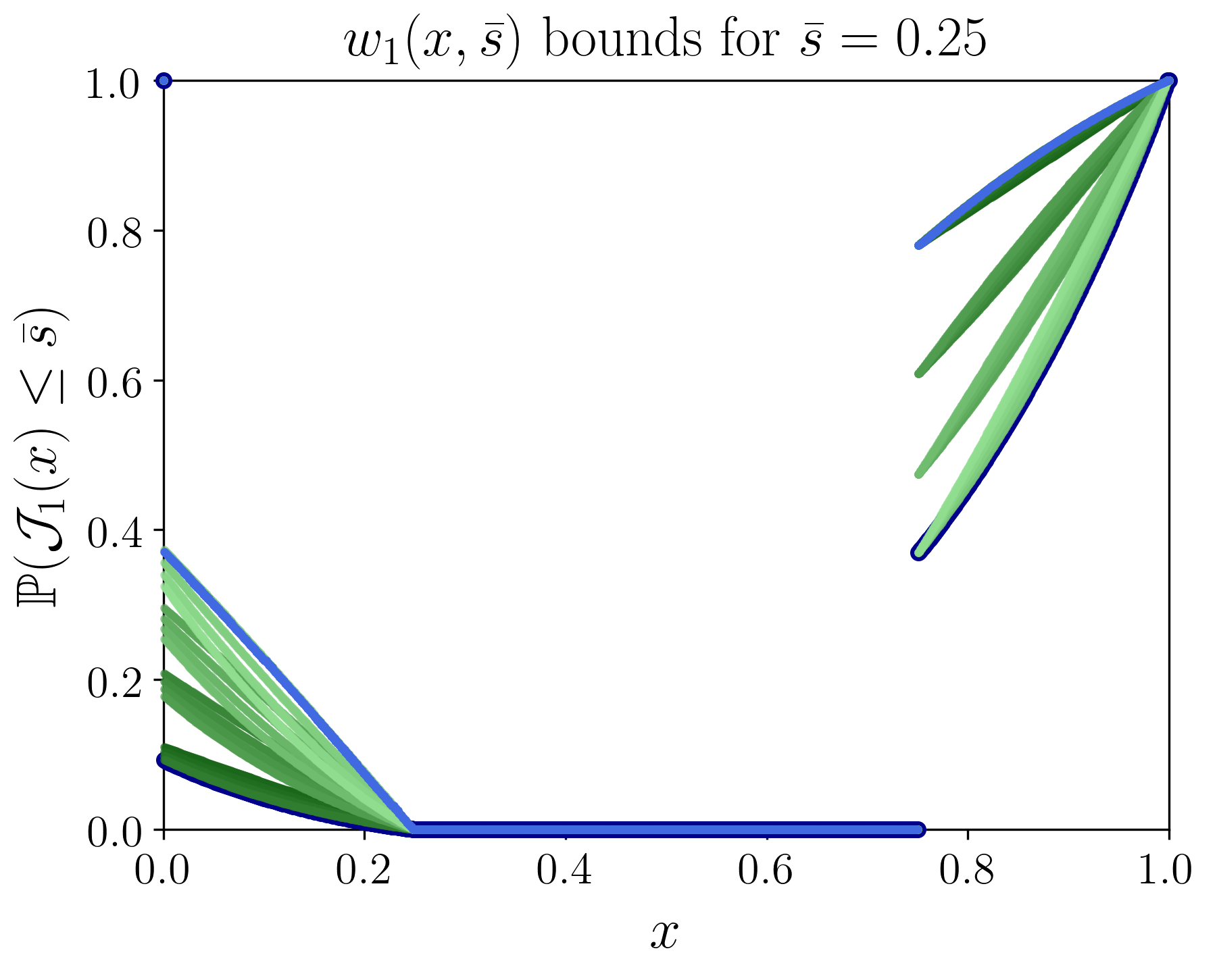}
	\end{subfigure}
	\begin{subfigure}{0.48\textwidth}
		\includegraphics[width=\textwidth]{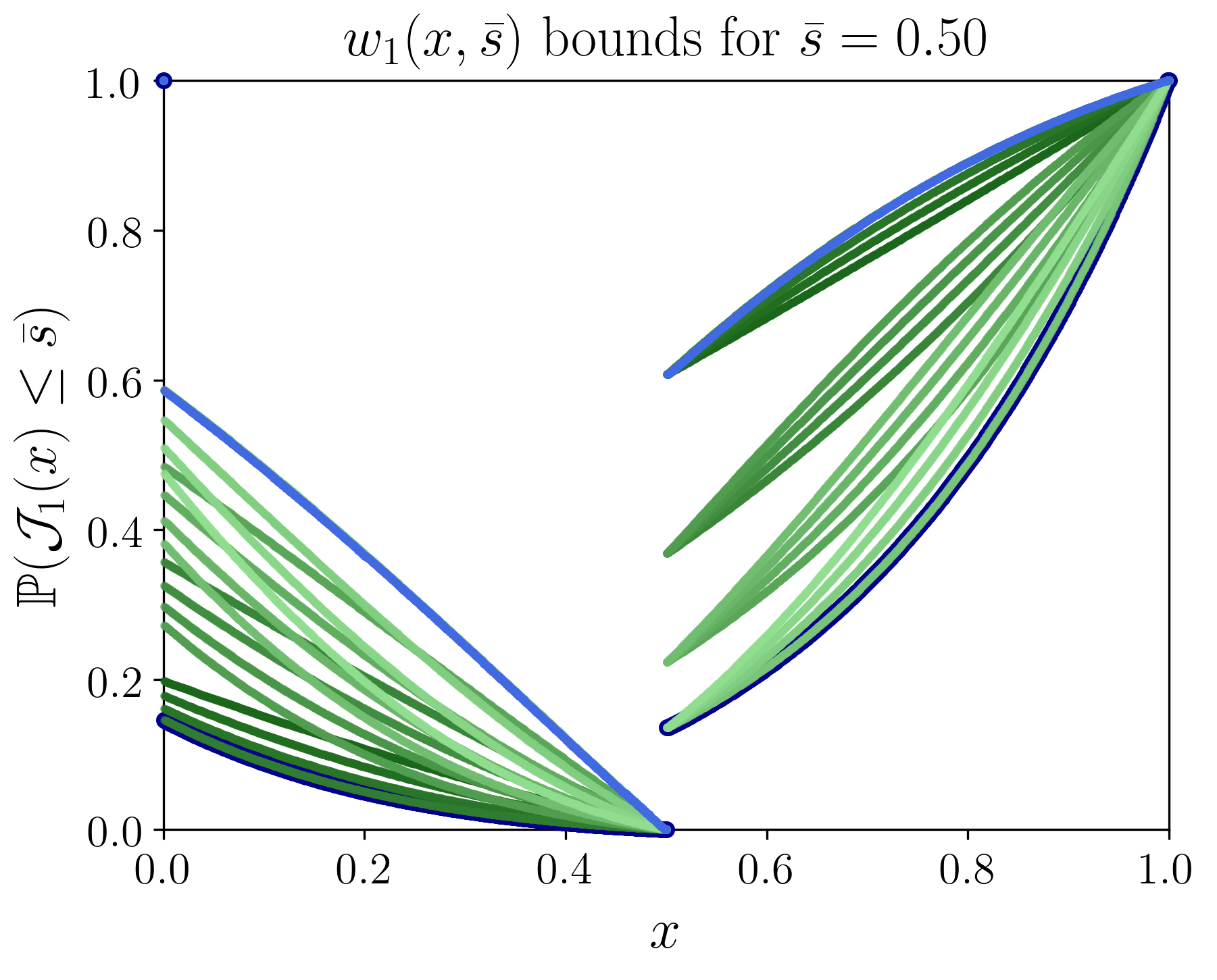}
	\end{subfigure}
	\begin{subfigure}{0.48\textwidth}
		\includegraphics[width=\textwidth]{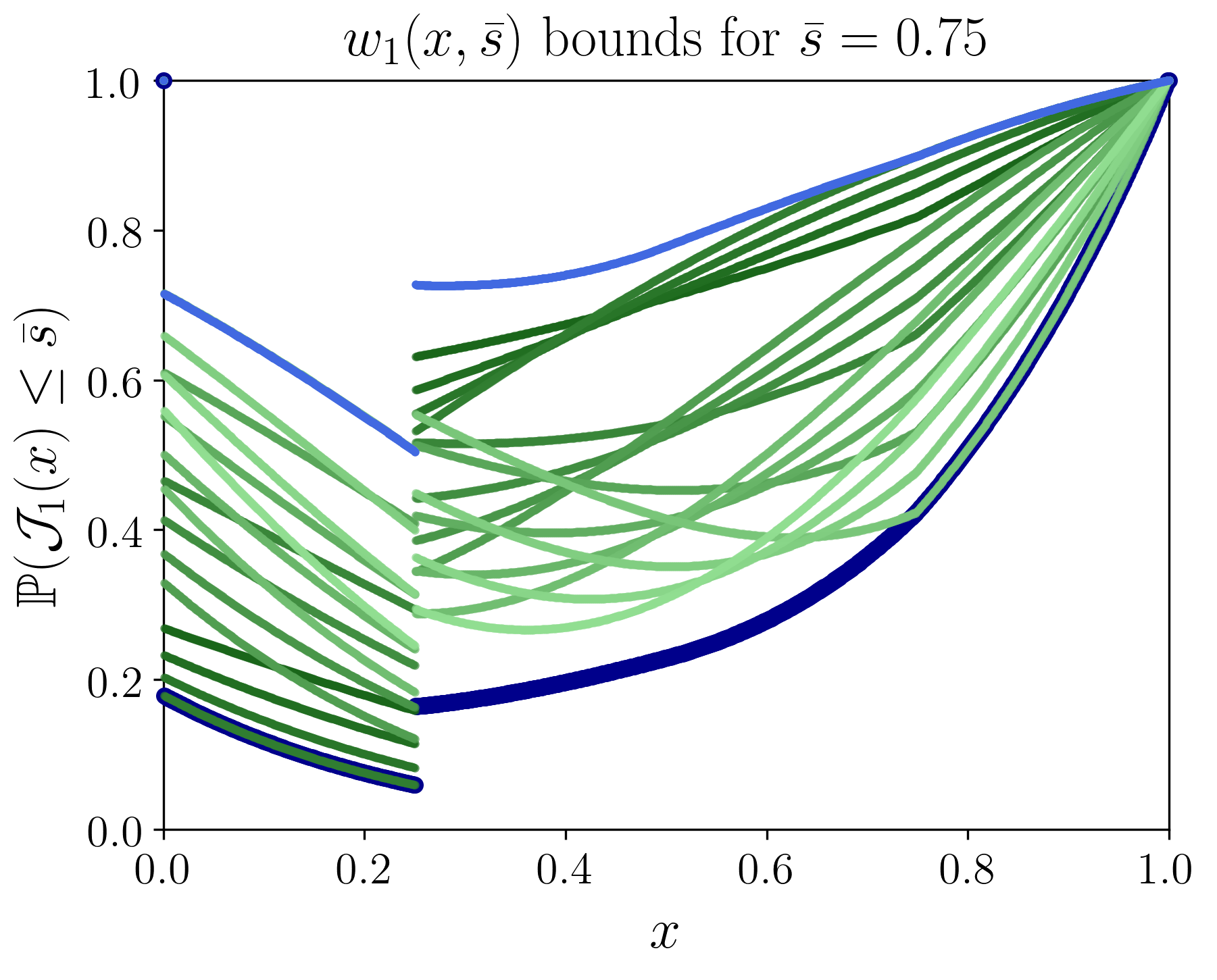}
	\end{subfigure}
	\begin{subfigure}{0.48\textwidth}
		\includegraphics[width=\textwidth]{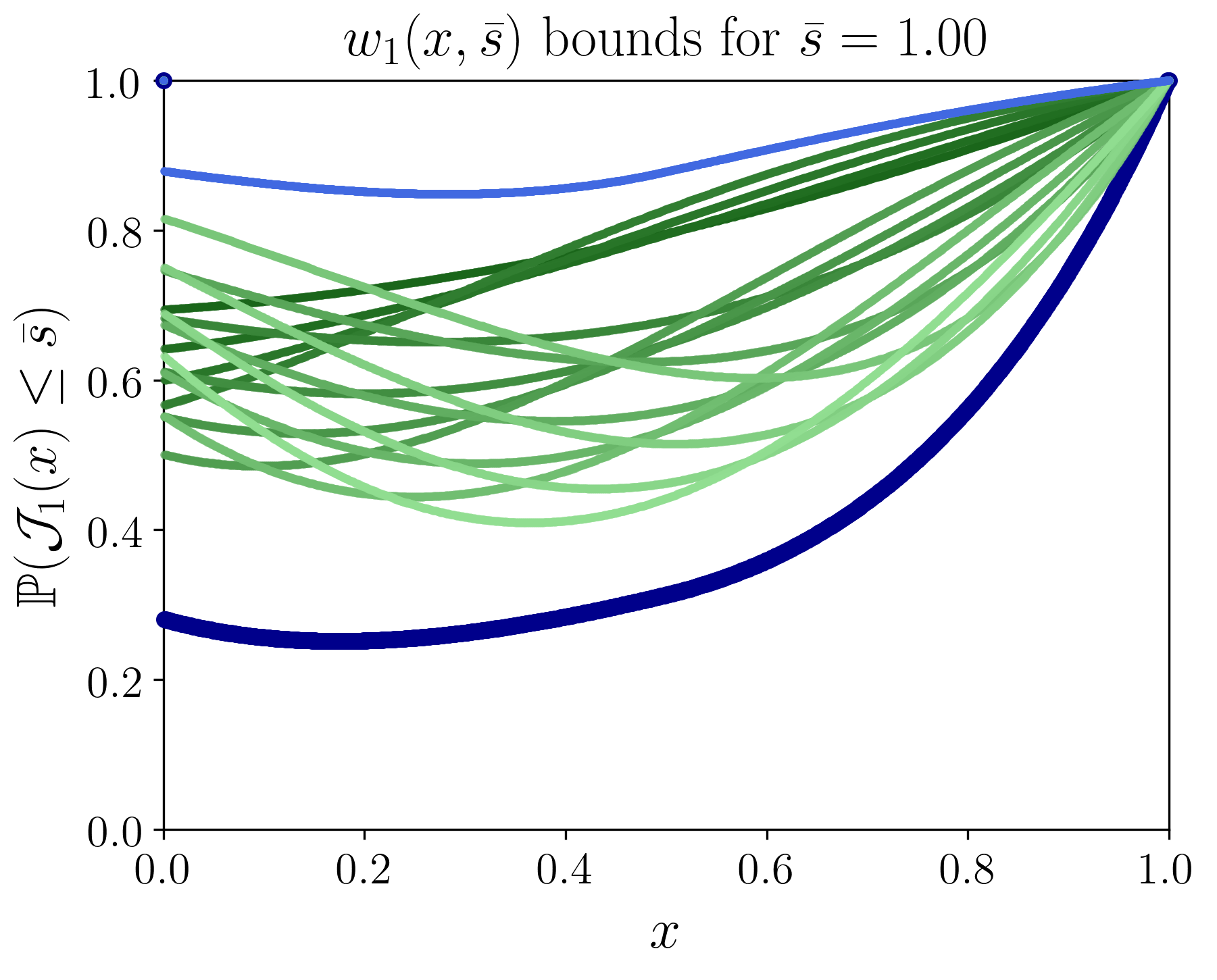}
	\end{subfigure}
	\caption{Example 4: bounds on the probability of timely exit (starting in Mode 1) for four different deadline values $\bar{s}.$ 
	(Probabilities and bounds for starting in Mode 2 can be obtained by a mirror symmetry relative to the line $x=0.5$.)
	Blue bounds are produced under the varying rates assumption by solving \cref{eq:w_plus}-\cref{eq:w_minus} for $\lambda_{12}, \lambda_{21} \in [1,4].$
	Green curves are produced under the fixed rates assumption by solving \cref{eq:cdf_pde}, 
	each corresponding to a specific $(\lambda_{12},\lambda_{21}) \in\{1, 2, 3, 4\}\times\{1, 2, 3, 4\}.$ 
	The darkest four curves are those associated with $\lambda_{12} = 1$, the next four are those associated with $\lambda_{12} = 2$, and so on. 
	Computed on $\domain \times \mS = [0,1]^2$ with $\Delta x = \Delta s = 0.001$.}
	\label{fig:bounds_unequal2}
\end{figure}

\begin{figure}
	\centering
	\begin{subfigure}{0.48\textwidth}
		\includegraphics[width=\textwidth]{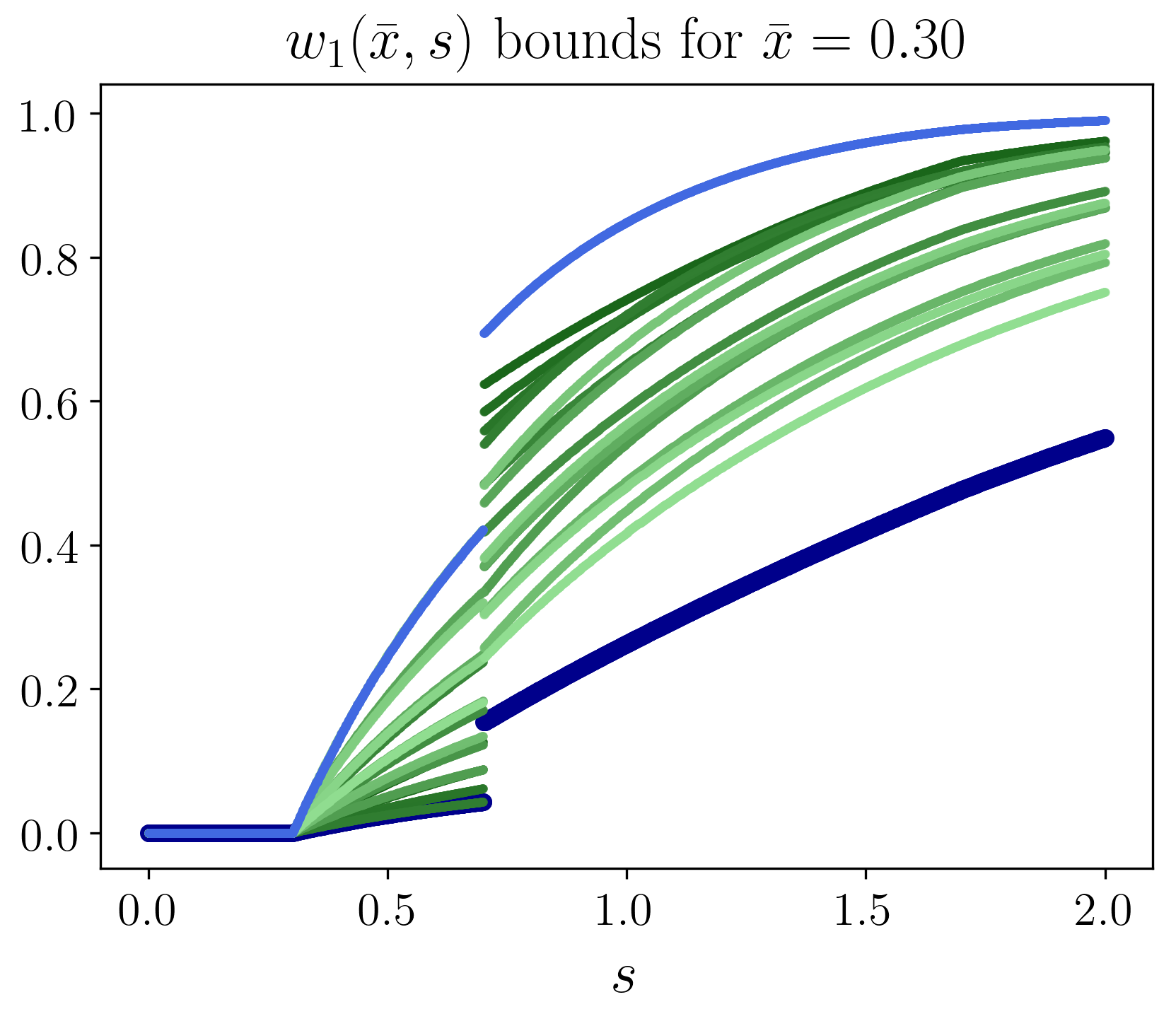}
	\end{subfigure}
	\begin{subfigure}{0.48\textwidth}
		\includegraphics[width=\textwidth]{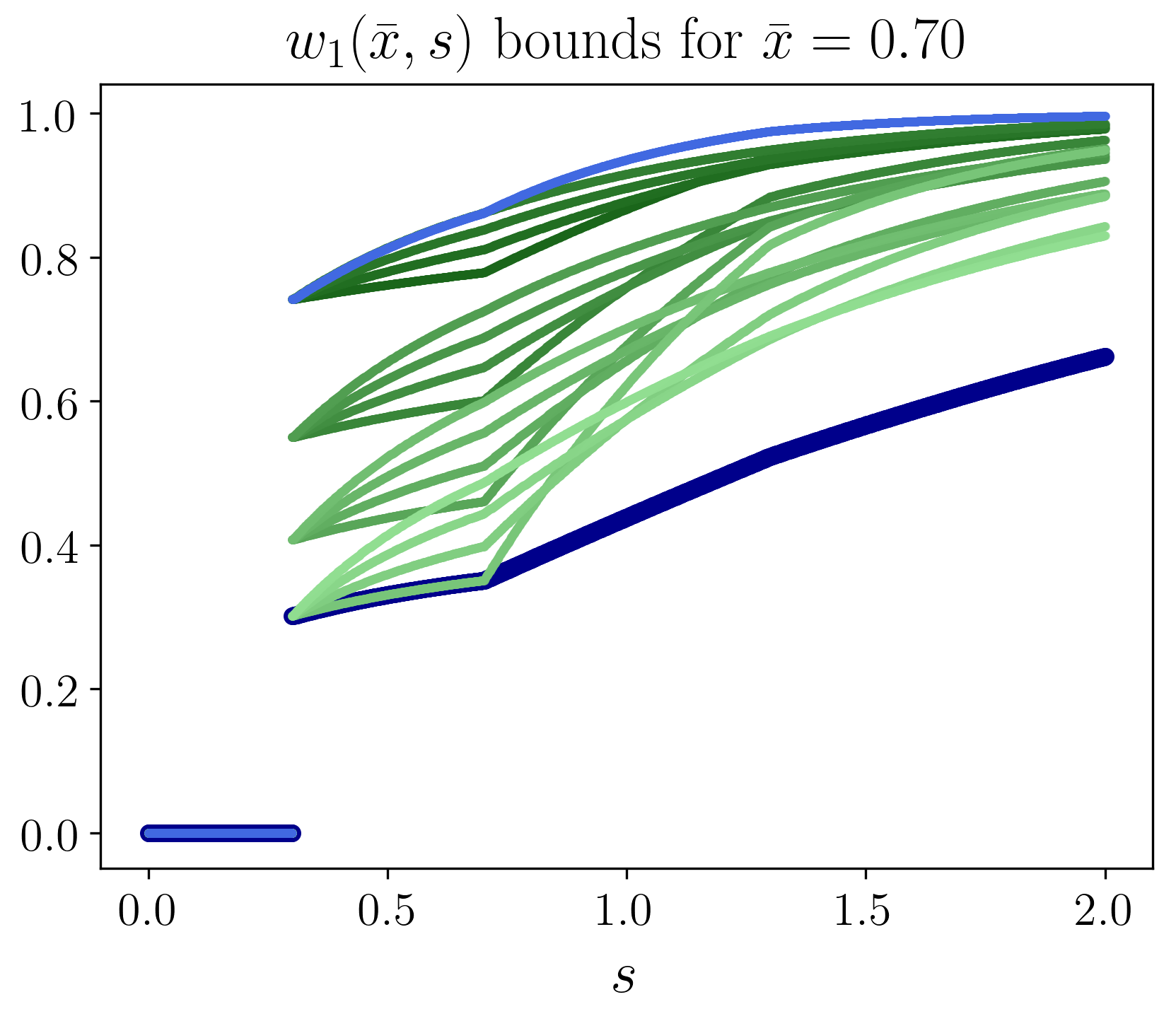}
	\end{subfigure}
	\caption{Example 4: bounds on CDF (starting in Mode 1) for two different initial  positions $\bar{x}=0.3$ and $\bar{x}=0.7.$ 
	See a detailed legend in caption of \cref{fig:bounds_unequal2}.}
	\label{fig:bounds_unequal_cdf2}
\end{figure}

We now generalize Example 1 from \cref{subsec:cdf_experiments} to consider epistemic uncertainty.
Recall that $\domain = [0,1], \, Q = \boundary, \, C \equiv 1,$  and $q \equiv 0$ so that the cumulative cost $\J$ corresponds to the exit time,
with $M=2$ and $f_i = (-1)^{i+1}$.  We will also assume that $\lambda_{ij} \in [1,4]$ for all $i \neq j.$
In \cref{fig:bounds_unequal2} we display our results for a particle that starts moving rightward (in Mode 1). 
The graphs shown in blue are the upper and lower bounds on the 
probability of a timely exit (i.e., before a specific deadline $\bar{s}$) for all initial positions $x$.
The bounds on CDF for two starting positions $\bar{x}$ are shown in \cref{fig:bounds_unequal_cdf2}.
All of these bounds are computed from \cref{eq:w_plus} and \cref{eq:w_minus} for the model of epistemic uncertainty where $\Lambda = (\lambda_{ij})$ is 
allowed to fluctuate within $L.$  Under this model, these bounds are sharp 
since they are computed by finding CDF-maximizing (and minimizing) sequences of \(\Lambda\)'s.

We can also compare the blue bounds to 
the corresponding timely-exit probabilities
for a process containing epistemic uncertainty via fixed and unknown (possibly asymmetric) transition matrix $\Lambda$. 
The green curves shown
in Figures \ref{fig:bounds_unequal2} and \ref{fig:bounds_unequal_cdf2} are computed 
by repeatedly solving \cref{eq:cdf_pde}
for a coarse grid of specific $\Lambda$'s in $L$. 
It should be noted that processes with this type of epistemic uncertainty are a subset of those previously discussed, and so the blue bounds will definitely hold
but will no longer be sharp.
This lack of sharpness is not surprising since
changing the transition rate 
can often
result in a ``better'' (higher or lower -- depending on the bound) probability of timely exit.
However, calculating tighter bounds for a ``fixed-unknown-$\Lambda$'' case is computationally expensive.
By inspection of the experimental data, it is clear that 
such sharp bounds
would have to be 
composed of many individual fixed-\(\Lambda\) CDFs. 

\section{Optimizing the CDF}
\label{sec:optimize}


The PDMPs considered in previous sections were not controllable in any way.  Since the dynamics are deterministic in every mode, each random trajectory was fully described by the initial (state, mode) pair and the discrete time sequence of mode switches.  The goal was to develop efficient methods for approximating the CDF of the cost accumulated up till termination.  We now turn to {\em controlled} PDMPs \cite{davis1984pdmp} 
-- a modeling framework useful in a wide range of applications, including 
production/maintenance planning \cite{boukas1990optimal},
control of manufacturing processes \cite{akella1986optimal, bielecki1988optimality, olsder1980time, sethi2012hierarchical},
multi-generational games \cite{haurie2005multigenerational},
economic growth \& climate change modeling\cite{haurie2006stochastic},
trajectory optimization for emergency vehicles \cite{andrews2013deterministic}, 
preventing the extraction of protected natural resources \cite{CarteeVlad_Poaching}, 
and robotic navigation \cite{ShenVlad, GeeVladimirsky2020}.  

We start with expectation-optimal controls considered in the above references, but then switch to selecting controls to manage the uncertainty in $\J$ and 
provide some notion of {\em robustness}.
Robust controls help practitioners to guard against both modeling errors and  prohibitively bad rare outcomes, which may result from random switches.
It might seem natural to mirror the robust approaches popular in traditional stochastic control, but we find them lacking in the PDMP context. 
$H_{\infty}$ controls are the mainstay of robustness for many processes with continuous perturbations \cite{bacsar2008h}, but they are not easily adaptable for discrete mode-switches.
Another popular idea is to minimize $\E [ \exp( \beta \J ) ]$, with the risk-sensitivity coefficient $\beta > 0$ reflecting our desire to avoid bad outcomes \cite{fleming2006controlled}.
For small $\beta$ values, this is roughly equivalent to minimizing a convex combination of $\E[J]$ and $Var[\J].$ 
While implementable with PDMPs, this method does not provide any guarantees on the likelihood of bad scenarios.
We thus develop a different approach to maximize the probability of not exceeding a specific cumulative cost threshold $\bar{s}$.
In subsection \ref{ss:opt_threshold} we develop PDEs to find such optimal policies for all initial configurations and all threshold values $\bar{s}$ {\em simultaneously}.
The numerical methods and computational examples (for $d=1$ and $d=2$) are covered in subsections \ref{ss:opt_numerics} and \ref{ss:opt_examples} respectively.

\subsection{Controlled PDMPs and expectation-optimal policies}
\label{ss:opt_exp}
To obtain a controlled PDMP, we will assume that both the running cost $C$ and velocity $\Bf$ also depend on additional control parameters, which can be 
changed dynamically while the system travels through $\domain \times \M$. 
We will assume that the set of available control values $A$ is a compact subset of $\R^n.$
Throughout this section, we will slightly overload the notation by using $\ba$ to refer to a generic element of $A$ and 
$\ba(\cdot)$ to refer to a generic {\em feedback-control} policy $\ba :  (\domain \times \M) \rightarrow A,$ which selects a control value based on the current 
system state. 
Once we select any specific $\ba(\cdot),$ 
we can define
\begin{equation}
\label{eq:policy_based}
\Bf(\x, i) = \Bf \left( \x, i, \ba(\x, i) \right)
\qquad \text{ and } \qquad C(\x, i) = C \left( \x, i, \ba(\x, i) \right),
\end{equation}
with equations \eqref{eq:dynamics}-\eqref{eq:J} describing the resulting trajectory and cumulative cost.  The latter will be denoted $\J^{\ba(\cdot)} (\x, i)  = \J^{\ba(\cdot)}_i (\x)$ to highlight the dependence on the chosen control policy.  The corresponding expected cost $u^{\ba(\cdot)}_i = \E\left[\J^{\ba(\cdot)}_i(\x)\right]$ 
and the CDF $w^{\ba(\cdot)}_i(\x,s) = \PP\left[\J^{\ba(\cdot)}_i(\x) \leq s \right]$ can then be found from equations \eqref{eq:expected_pde} and \eqref{eq:cdf_pde} respectively.
However, in controlled PDMPs literature the 
problem is usually first formulated as an optimization over a broader class of {\em piecewise open-loop policies}, 
and the dynamic programming argument is then used to show that an optimal policy can be actually found in feedback form.
We provide an overview of this construction below, but refer to  \cite{davis1984pdmp,verms1985optimal,DavisFarid1999} and related literature for technical details.

In deterministic setting, one considers the set of  measurable {\em open loop control functions} $\A^o = \left\{ \balpha^o : \R \rightarrow A \right\},$ with $\balpha^o(t)$ specifying the control value that will be used at the time $t$.
For PDMPs, a piecewise open-loop policy specifies a new open loop control function to be used after each mode switch. 
Adapting to our setting,  we can define a set of piecewise open-loop policies as 
$\A = \left\{ \balpha: (\R \times \domain \times \M) \rightarrow \A^o \right\}.$  The three inputs to $\balpha$ encode
all information about the last mode switch encountered before the current time $t$: the time and position ($t_{\#} \geq 0$ and $\x_{\#} \in \domain$) where that switch has happened and the resulting mode $i_{\#}.$  If no switch has occurred since we started, we will take $t_{\#} = 0$ and $(\x_{\#}, i_{\#})$ equal to the original (position, mode) pair.  
Assuming $\balpha(t_{\#},\x_{\#}, i_{\#}) = \balpha^o(\cdot) \in \A^o,$  the control value to be used at the time $t \geq t_{\#}$ will be specified by $\balpha^o(t - t_{\#})$ 
until we switch from the
mode $i_{\#}.$
Slightly abusing the notation, we can now replace \eqref{eq:dynamics} and \eqref{eq:J} by
\begin{align}
\label{eq:dynamics_p}
\y'(t)  &=  \; \Bf_{m(t)} \big( \y(t),   \balpha^o(t- t_{\#}) \big),\\
\label{eq:J_p}
\J^{\balpha(\cdot)}_i(\x) & = \; \int_{0}^{T_{\x, i}} C_{m(t)} \Big(\y(t),  \balpha^o(t- t_{\#}) \Big) \, dt \, + \, 
q \Big( \y\left( T_{\x, i} \right)\!, m\left(T_{\x, i} \right)\! \Big),
\end{align}
where $\balpha^o \in \A^o$ is the open loop control function currently in effect at the time $t$ based on the policy $\balpha \in \A$ and 
a sequence of mode switches that have occurred so far.
Recall from section \ref{sec:Intro} 
that $\y(0) = \x, \, m(0) = i$, and the changes in mode $m(t)$ are governed by the matrix of switching rates  $\Lambda.$
In this section, we will further assume that all $\Bf_i$'s and $C_i$'s are Lipschitz-continuous in both arguments.

The usual goal in controlled PDMPs literature is to minimize the expected total cost up to the termination time.  
The value function is thus defined as
\begin{equation}
\label{eq:uhat_def}
\uhat(\x, i) = \uhat_i(\x) = \inf\limits_{\balpha(\cdot) \in \A} \mathbb{E}[\J^{\balpha(\cdot)}_i(\x)],
\end{equation}
The existence of an {\em expectation-optimal} policy $\balpha_*(\cdot) \in \A$ such that $\uhat_i(\x) = u^{\balpha_*(\cdot)}_i(\x)$
is only guaranteed under additional assumptions; e.g., if the set 
\begin{equation*}
{\mathlarger\nu}(\x,i) \; = \; \big\{ \left( r, \Bf_i(\x, \ba) \right) \, \mid \, r \geq C_i(\x, \ba), \,  \ba \in A \big\}
\end{equation*}
is convex for every $\x$ and $i$.  (Alternatively, the existence of optimal policy is also assured if one allows {\em relaxed control functions},
with $\balpha^0$ taking values in the set of probability measures on $A$; see \cite{verms1985optimal,bardi2008optimal}.)
If such an optimal $\balpha_*(\cdot) \in \A$ exists,
it is easy to see that a ``tail'' of $\balpha^o_*$ must be also optimal for every $(\y(t), m(t))$ as long as the process continues.
Otherwise, we could obtain an improvement for the starting configuration $(\y(0),m(0)) = (\x,i)$ by concatenating $\balpha_*$ up to the time $t$ 
with whatever policy is optimal starting from $(\y(t), m(t)).$
A version of this {\em tail-optimality} property holds more generally, even when no expectation-optimal policy exists:
\begin{equation}
\label{eq:tail_opt}
\uhat(\x, i) \;= \; 
\inf\limits_{\balpha(\cdot) \in \A} 
\mathbb{E} \left[
 \int_{0}^{\tau} C_{m(t)} \Big(\y(t), \, \balpha^o(t- t_{\#}) \Big) \, dt \; + \; 
\uhat \Big( \y( \tau ), m( \tau ) \Big)
\right],
\end{equation}
for all $\tau >0$ sufficiently small to guarantee that $\y(t) \in \domain \backslash Q$ for all $t \in [0, \tau], \, \balpha(\cdot) \in \A.$
I.e., we assume that $\tau$ is small enough so that the system cannot reach $Q$ by $t= \tau$ regardless of the sequence of mode switches.
A standard argument 
based on Taylor-expanding \eqref{eq:tail_opt} in $\tau$ 
(e.g., see \cite[\S2]{ShenVlad})
shows that, if $\uhat_i$'s are sufficiently smooth, 
they must satisfy a 
system of Hamilton-Jacobi-Bellman PDEs:
\begin{equation}
\label{eq:uhat_pde}
\min_{\ba \in A} 
\left\{
\nabla \uhat_i(\x) \cdot \Bf_i(\x, \ba) + C_i(\x, \ba)
\right\}
\, + \, 
\sum_{j \not= i} \lambda_{ij}\bigl(\uhat_j(\x) - \uhat_i(\x)\bigr) 
 =  0,
\quad
\x \in \domain \backslash Q, \, i \in \M
\end{equation} 
with boundary conditions $\uhat_i(\x) = q_i(\x)$ for all $\x \in Q$.
In a non-smooth case, these value functions can be still interpreted as the unique viscosity solution \cite{DavisFarid1999}.
The system \eqref{eq:uhat_pde} is a natural non-linear generalization of \eqref{eq:expected_pde} and can be similarly discretized by semi-Lagrangian techniques.  
However, the coupling between different modes makes it difficult to solve the discretized system efficiently even in the case of simple/isotropic cost and dynamics.
A variety of Dijkstra-like non-iterative methods developed for deterministic problems (e.g., \cite{sethian1996fast, sethian2003ordered, chacon2012fast, alton2012ordered, mirebeau2014efficient}) will not be applicable for $\Lambda \neq 0$ and one has to resort to
slower iterative algorithms instead \cite{ShenVlad, GeeVladimirsky2020}.  

Once $\uhat_i$'s are computed, 
an optimal {\em feedback policy} $\ba_*(\x,i)$
can be defined pointwise (for all $\x$ and $i$ simultaneously) by utilizing $\argmin$ values\footnote{
Additional assumptions on $\Bf_i$'s and $C_i$'s can be imposed to ensure that this $\argmin$ is a singleton as long as $\uhat_i$ is differentiable \cite{bardi2008optimal}.
But even with these assumptions, the expectation-optimal policy will still be non-unique at the points where $\nabla \uhat_i$ does not exist, and a tie-breaking procedure 
(e.g., based on a lexicographic ordering) can be employed to avoid the ambiguity.}
 from \eqref{eq:uhat_pde}. 
 The $\ba_*(\cdot)$-determined running cost and dynamics defined by \eqref{eq:policy_based} will be only piecewise Lipschitz in $\x,$
 which is precisely the setting considered in section \ref{sec:compute}.
Finally, we note that
the above can be also viewed as an implicit definition for a piecewise open-loop optimal policy $\balpha_* \in \A$: 
 \begin{align*}
& \balpha_* (t_{\#},\x_{\#}, i_{\#}; t) \; = \;
\ba_*
 \; \in \; 
 \argmin_{\ba \in A} 
\left\{
\nabla \uhat_{i_{\#}} \left(\y(t) \right) \cdot \Bf_{i_{\#}} \left(\y(t), \ba \right) + C_{i_{\#}} \left(\y(t), \ba \right)
\right\}.
 \end{align*}

\subsection{PDEs for threshold-specific optimization}
\label{ss:opt_threshold}
In contrast to the above expectation-centric approach, our goal is to generalize the CDF-computation methods of section \ref{sec:compute} 
by choosing control policies that maximize the probability of desirable outcomes.  Two subtleties associated with this approach are worth pointing out before 
we start deriving the optimality equations.  First, the idea of ``generating the optimal CDF'' is misleading unless we state the goal more carefully.  
Given any fixed initial configuration $(\x,i)$ and two feedback control policies $\ba_1(\cdot)$ and $\ba_2(\cdot),$
it is entirely possible (and actually quite common!) that
$w^{\ba_1(\cdot)}_i(\x,s_1) > w^{\ba_2(\cdot)}_i(\x,s_1)$ while $w^{\ba_1(\cdot)}_i(\x,s_2) < w^{\ba_2(\cdot)}_i(\x,s_2).$
So, which of the resulting CDFs is preferable depends on which threshold is more important: is our priority to minimize the chances of 
the cumulative cost exceeding $s_1$ or $s_2$?  In this {\em threshold-specific optimization} setting, we will say that a policy 
$\ba(\cdot)$ is {\em $s$-optimal} if $w^{\ba(\cdot)}_i(\x,s) \, \geq \, w^{\bb(\cdot)}_i(\x,s)$ 
for all allowable control policies $\bb(\cdot).$ 

The second subtlety is in choosing the set of inputs used to define feedback control policies.  In threshold-specific optimization, 
the optimal actions are no longer fully defined by the current state $(\x,i)$.  In addition, they also depend on the cost incurred so far and the desired threshold for the cumulative cost up to the termination.  To handle this complication, we add an extra dimension to our state space, defining a new 
expanded set of piecewise open-loop control policies 
$\Aext = \left\{ \balpha :  (\R \times \domain \times \M \times \R) \rightarrow \A^o \right\},$ 
and the expanded PDMP dynamics:
\begin{align}
\label{eq:dynamics_expanded}
\y'(t)  &=  
\Bf_{m(t)} \Big(\y(t), \, \balpha^o(t- t_{\#})
\Big),\\
\nonumber
\y(0) &= \x \in \domain,\\
\nonumber
c'(t)  &=  
C_{m(t)} \Big(\y(t), \, \balpha^o(t- t_{\#})
\Big),\\
\nonumber
c(0) &= 0,\\
\nonumber
m(0) &= i \in \M.
\end{align}
Here $c(t)$ represents the total cost incurred so far and $m(t)$ is the current mode, evolving through
a continuous-time Markov process on $\M$.  Similarly, the last argument in $\balpha \in \Aext$ is $c_{\#},$
the total cost accumulated by the time of the last mode switch encountered so far.
Assuming $\balpha(t_{\#},\x_{\#}, i_{\#}, c_{\#}) = \balpha^o(\cdot) \in \A^o,$  the control value to be used at the time $t \geq t_{\#}$ will be specified by $\balpha^o(t - t_{\#})$ as long as we remain in mode $i_{\#}.$
We can now define our new threshold-aware value function:
\begin{equation*}
\what_i(\x,s) \; = \;
\sup\limits_{\balpha(\cdot) \in \Aext}
\PP\left[ \J^{\balpha(\cdot)}_i(\x) \leq s \right].
\end{equation*}
We note that a similar expansion of state space and policy class could also be used when defining $\uhat$, but 
it would not make any difference due to the linear properties of expectations.  In contrast, the $c$-dependence of policies is essential for writing down the tail-optimality property of $\what_i$'s:
\begin{equation}
\label{eq:tail_opt_w}
\what(\x, i, s) \;= \; 
\sup\limits_{\balpha(\cdot) \in\Aext} 
\E \left[
\what  \Big( \y( \tau ),  \, m( \tau ), \, s - c(\tau) \Big)
\right],
\end{equation}
for all $\tau >0$ sufficiently small to guarantee that $\y(t) \in \domain \backslash Q$ for all $t \in [0, \tau], \, \balpha(\cdot) \in \Aext.$
Similarly to the derivation of \eqref{eq:cdf_pde}, a Taylor expansion of \eqref{eq:tail_opt_w} yields a 
system of nonlinear PDEs satisfied by $\what_i$'s:
\begin{align}
\label{eq:control_pde_CDF} 
\nonumber
&\max\limits_{\ba \in A} 
\left\{
\nabla \what_i(\x,s) \cdot \Bf_i(\x,\ba) -  C_i(\x,\ba) \frac{\partial \what_i}{\partial s}(\x,s) 
\right\}
+  \sum\limits_{j\neq i}  \lambda_{ij}\left( \what_j(\x,s) - \what_i(\x,s) \right)
\, = \, 0,\\
&\forall \x \in \domain \backslash Q, \; i \in \M, \; s > 0; 
\end{align} 
with the same initial and boundary conditions previously specified for $w_i$'s in \eqref{eq:cdf_ic} and \eqref{eq:cdf_bc}.
We can also restrict the computational domain for $\what_i$'s (and decrease the numerical diffusion in the discretization) by generalizing equations 
\eqref{eq:cont_best_case}-\eqref{eq:cont_best_case_prob} and defining $\hat{s}^0$ and $\what^0_i$'s.

\subsection{Discretization of PDEs and control synthesis}
\label{ss:opt_numerics}
A semi-Lagrangian discretization of \eqref{eq:control_pde_CDF} can be obtained on a grid similarly to the treatment of an uncontrolled case in section \ref{subsec:cdf_numerics}:
\begin{equation}
\label{eq:What_update}
\What_{i,\bk}^n  \; =  \; 
\max\limits_{\ba \in A}
\bigg\{
\sum_{j =1}^M p_{ij}(\tau) \, \What_{j} \Big( \x_{\bk} + \tau \Bf_i(\x_{\bk}, \ba), \, s_n - \tau C_i(\x_{\bk}, \ba) \Big)
\bigg\}.
\end{equation}
Here $\What_{i,\bk}^n \approx \what_i(\x_{\bk}, s_n)$ is a grid function and $\What_{i}$ is its interpolated version defined on $\domain \times \R$.

Once all $\What_i$'s are computed, they can be used to approximate the optimal control not just on the grid but for all $(\x, i, s)$ by 
choosing control values from the set
\begin{equation}
\label{eq:A_argmax}
\hat{A}(\x,i,s) \, = \, \argmax\limits_{\ba \in A}
\bigg\{
\sum_{j =1}^M p_{ij}(\tau) \, \What_{j} \Big( \x + \tau \Bf_i(\x, \ba), \, s - \tau C_i(\x, \ba) \Big)  
\bigg\}.
\end{equation}
Wherever $\nabla \what_i$ is well-defined, we can also use $A(\x, i, s)$ to denote the $\argmax$ set in equation \eqref{eq:control_pde_CDF},
with  $\hat{A}(\x,i,s)$ interpreted as its grid approximation.
 
We note that this procedure allows synthesizing a policy (approximately) optimal with respect to {\em any} desired threshold value.
To obtain an $\bar{s}$-optimal feedback policy 
$\ba(\cdot)$, 
we would simply need to select $\ba(\x,i,c) \in \hat{A}(\x, i, \bar{s}-c).$
However, such policies will be generically non-unique since $\nabla \what_i(\x,s) = 0$ might hold on a large part of $\domain \times (0, +\infty).$
For example, there is always a ``hopeless region'' 
$H = \left\{ (\x,s) \mid \, \what_i(\x,s) = 0, \forall i \in \M \right\}$ since this equality holds by definition whenever $s < \hat{s}^0(\x)$.  
If $C_i$'s do not depend on $\ba$, then $\nabla \what_i(\x,s) = 0$ implies $A(\x,i,s) = A.$
If we start from $(\x_0, i_0)$  such that $\what_{i_0}(\x_0,\bar{s}) < 1,$  then every policy will have a non-zero probability of exceeding the threshold $\bar{s}.$
Which control values are used on $H$ does not change the probability of ``success'' 
$\left(\J^{\ba(\cdot)}_{i_0}(\x_0) \leq \bar{s} \right),$
but it can significantly impact the overall CDF of  that policy. 
In many problems there is also an ``unconditionally successful''  region
$U = \left\{ (\x,s) \mid \, \what_i(\x,s) = 1, \forall i \in \M \right\}.$
If  $\what_{i_0}(\x_0,\bar{s}) = 1$ then an optimal policy will never exceed the threshold 
$\bar{s}$ regardless of the timing of mode switches. If $(\x_0, \bar{s})$ is in the interior of $U$, then the success is guaranteed 
regardless of  control values chosen until we reached $\partial U$, but these choices will generally affect the CDF.
To resolve these ambiguities, we use a tie-breaking procedure in defining optimal policies: whenever $\hat{A}(\x,i,s)$ is not a singleton, we select 
its element that minimizes the expectation.  
(On $H$ this will coincide with an expectation-optimal policy $\ba_*(\cdot).$ 
But on $U$ this need not be the case 
since our optimal policy is $c$-dependent and we need to account for expected values on $\partial U$.) 

Assuming that $\Vhat_{i,\bk}^n$ is a grid function approximating the expected outcome and  $\Vhat_i$ is its interpolated version, we
can summarize the computational process as follows:
\begin{align}
\label{eq:What_update_1}
\What_{i,\bk}^n  &=  \; \sum_{j =1}^M p_{ij}(\tau) \, \What_{j} \Big( \x_{\bk} + \tau \Bf_i(\x_{\bk}, \hat{\ba}_{i,\bk}^n), \, s_n - \tau C_i(\x_{\bk}, \hat{\ba}_{i,\bk}^n) \Big);
\\
\label{eq:Vhat_update_1}
\Vhat_{i,\bk}^n  &=  \; \tau C_i(\x_{\bk}, \hat{\ba}_{i,\bk}^n) \, + \, 
\sum_{j =1}^M p_{ij}(\tau) \, \Vhat_{j} \Big( \x_{\bk} + \tau \Bf_i(\x_{\bk}, \hat{\ba}_{i,\bk}^n), \, s_n - \tau C_i(\x_{\bk}, \hat{\ba}_{i,\bk}^n) \Big);
\\
\label{eq:ahat_full}
\hat{\ba}_{i,\bk}^n & \in \, \argmin\limits_{\ba \in \hat{A}(\x_{\bk},i,s_n)}
\bigg\{
\tau C_i(\x_{\bk}, \ba) \, + \, 
\sum_{j =1}^M p_{ij}(\tau) \, \Vhat_{j} \Big( \x_{\bk} + \tau \Bf_i(\x_{\bk}, \ba), \, s_n - \tau C_i(\x_{\bk}, \ba) \Big)
\bigg\}. 
\end{align}
The above description removes almost all ambiguity from the synthesis of threshold-optimal policies, 
but the $\argmin$ in \eqref{eq:ahat_full} might still have multiple elements on a set of measure zero in $\domain \times \mS$. 
In such rare cases, additional tie-breaking can be used based on another criterion (e.g., a lexicographic ordering).  

\subsection{Numerical experiments}
\label{ss:opt_examples}
We first illustrate these subtleties of policy synthesis with a simple example on a one-dimensional state space $\domain = [0,1]$ and two modes, each with its own preferred (faster) direction of motion.

\begin{figure}
	\centering
	\includegraphics[width=\linewidth]{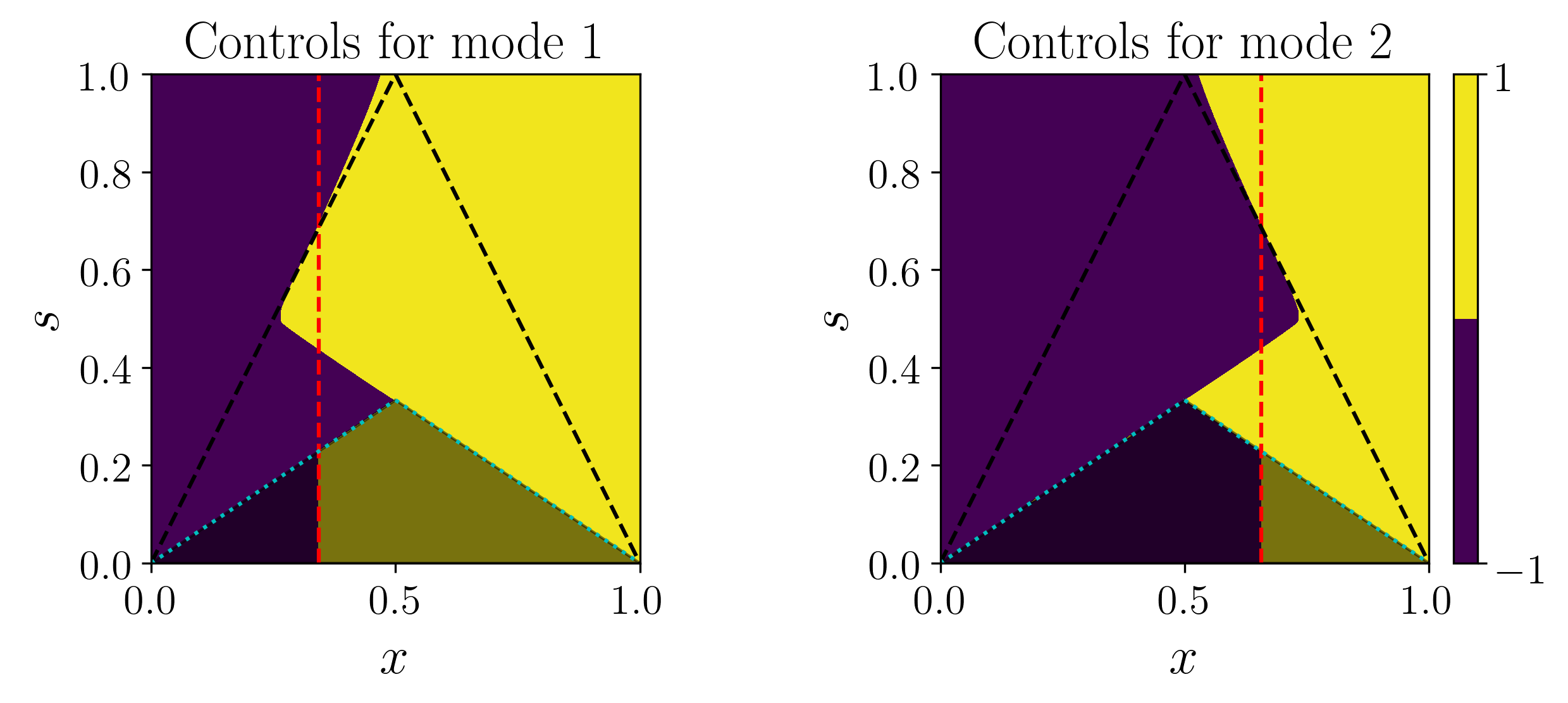}
	\caption{Example 5: a map of threshold-optimal control values with position on the horizontal axis and time remaining until the deadline on the vertical axis. The purple color represents the optimal choice of moving to the left, and the yellow color represents the optimal choice of moving to the right. The shaded area with the cyan border represents the ``hopeless region'' $H,$ where $w_i$'s are uniformly zero and the threshold-specific optimal policies coincide with the expectation-optimal policy.  The ``unconditionally successful'' region $U$ is shown above the black dashed line. 
Under grid refinement, everything in the left part of $U$ becomes purple and everything in the right part of $U$ becomes yellow in both modes. The red dashed vertical lines show the point of direction-switching for the expectation-optimal policy. 
Computed with $\Delta x = 1.25 \cdot 10^{-4}, \, \Delta s = 0.625 \cdot 10^{-4}.$}
\label{fig:controlplot}
\end{figure}

\begin{figure}
	\centering
	\includegraphics[width=\linewidth]{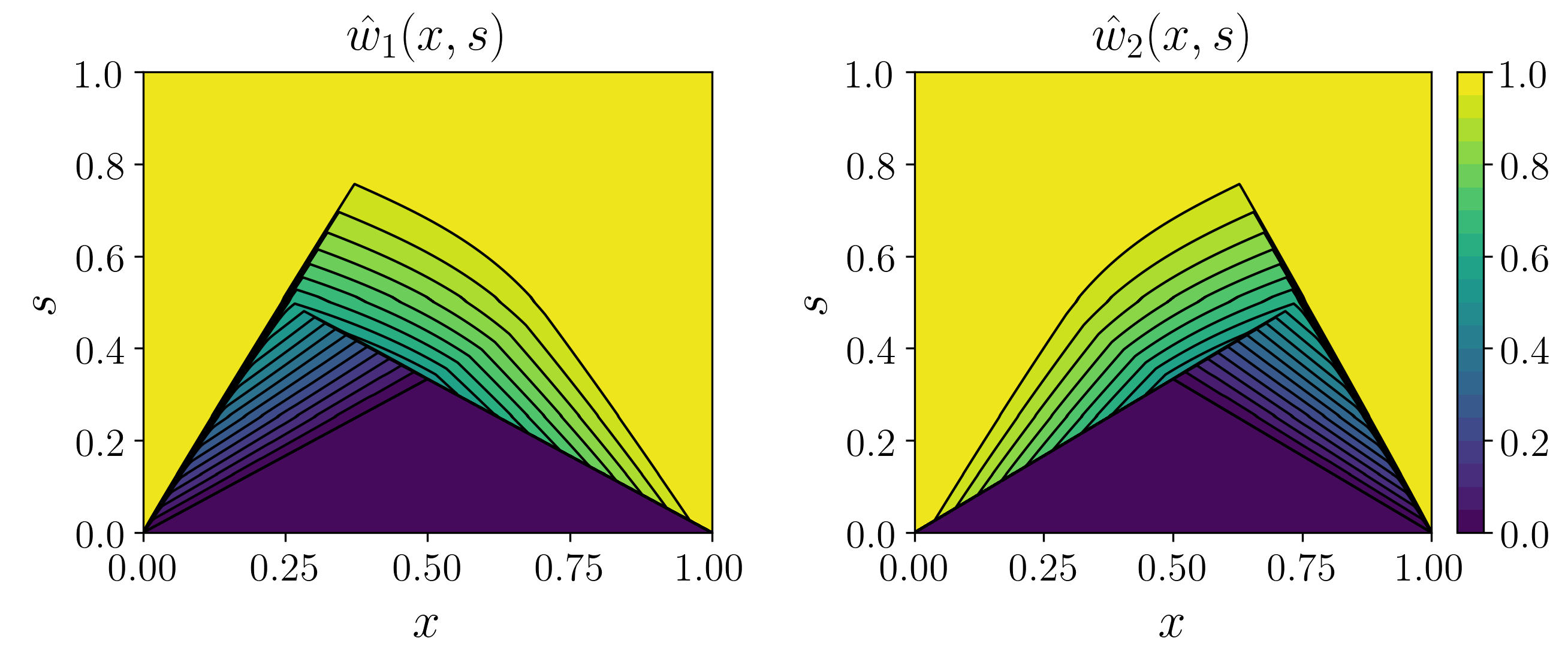}
	\caption{Contour plot of $\what_i(\x,s)$ for Example 5.
	Computed with $\Delta x = 1.25 \cdot 10^{-4}, \, \Delta s = 0.625 \cdot 10^{-4}.$}
	\label{fig:what_plot}
\end{figure}

\begin{figure}
	\centering
	\includegraphics[width= \linewidth]{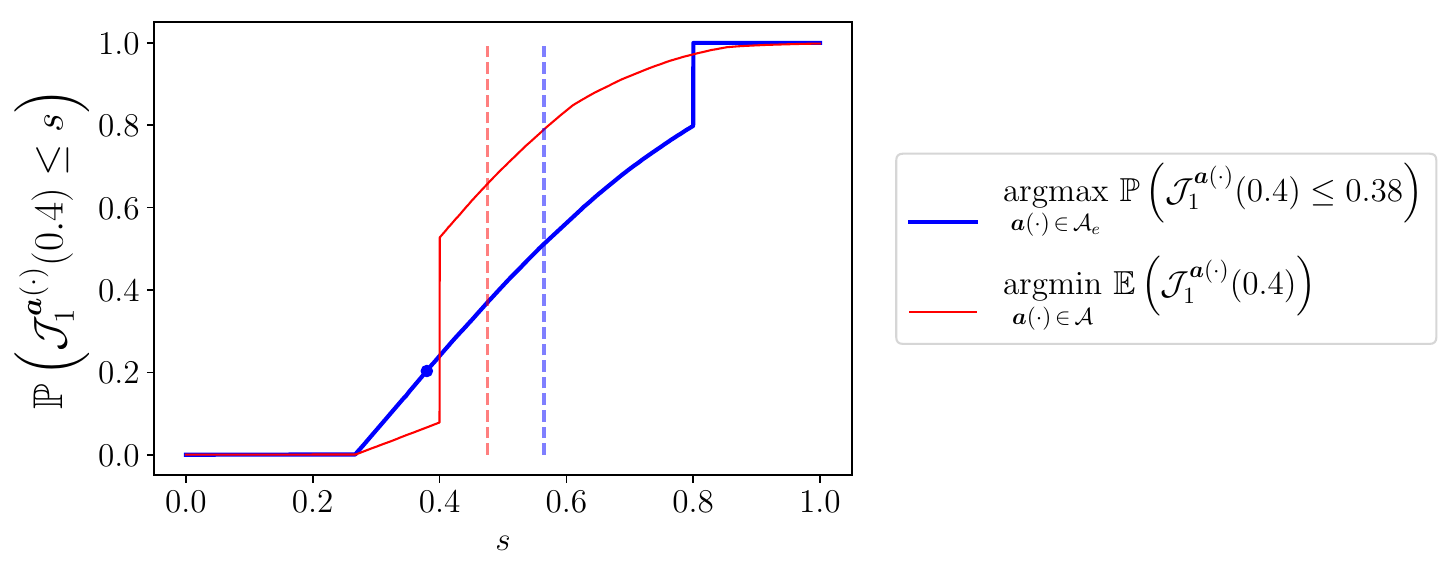}
	\caption{Example 5: CDF of an expectation-optimal policy (in red) and CDF of a threshold-specific optimal policy computed for $\bar{s} = 0.38$ (in blue).
	In both cases, the starting configuration is $(\x_0, i_0) = (0.4, 1).$
	The value of the CDF at the threshold $\bar{s} = 0.38$ is marked by a blue dot.
	The vertical dashed lines indicate the expected value of each policy.}
	\label{fig:cdf_comparison}
\end{figure}

{\bf Example 5}: 
More precisely, the control value $\ba \in A = \{-1, 1\}$ specifies the chosen direction of motion, and the dynamics are $\Bf_i(\x, \ba) = \ba + (-1)^{i-1}\frac{1}{2}$ with $i=1,2.$
In other words, in mode 1 we can move right with speed $3/2$ and left with speed $1/2,$ while in mode 2 it is the opposite. 
We use $q \equiv 0$ on $Q=\boundary$ and $C_1 \equiv C_2 \equiv 1$, ensuring that the cumulative cost $\J$ is just the time to target.
For simplicity, we also use symmetric switching rates $\lambda_{12} = \lambda_{21} = 2.$ 
The resulting optimal policies and the contour plots of $\what_i(\x,s)$ are shown in Figures \ref{fig:controlplot} and \ref{fig:what_plot} respectively.
In Figure \ref{fig:cdf_comparison} we fix a starting configuration and compare the CDFs of two different policies.
The expectation-optimal feedback policy $\ba_*(\cdot)$ 
is obtained by solving \eqref{eq:uhat_pde} and its CDF is then found by solving \eqref{eq:cdf_pde}.
Unfortunately, the same approach is not available for threshold-specific optimal policies: for an $\bar{s}$-optimal feedback policy 
$\ba(\cdot),$ 
there is no reason to expect $\what_i(\x, s) = \PP\left[\J^{\ba(\cdot)}_i(\x) \leq s \right]$ unless $s = \bar{s}$.
Instead, we approximate their CDF using 100,000 Monte-Carlo simulations\footnote{While we do not pursue this alternative here, one could also approximate this CDF by solving the Kolmogorov Forward Equation with initial conditions chosen based on this specific starting configuration $(\x,i)$.}.
Not surprisingly, the threshold-specific policy reduces the probability of missing the deadline $\bar{s}$ but at the expense of increasing the expected time to target.  

Moreover, threshold-specific optimal policies (and their respective CDFs) may also vary significantly depending on the chosen threshold $\bar{s}.$
To illustrate this, we now consider an example on a two-dimensional state space $\domain = [0,1] \times [0,1]$ with four modes, each with its own faster direction of motion. 

\begin{figure}
	\centering
	\includegraphics[width= \linewidth]{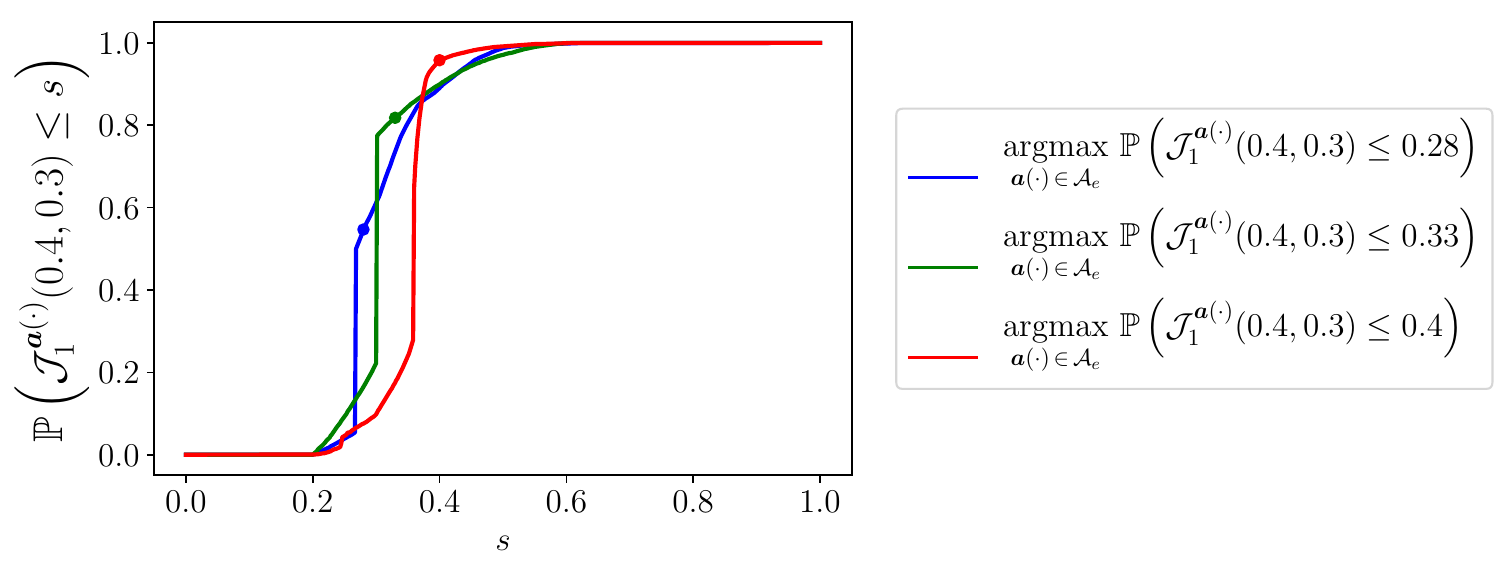}
	\caption{Example 6: CDFs of threshold-specific optimal policies computed for $\bar{s} = 0.28$ (in blue), 
	$\bar{s} = 0.33$ (in green), and $\bar{s} = 0.40$ (in red).
	The value of the CDFs at each threshold are denoted by dots of the corresponding color.
	In all cases, the starting configuration is $(x_0,y_0,  i_0) = (0.4, 0.3, 1).$}
	\label{fig:2d_cdf_comparison}
\end{figure}

\begin{figure}
\begin{tabular}{>{\centering\arraybackslash}m{0.65in}>{\centering\arraybackslash}m{0.85in}>{\centering\arraybackslash}m{0.85in}>{\centering\arraybackslash}m{0.85in}>{\centering\arraybackslash}m{0.85in}}
& $\what_1  \quad \leftarrow$ & $\what_2 \quad \uparrow$ & $\what_3 \quad \rightarrow$ & $\what_4 \quad \downarrow$ \\
s = 0.125 & \includegraphics[height=0.16\textwidth]{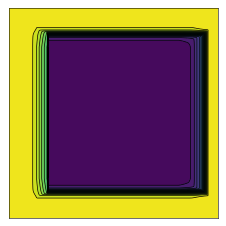} & \includegraphics[height=0.16\textwidth]{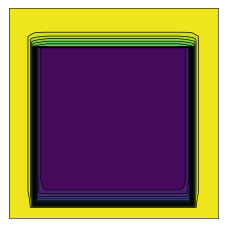} & \includegraphics[height=0.16\textwidth]{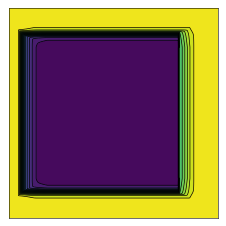} & \includegraphics[height=0.16\textwidth]{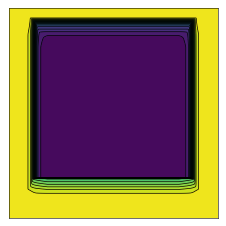} \\
s = 0.250 & \includegraphics[height=0.16\textwidth]{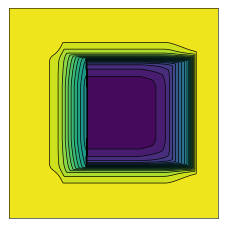} & \includegraphics[height=0.16\textwidth]{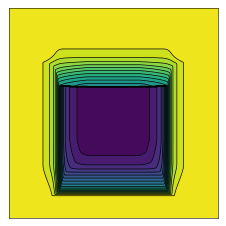} & \includegraphics[height=0.16\textwidth]{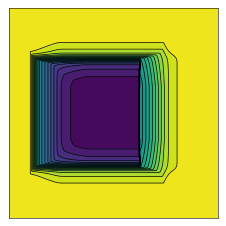} & \includegraphics[height=0.16\textwidth]{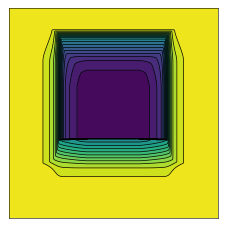} \\
s = 0.375 & \includegraphics[height=0.16\textwidth]{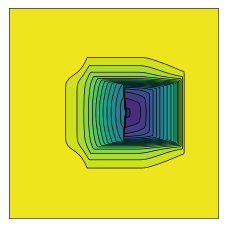} & \includegraphics[height=0.16\textwidth]{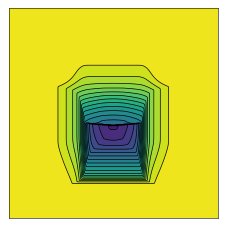} & \includegraphics[height=0.16\textwidth]{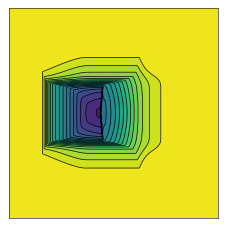} & \includegraphics[height=0.16\textwidth]{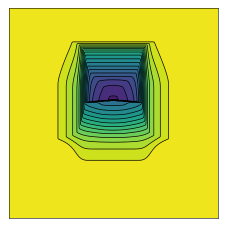} \\
s = 0.500 & \includegraphics[height=0.16\textwidth]{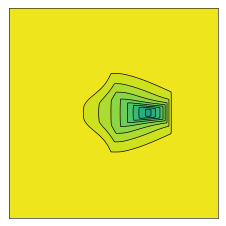} & \includegraphics[height=0.16\textwidth]{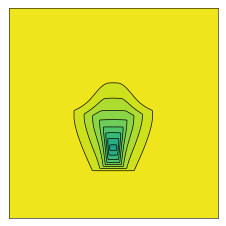} & \includegraphics[height=0.16\textwidth]{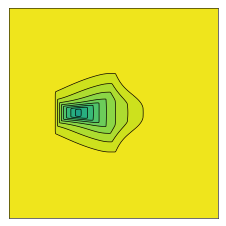} & \includegraphics[height=0.16\textwidth]{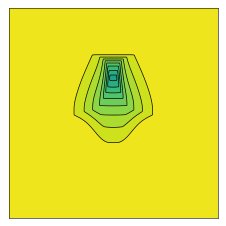}
\end{tabular}
\begin{center}
\hspace{0.6in}\includegraphics[width=0.5\textwidth]{Images/2D/colorbar.pdf}
\end{center}
\caption{Contour plots of $\what_i(\x,s)$ for Example 6. Transition rate between all modes is $\lambda = 1$. Each subplot contains a snapshot of $\what_i(\x,s)$. Each row has a fixed $s$ value, and each column has a fixed mode $i$. Dynamics are given by \cref{eq:2d_dynamics}. Computed on $\domain \times \mS = [0,1]^3$ with $\Delta x = \Delta y = \Delta s = 0.0025$.}
\label{fig:2d_what_plot}
\end{figure}

{\bf Example 6}: 
The control values $\ba$ now reside in $A = 
\{ \ba \in \R^2 \, \mid \, |\ba| = 1 \},$ and the dynamics are given by
\begin{align}
\label{eq:2d_dynamics}
&\Bf_1(\x,\ba) = \ba + \begin{bmatrix} -0.5 \\ 0 \end{bmatrix}, \hfill &\Bf_2(\x,\ba) &= \ba + \begin{bmatrix} 0 \\ 0.5 \end{bmatrix}, \\ 
\nonumber
&\Bf_3(\x,\ba) = \ba + \begin{bmatrix} 0.5 \\ 0 \end{bmatrix}, \hfill &\Bf_4(\x,\ba) &= \ba + \begin{bmatrix} 0 \\ -0.5 \end{bmatrix}.
\end{align}
Again, we use $q \equiv 0$ on $Q=\boundary$ and $C_i \equiv 1$ for all $i$, ensuring that the cumulative cost $\J$ is just the time to $\boundary.$
The switching rates are $\lambda_{ij} = 1$ for all $i \neq j$.
In Figure \ref{fig:2d_cdf_comparison}, we show the CDFs 
(each approximated using 10,000 Monte-Carlo simulations)
for three different threshold-specific optimal policies with the same starting location.
Not surprisingly, each of these policies is strictly better than others with respect to its particular threshold value.
The contour plots of $\what_i(\x,s)$ at various $s$-slices are also shown in Figure \ref{fig:2d_what_plot}.

\section{Conclusion}
\label{sec:Conclude}

The versatility of Piecewise-Deterministic Markov Processes (PDMPs) makes them a useful modeling framework for applications with non-diffusive random perturbations.
In prior literature on PDMPs, the focus has been mostly on the average/expected performance.  Unfortunately, this ignores the practical 
importance of relatively rare yet truly bad outcomes.  The primary goal of our paper is to address this shortcoming and fully characterize the aleatoric uncertainty in a broad class of discrete and continuous PDMPs.  We have accomplished this in section \ref{sec:compute}, approximating the Cumulative Distribution Function (CDF) for their outcomes by solving a system of linear hyperbolic PDEs.  
Although we did not pursue this here, it would be easy to adapt our approach to compute the CDF of {\em hitting times} in discrete time Markov chains. 
In continuous setting, similar ideas could be also extended to stochastic switching in diffusive systems. 
Despite our focus on time-till-exit examples, the presented approach is suitable for a broader class of running costs and PDMP performance measures. 
We illustrate this with a bioeconomic sample problem described in the Appendix.   

For simplicity of exposition, we have assumed the mode-switching rates $\lambda_{ij}$ to be constant, but it should be easy to extend our framework to state-dependent switching rates $\lambda_{ij}(\x)$.  The case of $\lambda_{ij}$'s deterministically evolving in time can be treated similarly by increasing the dimension of our state space.
But random changes in rates present a more serious challenge, which is also related to handling model uncertainties.  The latter is particularly important in PDMPs
since in many practical applications 
these 
rates are not known precisely and are instead estimated based on historical data.  It is thus useful to characterize the {\em range} of possible CDFs -- a task accomplished in section \ref{sec:bounds},
where tight CDF bounds are developed under the assumption that each (state-independent) transition rate $\lambda_{ij}$ has known bounds but does not necessarily remain constant throughout the process.  

Finally, in section \ref{sec:optimize} we have extended our methods to {\em control} the PDMP dynamics,  showing how to maximize the probability of not exceeding a specific cumulative cost threshold $\bar{s}$.  
Our approach 
is also related to the Stochastic On-Time Arrival (SOTA) formulation, developed in discrete setting by transportation engineers to optimize the routing on stochastic networks \cite{Fan2006, SamaranayakeBlandinBayen2012, ErmonGSV12}.  In the context of SOTA, there is only one ``mode,'' but the running cost is random.   While we do not pursue it here, our method can be similarly adapted to optimize the CDF for a subclass of Markov Decision Processes with deterministic running cost and random successor nodes.

Several generalizations of the described methods will broaden their appeal to practitioners.  
First, all PDMPs considered here were {\em exit-time} problems, with the process terminating as soon as the system enters a specific subset $Q$ of the state space $\domain$.  It will be easy to extend our approach to finite horizon problems, but the extensions to infinite horizon (with time discounting of running cost) or ergodic (time-averaged cumulative cost) problems will be more challenging. 

Second, the classical controlled PDMP models in \cite{davis1984pdmp, verms1985optimal} were more general than the setting presented here: instead of our ``mode switching'' they considered ODE trajectories punctuated by random jumps in state space, with both the rate of jumping and the distribution over the set of post-jump positions generally dependent on the pre-jump state and the chosen control value.  It would be clearly useful to extend our methods to this broader setting.  Our section \ref{sec:bounds} can be viewed as a small step in this direction, since we are essentially {\em controlling} mode-transition rates to either maximize or minimize the CDF.

Third, there are many potential ways to improve the accuracy and computational efficiency.
Our approach relies on solving 
systems of hyperbolic PDEs, whose solutions are typically piecewise continuous. 
While the described implementation is based on a first-order accurate semi-Lagrangian discretization, it would be useful to replace these with higher-order accurate schemes \cite{falcone2013semi}.  Our preliminary experimental results based on ENO/WENO \cite{shu1998essentially} spatial discretization in one dimension seem promising, but we have decided to omit them here due to length constraints.    We have also developed a technique restricting the computational domain by pre-computing the minimal attainable cumulative cost.  In controlled PDMPs with an ``unconditionally successful'' region, further domain restriction techniques might be used to maximize the probability of desirable outcomes while also imposing a hard constraint on the worst-case performance.  This would mirror the approach previously developed for routing on stochastic networks \cite{ErmonGSV12}.


For controlled processes, another interesting challenge is to carefully evaluate all trade-offs between conflicting objectives.  This is usually done by computing
non-dominated (or {\em Pareto-optimal} controls), for which any improvement in one of the objectives must come at the cost of decreased performance based on some other objective(s).  With PDMPs, the natural objectives would include traditional minimization of the expected cumulative cost and maximizing the probability of not exceeding a threshold (possibly for several different threshold values).
In the fully deterministic case, several methods for multiobjective 
optimal control have been developed in the last ten years \cite{KumarVlad, guigue2014approximation, DesillesZidani}.  
It will be useful but more challenging to extend these to the piecewise-deterministic setting.

It would be also very interesting to explore additional notions of robustness for PDMPs.
Our approach can be viewed as a dual of optimizing the Value-at-Risk (VaR), 
in which the goal is to minimize a specific percentile of the random outcome.  
We minimize the probability of exceeding a specific threshold, but similarly to VaR, we provide no guarantees on how bad the outcomes can be once that threshold is exceeded.  The Conditional Value-at-Risk (CVaR) is an extended risk measure which addresses this limitation.  A method based on CVaR optimization has been developed for Markov Decision Processes in \cite{chow2015risk}.  It would be useful to extend it to PDMPs and compare with the threshold-optimal policies described here.


\section*{Acknowledgements}
The authors would like to thank Tristan Reynoso and Shriya Nagpal for their help in the initial stages of this project during the summer REU-2018 program at Cornell University.
The authors are also grateful to anonymous reviewers
whose suggestions have helped us improve this paper. 

\section{Appendix: a fish harvesting example}
\label{sec:harvesting}

To show that our general methods are broadly useful beyond the set of illustrative first-exit-time problems considered above,
we include an example based on a PDMP with non-constant dynamics and non-constant running cost in each mode.
We quantify 
the uncertainty in harvesting fish population (whose changing level is encoded by $y(t)$) in the environment with randomly switching carrying capacity $K$. 
As a motivation for such switching, we note that the fish population in the tropical Pacific depends on upwelling of nutrients due to the common easterly winds.
In El Niño years, these winds weaken, temporarily reducing both the upwelling and the carrying capacity.

\begin{figure}
\centering
	 \includegraphics[width = 0.625\textwidth]{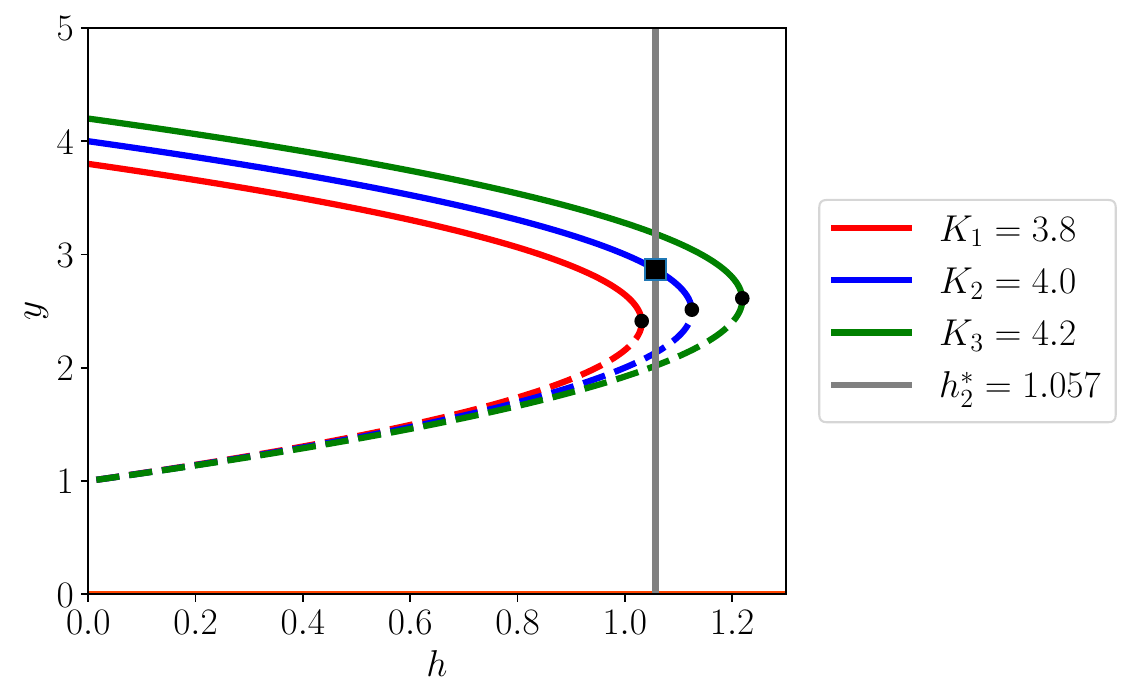}
\caption{Bifurcation diagrams corresponding to the dynamics in each mode.
The orange $y=0$ line is a stable equilibrium for each value of $K$.  The other two equilibria exist for a range of $h$ values and are shown in red, green, and blue
for $K_1, \, K_2,$ and $K_3$ respectively.  The stable $y_+(h)$ is shown by solid lines while the unstable $y_-(h)$ is shown by dashed lines.
The saddle node bifurcations for each $K_i$ are indicated by black dots, while the $K_2$-deterministically-optimal rate $h_2^*$ is shown by the gray line.}
\label{fig:bifurcation}
\end{figure}

The usual logistic population growth model $y' = r(1 - \frac{y}{K}) y$ assumes that the per capita growth rate decreases linearly with the current population size, starting from the rate $r$ when $y=0$
and decreasing to zero if $y$ reaches the carrying capacity $K$.  
This logic reflects the ideas of increasing competition for limited resources when the population grows.  But at low population sizes, other considerations might be more important -- having more individuals might make it easier to find partners for mating, cooperate in finding food, or fend off predators.  This ``Allee effect'' \cite{allee1932} is reflected by having per capita growth rate that first increases (until some threshold value $y=A < K$) and only then decreases (until it reaches $0$ at $y=K$).  Perhaps the simplest model that captures this and includes harvesting is 
\begin{equation}
\label{eq:Allee}
\frac{dy}{dt} \; = \; r \, y \, \left( \frac{y}{A} - 1 \right)  \, \left( 1 - \frac{y}{K} \right) \, - \, h y
\end{equation}
Here, $h \geq 0$ is the effective fishing efforts coefficient, which we will assume to be constant.
Note that $y=0$ is a stable equilibrium for all $h \geq 0$, including the no-fishing case $h=0$.
(This is because we are modeling a {\em strong} Allee effect and $y' <0$ for all $y \in (0, A).$)
But for sufficiently small $h$, the system has two more equilibria: an unstable one at $y_{-}(h)$ and
a stable one at a higher $y_{+}(h)$: 
\begin{equation}
y_{\pm}(h) \; = \; \frac{K+A}{2} \, \pm \, \frac{\sqrt{(K-A)^2 - h\frac{4K}{r}}}{2},
\qquad \qquad
\text{ for } \; h \leq \frac{r(K-A)^2}{4K} = h^{\#}.
\end{equation}
In this regime, the asymptotic behavior depends on the initial conditions:
$\lim\limits_{t \to +\infty} y(t) = y_{+}(h)$ if $y(0) > y_{-}(h)$ and 
$\lim\limits_{t \to +\infty} y(t) = 0$ if $y(0) < y_{-}(h)$.
As shown in Figure \ref{fig:bifurcation},
this bi-stability disappears in a saddle-node bifurcation at  
$h= h^{\#},$
marking the onset of population collapse.
However, we make a distinction between two stages of collapse:
for all $h > h^{\#}$ the collapse of fish population is {\em imminent} since $\lim\limits_{t \to +\infty} y(t) = 0$ regardless of $y(0).$
It can be still reversed by reducing $h$ sufficiently quickly, but becomes {\em irreversible} as soon as $y(t) < A.$ 
We will view this irreversible collapse as a terminal event, motivating our choice of the exit set $Q= \{A\}.$

Another value of obvious relevance is
the optimal level of fishing efforts $h^*$ that maximizes the sustainable yield $h  y_{+}(h)$ over all $h \in [0, h^{\#}].$
A straightforward calculation shows that 
\begin{equation}
h^* \; = \; 
\frac{r}{9K} 
\left[
K^2 + A^2 - 4 K A + (K+A)  \sqrt{ (K-A)^2 + K A }
\right].  
\end{equation}

\begin{figure}
\centering
\begin{subfigure}{0.32\textwidth}
	 \includegraphics[width = \textwidth]{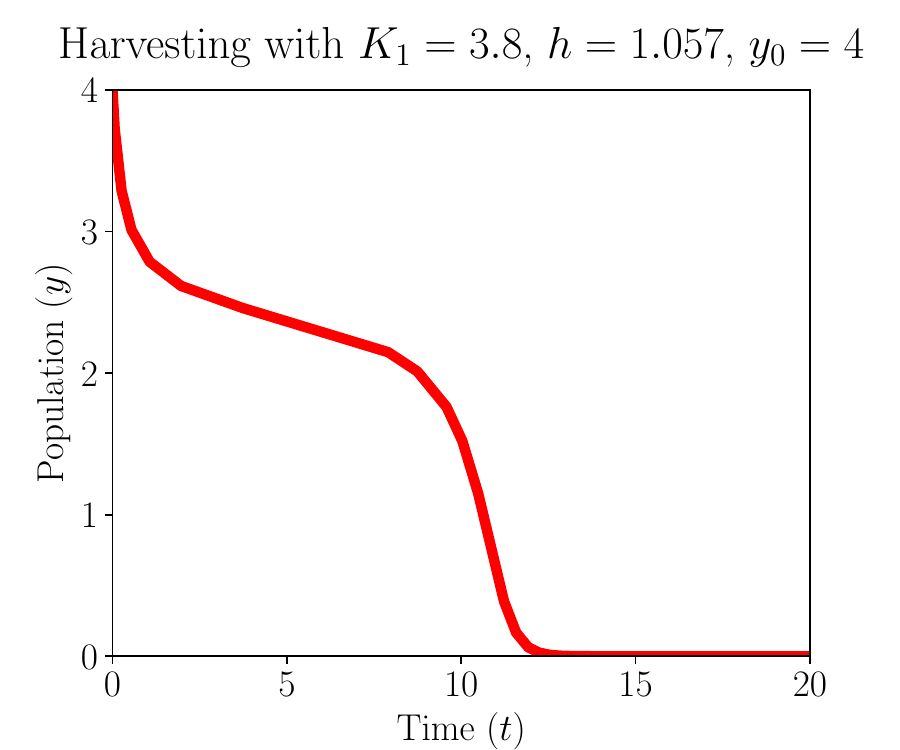}
	 \caption{Mode 1 ($K = 3.8$)}
\end{subfigure}
\begin{subfigure}{0.32\textwidth}
	 \includegraphics[width = \textwidth]{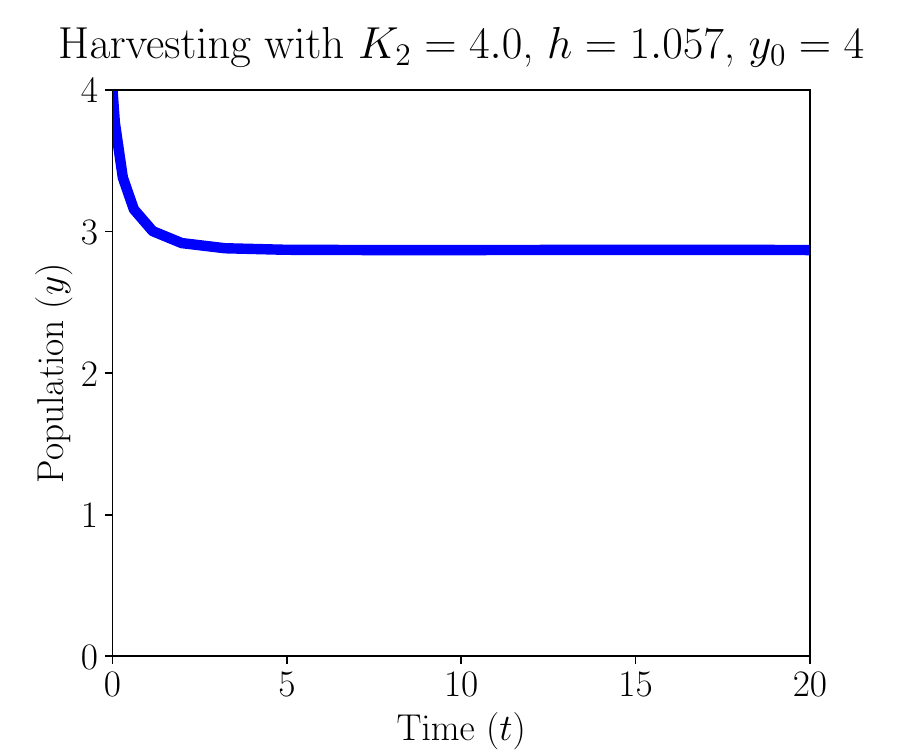}
	 \caption{Mode 2 ($K = 4.0$)}
\end{subfigure}
\begin{subfigure}{0.32\textwidth}
	 \includegraphics[width = \textwidth]{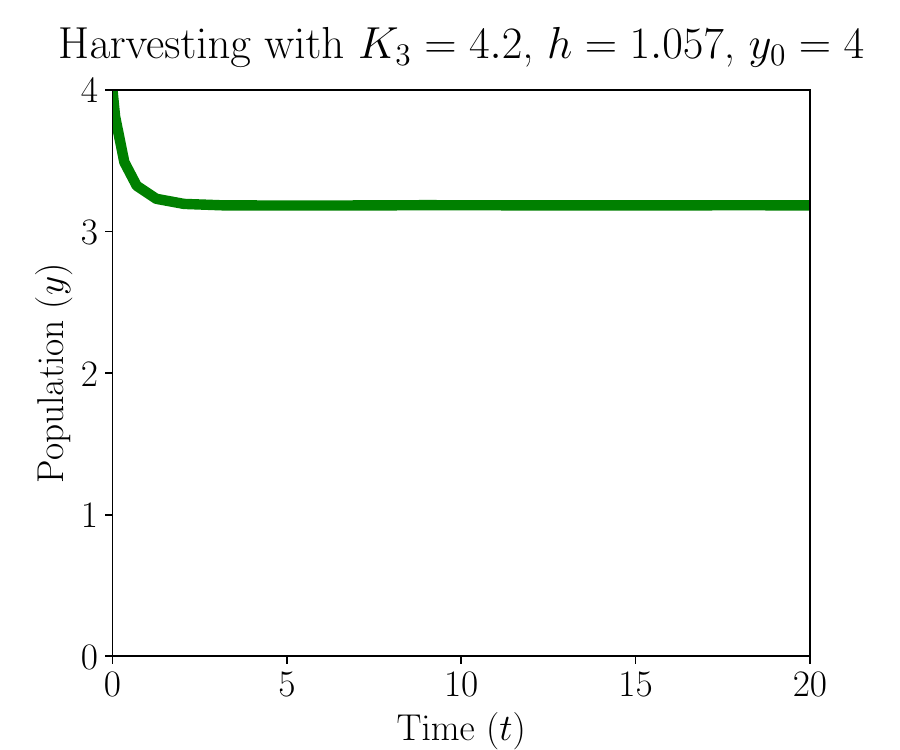}
	 \caption{Mode 3 ($K = 4.2$)}
\end{subfigure}
\caption{Deterministic dynamics for modes 1, 2 and 3 with harvesting rate $h_2^* \approx 1.057$}
\label{fig:determ_harvesting}
\end{figure}

Until now, we have treated all other parameters as fixed and considered the changes to asymptotic behavior as a function of the chosen $h$.
We now turn to a PDMP model with 3 modes, each with its own carrying capacity ($K_1 = 3.8, \, K_2 = 4,$ and $K_3 = 4.2$) and with the other two parameters held constant ($r = 2$ and $A=1$).
For notational convenience, we will use $h^{\#}_i$ and $h^*_i$ to refer to the corresponding maximal sustainable and the yield-optimal fishing rates for each value $K_i$ in the deterministic setting (i.e., if you start in Mode $i$ and there are no mode switches).
But in our computational experiments, we will assume the switching rates $\lambda_{12} = \lambda_{32} = 0.1, \, \lambda_{21} = \lambda_{23} = 0.05, \,$
and $\lambda_{13} = \lambda_{31} = 0.$  As a result, the system spends on average $50\%$ of time in Mode 2 and
$25\%$ of time in each of the Modes 1 and 3.
So, it is natural to view $K_2$ as ``usual'' (or at least as ``average'') and it might be tempting to select its optimal fishing rate $h^*_2 \approx 1.057$.
But if we stay in Mode 1 for a sufficiently long time, this rate will lead to a population collapse since $h^*_2 > h^{\#}_1 \approx 1.032$; see Figures \ref{fig:bifurcation} and 
\ref{fig:determ_harvesting}(a).
In fact, if we stick to the same harvesting rate $h^*_2$ in all modes, this collapse eventually happens with probability one 
as long as $\lambda_{21}, \lambda_{32} > 0.$

One natural question is to quantify the uncertainty in the time until this collapse becomes irreversible; i.e., 
a random time $T_x $ until $y(t) = A = 1$ for a trajectory starting from $y(0) = x$ in Mode $i.$  This could be interpreted as another ``first-exit time problem'', similar to those
considered in \cref{subsec:cdf_experiments}, but with $x$-dependent dynamics in each mode.  
Instead, we have chosen to focus on the CDF for {\em the total amount harvested} before the collapse becomes irreversible:
\begin{equation}
\J_i(x) \; = \; \int_0^{T_{x}} h^*_2 \, y(t) \, dt, 
\qquad \qquad
w_i(x,s) \; = \; \PP(\J_i(\x) \leq s).
\end{equation}
which requires using a non-constant running ``profit''.
To map this back to the notation of section \ref{subsec:continuous}, we will take 
\begin{equation}
f_i(x) =  r \, x \, \left( \frac{x}{A} - 1 \right)  \, \left( 1 - \frac{x}{K_i} \right) \, - \, h^*_2 x
\qquad \text{ and } \qquad C_i(x) = h^*_2 x.
\end{equation}
Our exit set is $Q = \{1\};$  so, we set all $q_i(1)=0$ and note that no boundary condition is needed at the other endpoint since $f_i(4) < 0$ for all $i.$

We solve the three coupled PDEs \eqref{eq:cdf_pde} for $\x \in \domain = [A, K_2] = [1, 4]$ and $s \in [0, 200].$ 
The approximate solution is computed through a semi-Lagrangian scheme \cref{eq:cdf_update}
on a 101 by 120,001 grid, corresponding to $\Delta x = 0.03$ and $\Delta s  =  1/600.$
In this example, $\max_{(x,i)} \lvert f_i(x) \rvert \approx 5.4912$ and $\min_{(x,i)} C_i(x) \approx 1.057.$
So, we use a pseudo-timestep $\tau = 1/600$ to satisfy the inequalities \cref{eq:tau_ineq1} and \cref{eq:tau_ineq2}. 

An obvious lower bound for $\J_i(x)$ is $(x-A)$, but this does not include all the fish born and harvested before $T_{x}.$
The sharp lower bound $s^0(x)$ can be computed by noting that the fastest collapse happens if we stay in Mode 1 throughout.
For the initial condition $x=4$ depicted in Figure \ref{fig:determ_harvesting}(a), this quantity is approximately $s^0(4) \approx 14.87$. 
If starting in Mode 1, the probability of such an outcome is $w_1^0(x)  = \exp{ \big( -\lambda_{12} \mathcal{T} \big)},$ where $\mathcal{T}$ is the time taken by this ``deterministically fastest'' collapse.
If starting in any other mode, this outcome would require an instantaneous transition to Mode 1; so, $w_2^0(x) = w_3^0(x) = 0.$

\begin{figure}
\centering
\begin{subfigure}{0.35\textwidth}
        \includegraphics[width = \textwidth]{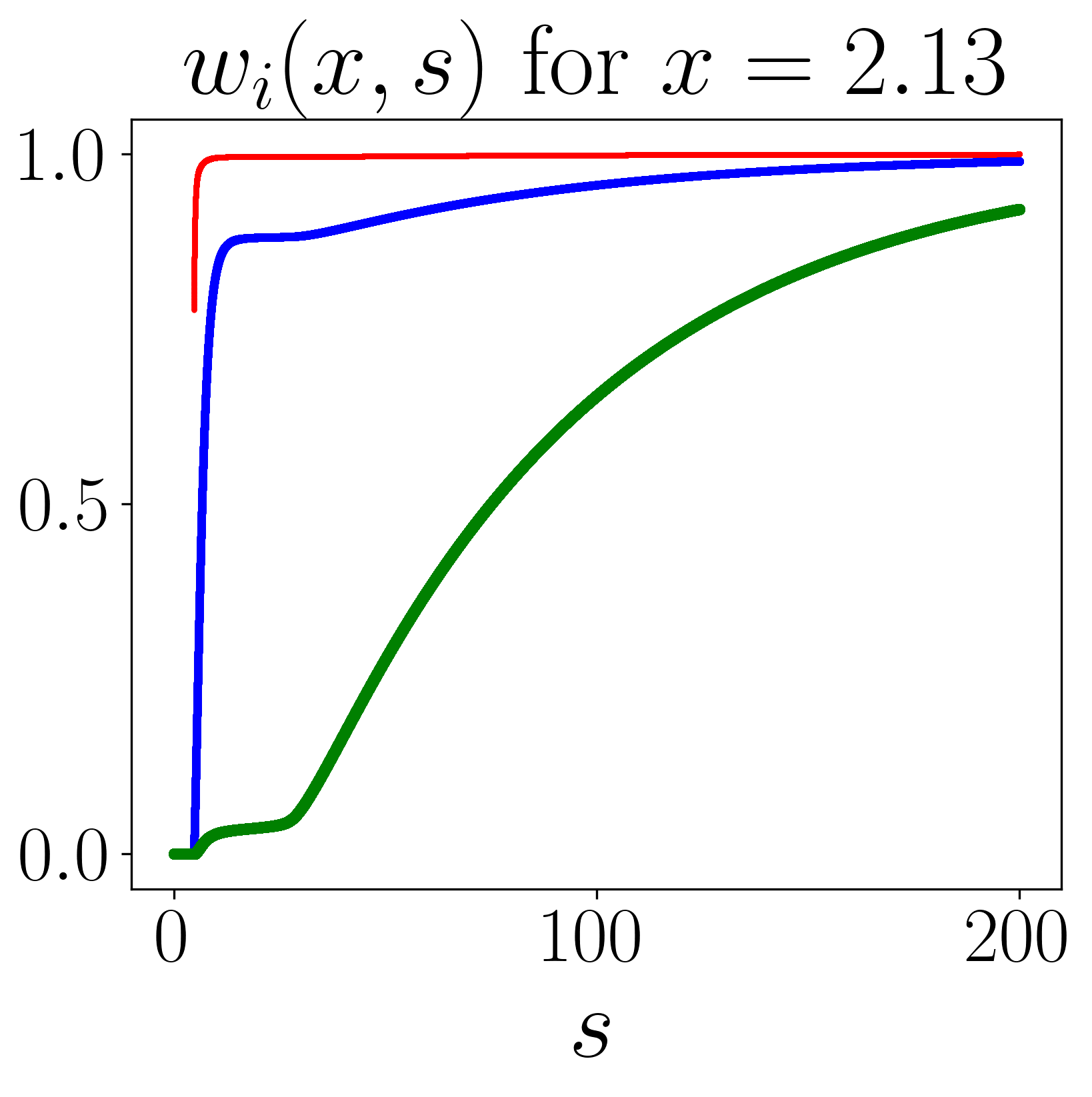}
\end{subfigure}
\hspace*{1.5cm}
\begin{subfigure}{0.35\textwidth}
        \includegraphics[width = \textwidth]{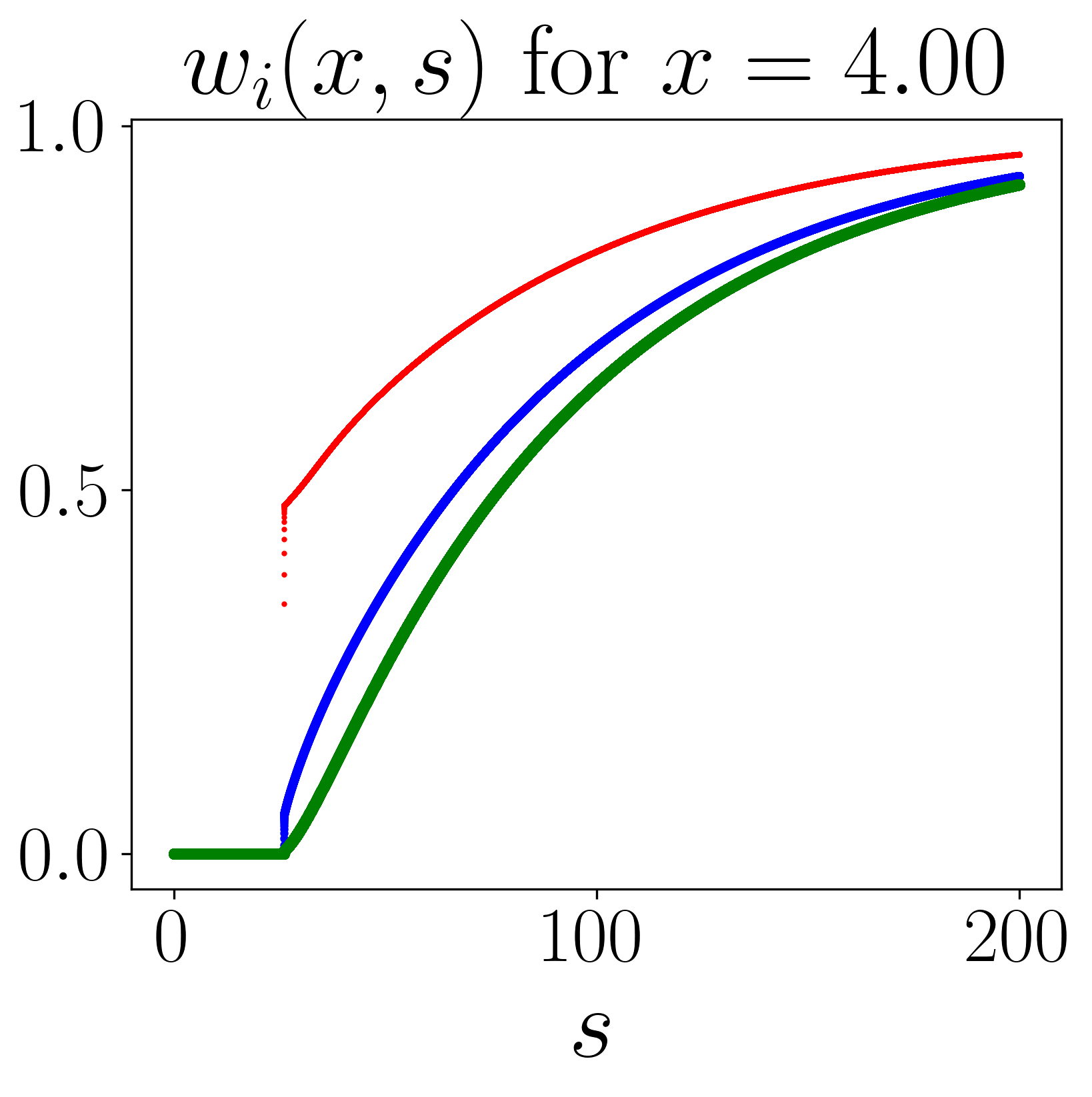}
\end{subfigure}
\caption{CDFs 
computed up to $s=200$ for two starting population levels: $x=2.13$ (left) and $x=4$ (right). Red, blue, and green are starting in Modes 1, 2, and 3 respectively. Plots represent probability of irreversible population collapse at or before a harvest of size \(s\). 
Note that $y_{-}(h^*_2) \approx 2.1312$; so, with $x=2.13$ a quick collapse can only be prevented by an early switch to Mode 3. 
If starting at $x=4$ in Mode 2 (blue curve in the right subfigure) the 25th, 50th, and 75th percentiles are $s \approx 39.09$, $s \approx 66.10$, and $s \approx 112.95$ respectively.}
\label{fig:fishing_CDFs}
\end{figure}

\begin{figure}
\centering
\begin{subfigure}{0.32\textwidth}
        \includegraphics[width = \textwidth]{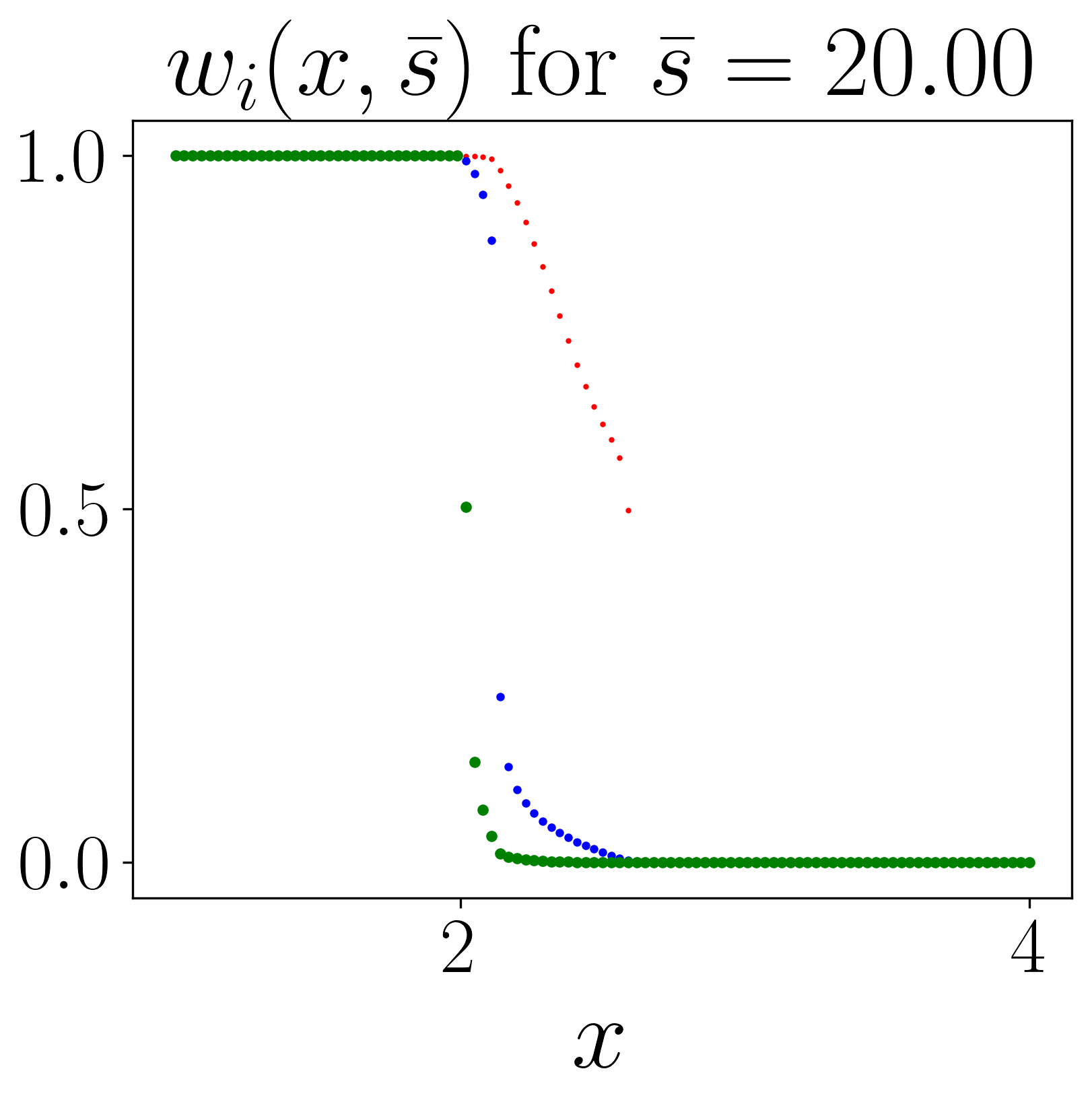}
\end{subfigure}
\begin{subfigure}{0.32\textwidth}
        \includegraphics[width = \textwidth]{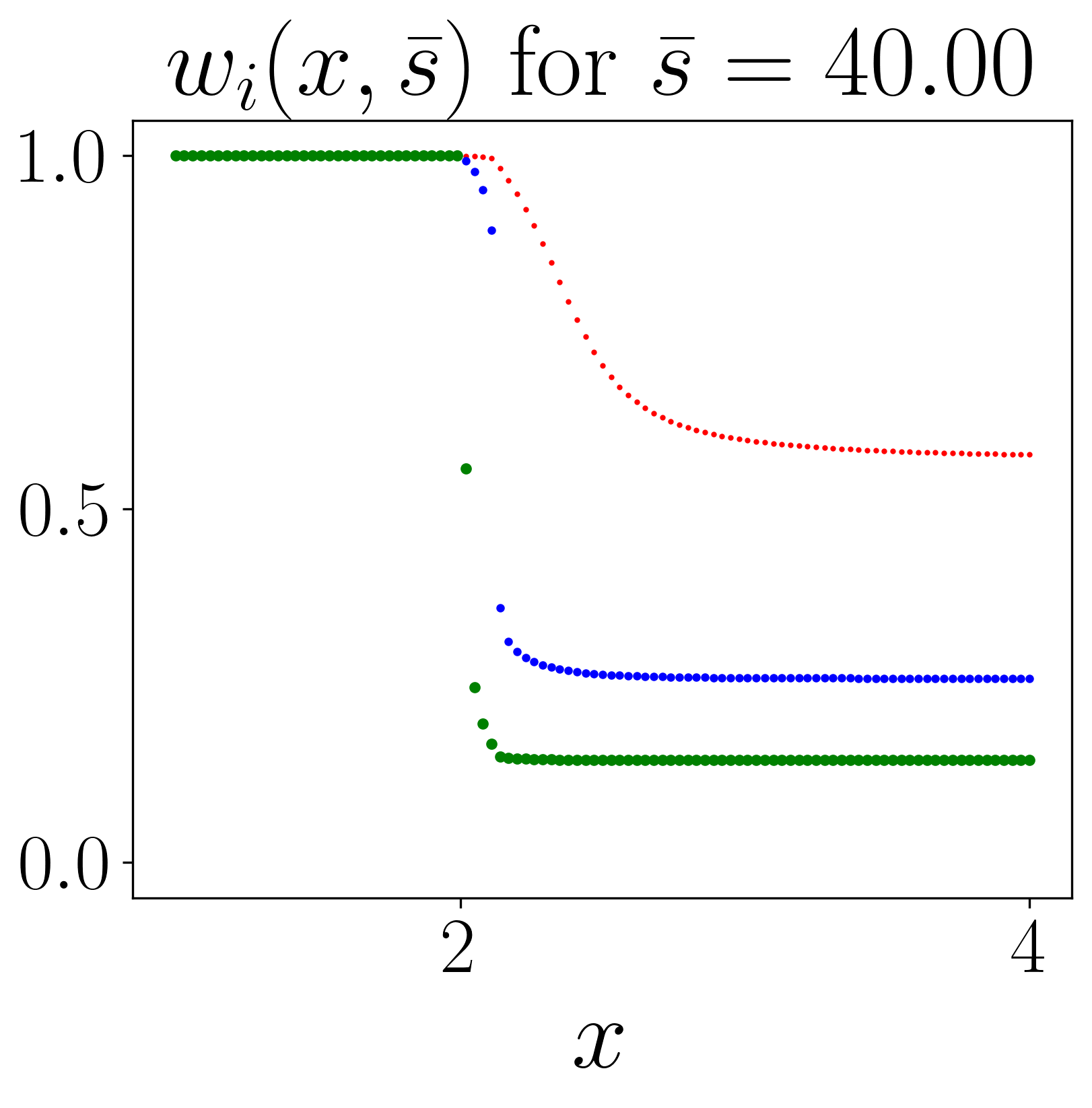}
\end{subfigure}
\begin{subfigure}{0.32\textwidth}
        \includegraphics[width = \textwidth]{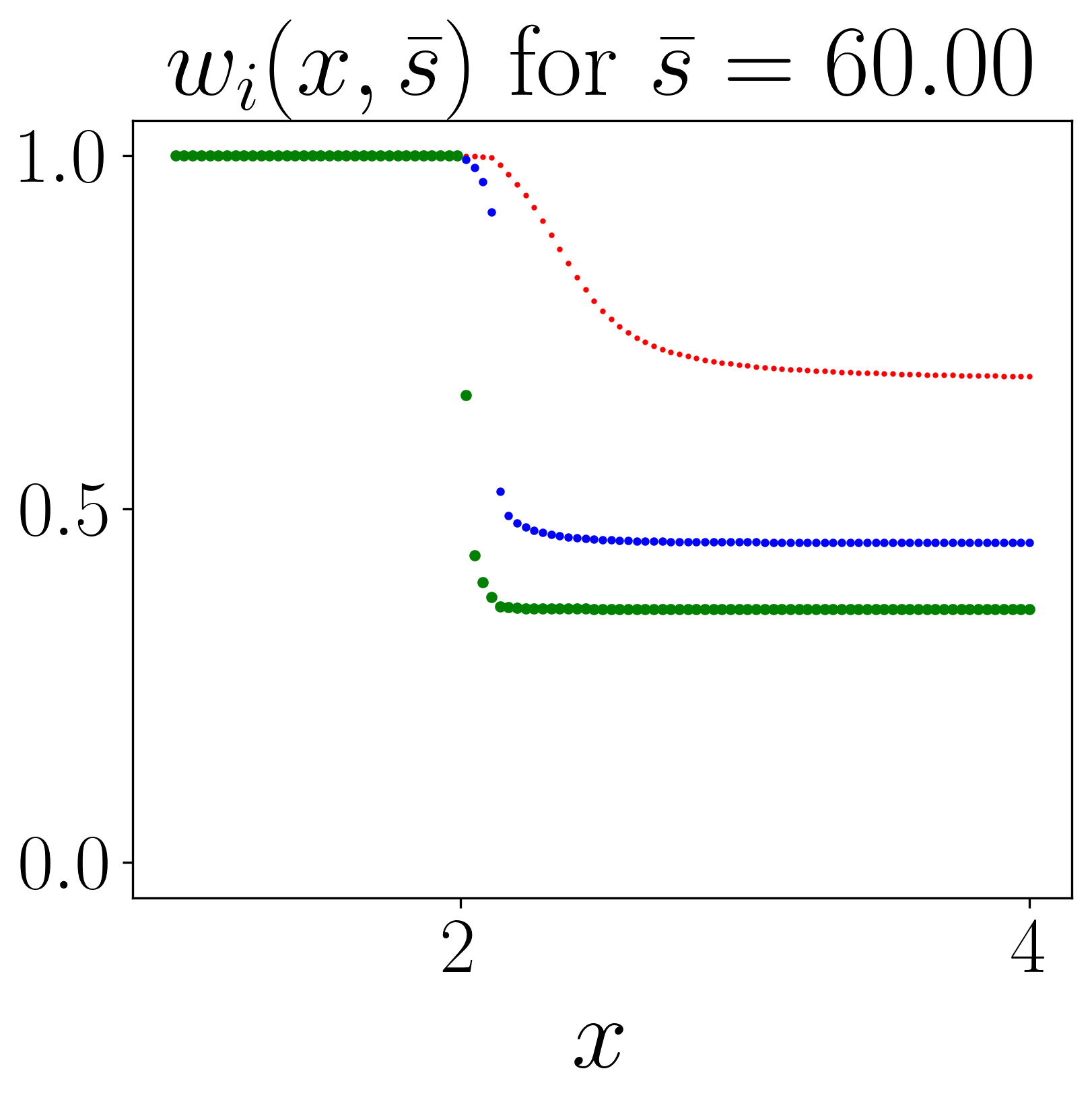}
\end{subfigure}
\caption{Probability of irreversible population collapse before a harvest of size \(\bar{s}\) computed for all starting populations \(x\) and three specific  \(\bar{s}\) values. 
Red, blue, and green represent starting in Modes 1, 2, and 3 respectively.} 
\label{fig:fishing_ws_fixed_sbar}
\end{figure}

Figure \ref{fig:fishing_CDFs} presents the corresponding CDFs $w_i(\bar{x}, s)$ for the initial populations $\bar{x}=2.2$ and $\bar{x}=4$
while the graphs of $w_i(x, \sbar)$ for 3 different values of $\sbar$ are shown in Figure \ref{fig:fishing_ws_fixed_sbar}.
We note that in a deterministic scenario of $K=K_2$, the sustainable equilibrium would be
$y_{+}( h^*_2) \approx 2.8685$
and the sustained optimal yield
(i.e., the amount perpetually harvested per unit time) would be
$
R \; = \; h^*_2  \, \, \, y_{+} (h^*_2) \approx 3.0323.
$
This provides a natural yardstick for thinking about the argument $s$ used in our CDFs.  E.g., based on Figure \ref{fig:fishing_CDFs}, if we start with $y=K_2$ in mode 2 under the random switches defined by $\lambda_{ij}$'s, there is an approximately 50\% chance of harvesting at least $21.8 \times R \approx 66.1$ 
before the collapse becomes irreversible.


The above story is based on the assumptions that mode transitions are not observed and $h$ is chosen once and for all.  In reality, declining catch would provide an advance warning 
that the population starts to collapse and one could reduce $h$ adaptively. 
Selecting $h \leq h_i^{\#}$ (or $h \leq \min_i h_i^{\#},$ if mode switches are not directly observable) would make harvesting indefinitely sustainable.   
One could also use the methods of \cref{sec:optimize} to find the CDF-optimizing harvesting rates in feedback form.
So, our model described above is a vast simplification, but it already illustrates the promise of presented techniques for quantifying uncertainty in bioeconomic applications.

\bibliographystyle{siam}
\bibliography{Ref}

\end{document}